\documentclass[a4paper]{amsart}
\usepackage[utf8]{inputenc}
\usepackage{amssymb,amsmath,enumerate,amsthm,url}
\usepackage[all,cmtip]{xy}
\usepackage{graphicx}
\usepackage[numbers]{natbib}
\usepackage{bibentry}
\usepackage{hyperref}
\usepackage{comment}

\providecommand{\eqnref}[1]{(\ref{e:#1})}

\theoremstyle{plain}
\newtheorem{theorem}{Theorem}[section]
\newtheorem{lemma}[theorem]{Lemma}
\newtheorem{proposition}[theorem]{Proposition}

\newtheorem{fact}[theorem]{Fact}
\newtheorem*{fact*}{Fact}
\newtheorem{corollary}[theorem]{Corollary}
\newtheorem{claim}[theorem]{Claim}
\newtheorem*{claim*}{Claim}
\theoremstyle{definition}
\newtheorem{definition}[theorem]{Definition}
\newtheorem*{definition*}{Definition}

\newtheorem*{notation*}{Notation}
\theoremstyle{remark}
\newtheorem{remark}[theorem]{Remark}
\newtheorem*{remark*}{Remark}
\newtheorem{example}[theorem]{Example}
\newtheorem*{example*}{Example}

\newtheorem*{note*}{Note}
\newtheorem{question}[theorem]{Question}
\newtheorem*{question*}{Question}
\newtheorem{assumption}[theorem]{Assumption}
\begin{document}
\providecommand{\defn}[1]{{\bf #1}}
\providecommand{\bdl}{\boldsymbol\delta}
\renewcommand{\o}{\circ}
\newcommand{\A}{{\mathbb{A}}}
\newcommand{\C}{{\mathbb{C}}}
\newcommand{\R}{{\mathbb{R}}}
\newcommand{\N}{{\mathbb{N}}}
\newcommand{\Z}{{\mathbb{Z}}}
\newcommand{\Q}{{\mathbb{Q}}}
\newcommand{\K}{{\mathbb{K}}}
\newcommand{\G}{{\mathbb{G}}}
\newcommand{\U}{{\mathcal{U}}}
\newcommand{\M}{{\mathcal{M}}}
\newcommand{\fin}{{\operatorname{fin}}}
\newcommand{\tp}{{\operatorname{tp}}}
\newcommand{\Th}{{\operatorname{Th}}}
\newcommand{\cl}{{\operatorname{cl}}}
\newcommand{\trd}{{\operatorname{trd}}}
\newcommand{\acl}{{\operatorname{acl}}}
\newcommand{\dcl}{{\operatorname{dcl}}}
\newcommand{\cb}{{\operatorname{cb}}}
\newcommand{\locus}{{\operatorname{locus}}}
\newcommand{\loc}{{\operatorname{loc}}}
\newcommand{\Aut}{{\operatorname{Aut}}}
\newcommand{\Gal}{{\operatorname{Gal}}}
\newcommand{\Stab}{{\operatorname{Stab}}}
\newcommand{\ACF}{{\operatorname{ACF}}}
\newcommand{\GL}{{\operatorname{GL}}}
\newcommand{\AGL}{{\operatorname{AGL}}}
\newcommand{\PSL}{{\operatorname{PSL}}}
\newcommand{\st}{{\operatorname{st}}}
\newcommand{\x}{{\bar{x}}}
\newcommand{\y}{{\bar{y}}}
\renewcommand{\a}{{\bar{a}}}
\renewcommand{\b}{{\bar{b}}}
\renewcommand{\c}{{\bar{c}}}
\renewcommand{\d}{{\bar{d}}}

\let\polishL\L
\newcommand{\Los}{\polishL{}o\'s}
\renewcommand{\L}{{\mathcal{L}}}
\newcommand{\Lring}{\L_{\operatorname{ring}}}

\newcommand{\dind}{\ind^{\bdl}}
\newcommand{\ndind}{\nind^{\bdl}}

\renewcommand{\.}{\;}

\setcounter{tocdepth}{1}

\providecommand{\FIXME}[1]{({\bf \small FIXME: #1})}
\providecommand{\TODO}[1]{({\it \small TODO: #1})}

\newcommand\subqed[1]{\qedhere$_{\mathit{\ref{#1}}}$}
\newcommand\subqedinmath[1]{\qedhere_{\mathit{\ref{#1}}}}

\providecommand{\calS}{{\mathcal{S}}}
\providecommand{\calN}{{\mathcal{N}}}
\providecommand{\pow}{\mathcal{P}}
\providecommand{\eps}{\varepsilon}

\def\Ind#1#2{#1\setbox0=\hbox{$#1x$}\kern\wd0\hbox to 0pt{\hss$#1\mid$\hss}
\lower.9\ht0\hbox to 0pt{\hss$#1\smile$\hss}\kern\wd0}
\def\ind{\mathop{\mathpalette\Ind\emptyset }}
\def\notind#1#2{#1\setbox0=\hbox{$#1x$}\kern\wd0\hbox to 0pt{\mathchardef
\nn=12854\hss$#1\nn$\kern1.4\wd0\hss}\hbox to
0pt{\hss$#1\mid$\hss}\lower.9\ht0 \hbox to
0pt{\hss$#1\smile$\hss}\kern\wd0}
\def\nind{\mathop{\mathpalette\notind\emptyset }}

\newcommand{\grpconf}[6]{
 \[ \xymatrix@=0em@!{
   & & &&&& #3 \ar@{-}'[ddll][dddlll] \ar@{-}'[dlll][ddllllll] \\
   & & & #2  & & & \\
   #1 & & & & #6  & & \\
   & & & #4 && & \\
   && & & & & #5\;\; \ar@{-}'[uull][uuulll] \ar@{-}'[ulll][uullllll] } \]
 }

\title{An asymmetric version of Elekes-Szabó via group actions}
\author{Martin Bays and Tingxiang Zou}
\address{Martin Bays, Mathematical Institute, University of Oxford, Andrew Wiles Building, Radcliffe Observatory Quarter, Woodstock Road, Oxford OX2 6GG, UK}
\email{mbays@sdf.org}
\address{Tingxiang Zou, Mathematical Institute, University of Bonn, Endenicher Allee 60, 53115 Bonn, Germany}
\email{tzou@math.uni-bonn.de}
\thanks{Both authors were partially supported by DFG EXC 2044–390685587.}

\begin{abstract}
  We consider when finite families $F \subseteq \C[t]$ of bounded degree 
  polynomials, or more generally of bounded complexity finite-to-finite 
  correspondences on $\C$, can exhibit non-expansion of the form $|F(A)| = 
  O(|A|^{1+\eta})$ in their actions on finite sets $A \subseteq \C$ with $|F| 
  \gg |A|^\eps \gg 1$, for a fixed $\eps>0$ and arbitrarily small $\eta>0$.
  Our conclusions generalise the Elekes-Rónyai and Elekes-Szabó theorems, 
  which correspond to the case that $F$ is parametrised by a single complex 
  variable and $|F|=|A|$. Our result also applies to families of 
  correspondences between varieties of arbitrary dimension if we impose a 
  general position assumption on $A$. In all cases, the conclusion is that a 
  commutative algebraic group structure is responsible. As a special case, we 
  obtain asymmetric versions of Elekes-Rónyai and Elekes-Szabó, with explicit 
  bounds on exponents. Our methods originate in model theory.
\end{abstract}

\maketitle
\tableofcontents
\section{Introduction}
\label{s:intro}
A polynomial $f \in \R[x,y]$ in two variables is \defn{additive} if it is of 
the form $f(x,y) = g(h(x)+k(y))$ for some $g,h,k \in \R[t]$, and 
\defn{multiplicative} if it is of the form $g(h(x)\cdot k(y))$. The Elekes-Rónyai 
Theorem states that any other polynomial is \textit{expanding}, in the 
following sense:
\begin{fact}[Elekes-Rónyai \cite{ER},\cite{dZ-ER}]
    Let $f \in \R[x,y]$ be neither additive nor multiplicative.
    Then there exist $c,\eta>0$ (depending only on $\deg(f)$) such that for all
    $A,B \subseteq \R$ with $|A| = |B| = n \in \N$,
    $$|f(A,B)| = |\{ f(a,b) : a \in A,\; b \in B\}| \geq cn^{1+\eta}.$$
\end{fact}
Subsequently, explicit exponents $\eta$ were obtained, with $\eta=\frac12$ the 
current record \cite[Corollary~1.5]{SZ-ES}, and it was observed (see \cite[\S~1.2]{dZ-ER})
that $\R$ can be replaced with $\C$ (with $\eta=\frac13$).

A related theorem of Elekes and Szabó treats the more general setting where 
an arbitrary algebraic surface replaces the graph $z=f(x,y)$:
\begin{fact}[Elekes-Szabó \cite{ES-groups}]
    There is $\eta>0$ ($\eta=\frac16$ works) such that
    if $S \subseteq \C^3$ is an irreducible complex algebraic surface which does 
    not project in any two co-ordinates to a curve,
    and if for all $c>0$ there are $A_i \subseteq \C$ with $|A_i| 
    = n \in \N$ and $|S \cap (A_1\times A_2\times A_3)| \geq cn^{2-\eta}$,
    then $S \subseteq \C^3$ is in co-ordinatewise correspondence (see
    Definition~\ref{d:corr}) with the graph of the group operation of a
    1-dimensional algebraic group.
\end{fact}

A natural higher-dimensional generalisation would consider a subvariety $V 
\subseteq W_1\times W_2 \times W_3$ of a product of three complex varieties, 
and look for conditions under which $V$ admits  asymptotically large
intersections with products $A_1\times A_2 \times A_3$ of finite subsets $A_i
\subseteq W_i$.
One version of this is achieved in \cite[Theorem~27]{ES-groups}, but only 
under a strong \textit{general position} assumption on the $A_i$ (slightly 
weakened in \cite{BB-cohMod}).

The main result of this paper is a version which removes the general position 
assumption in one co-ordinate, say for $A_2$, and also allows $|A_2|$ to be 
much smaller than $|A_1| = |A_3|$.
To motivate this, we first present a special case generalising the 
Elekes-Rónyai theorem.

Given a polynomial $f(x,y) \in \C[x,y]$ and finite sets $A,B \subseteq \C$, we 
can view $B$ as defining a finite set $F := \{ f(x,b) : b \in B\} \subseteq 
\C[x]$ of polynomials in one variable, so $f(A,B) = F(A) := \{ f(a) : f 
\in F,\; a \in A \}$. This suggests a generalisation: given a finite set of 
polynomials $F \subseteq \C[x]$ and a finite set $A \subseteq \C$, when does 
$|F(A)|$ not expand? From this point of view, a constraint $|F| = |A|$ is 
unnatural, and we are led to an asymmetric view of expansion as meaning 
$|F(A)| \geq |A|^{1+\eta}$.
In the case $F = \{f(x,b) : b \in B\}$, the polynomials in $F$ are of degree 
bounded independently of $B$, and it turns out that this is all that is needed 
to obtain expansion from polynomials which are neither additive nor 
multiplicative.
More precisely, define $F$ to be \defn{$\eps$-additive}
if there are $g,h \in \C[t]$ such that
$| F \cap \{ g(h(x) + a) : a \in \C \} | \geq |F|^{1-\eps}$,
and define \defn{$\eps$-multiplicative} similarly.
Observe that if $f(x,y)$ is additive, then $F = \{f(x,b) : b \in B\}$
is 0-additive.
We show:

\begin{theorem}[Theorem~\ref{t:higherER}] \label{t:higherER_intro}
  For all $\eps>0$, there exist $c,\eta > 0$
  such that
  for any finite set $F \subseteq \C[x]$ of 
  non-constant polynomials of degree $\leq 1/\eps$ which is neither 
  $\eps$-additive nor $\eps$-multiplicative,
  and any finite set $A \subseteq \C$ with $|F|\geq |A|^\eps$,
  we have
  \[|F(A)| = |\{ f(a) : f \in F,\; a \in A\}|\geq c|A|^{1+\eta}.\]
\end{theorem}

In this generality, we prove only existence of $\eta$ without any bounds. This 
is a limitation due to our rather indirect proof method.

However, Theorem~\ref{t:higherER_intro} gives new information in the original 
Elekes-Rónyai context of $F = \{ f(x,b) : b \in B \}$, in the ``extremely 
unbalanced'' case that $B$ is very small compared to $A$, and in this case we 
do obtain bounds on $\eta$:
\begin{theorem}[Theorem~\ref{t:unbalancedERMany}] \label{t:unbalancedER_intro}
    Suppose $f\in\C[x,y]$ is neither additive nor multiplicative.
    For all $\eps>0$ there exists $C$
    such that for all finite $A,B\subseteq \C$ with $|A| \geq |B|\geq |A|^\eps \geq C$,
    we have $|f(A,B)|\geq |A|^{1+2^{-c\eps^{-1}}}$,
    where $c$ is an absolute constant.
\end{theorem}
Previous such results \cite{RSS-ER,SZ-ES} obtain this when $\eps>\frac12$
(with a better exponent when $\eps-\frac12 \gg 0$).

In Theorem~\ref{t:unbalancedERMany}, we moreover prove a multivariate version 
of Theorem~\ref{t:unbalancedER_intro}:
for a polynomial $f\in\C[x,\y]$, we obtain $|f(A,B_1,...,B_n)| \geq 
|A|^{1+n^{-2}2^{-c\eps^{-1}}}$ for $|A| \geq |B_i| \geq |A|^\eps\geq C$, 
unless $f$ is of the form $g(h(x)+k(\y))$ or $g(h(x)\cdot k(\y))$.
This can be seen as an asymmetric version of the main result of \cite{RazTov}, 
which with $|B_i|=|A|$ obtains bounds of the form 
$\Omega(|A|^{1+\frac12})$ unless $f$ is of the form $h(\sum_i g_i(x_i))$ or 
$h(\prod_i g_i(x_i))$.

We also obtain a similarly unbalanced version of Elekes-Szabó:
\begin{theorem}[{Theorem~\ref{t:main_fin_1d}}]
  Suppose $S \subseteq \C^3$ is an irreducible complex algebraic surface which does 
  not project in any two co-ordinates to a curve,
  and $S \subseteq \C^3$ is not in co-ordinatewise correspondence
  with the graph of the group operation of a 1-dimensional algebraic group.
  Then for all $\eps > 0$ there exists $C>0$ such that
  for all finite subsets $A_i \subseteq \C$ with
  $|A_1| = |A_3|$ and $|A_1|^{\frac1\eps} \geq |A_2| \geq |A_1|^\eps \geq C$, 
  we have
  $$|S \cap (A_1\times A_2\times A_3)| \leq (|A_1||A_2|)^{1-\eta_0},$$
  where
  $\eta_0 = \eps(1+\frac1\eps)^{-1}2^{-\frac4{c\eps}+7}$ and $c$ is an 
  absolute constant.
\end{theorem}

We deduce the above results as consequences of our main theorem, which is a 
similarly asymmetric and partially higher dimensional generalisation of the 
Elekes-Szabó theorem.
Our approach to this is model-theoretic, continuing the line of investigation 
developed in \cite{HrW-LP}, \cite{Hr-psfDims}, \cite{BB-cohMod}, \cite{bgth}, 
and \cite{cubicSurfaces}. This involves finding, as a witness to the 
asymptotically large intersections considered in Elekes-Szabó, a type in a 
certain expansion of the field language by pseudofinite sets.
Our result is then most naturally formulated in terms of the properties of a 
realisation of this type, as follows; the technical terms and notation
used here
will be defined in the body of the paper. An elementary finitary statement 
appears as Theorem~\ref{t:main_fin}.
\begin{theorem}[Theorem~\ref{t:main}] \label{t:main_intro}
  Suppose $a$, $d$, and $b$ are broad tuples in an adequate expansion of the field $\C^\U$ such that:
\begin{enumerate}[(i)]
    \item $\tp(a)$ is in coarse general position;
    \item $\tp(d)$ is in weak general position;
    \item $a \sim_d b$ (meaning $a$ is field-theoretically interalgebraic 
      with $b$ over $d$);
    \item $d \sim \cb^{\ACF}(ab/d)$ (i.e.\ $d$ is (almost) canonical for the 
      interalgebraicity);
    \item $a \dind d$ and $d \dind b$.
  \end{enumerate}
  Then there exists a connected commutative algebraic group $G$ and a generic 
  pair $(g_a,g_d) \in G^2$, such that
  $a \sim g_a$ and $d \sim g_d$ and $b \sim g_d+g_a$.
\end{theorem}

A key idea in our proof is that the conditions (iii)-(v), with algebraic 
independence in place of $\dind $, are what we expect not just of a 
triple $b = d+a$ in a commutative algebraic group as in the conclusion, but 
more generally of a generic triple $b = d*a$ where $*$ is a faithful 
transitive action of an algebraic group on a variety. Our proof goes via 
invoking the suitable version of the group configuration theorem to find such 
a group action, then using the coarse general position assumption to show that 
the group must actually be abelian, and the action therefore trivial.
When $\tp(d)$ is in coarse general position (as in
Theorem~\ref{t:unbalancedER_intro}), we obtain the group configuration using
similar ideas to those applied in \cite{BB-cohMod} to the original
Elekes-Szabó situation, but with a finite iteration required to get to the
group configuration -- see Proposition~\ref{p:maincgp}. More generally,
without this coarse general position assumption, we first obtain a family of
such groups by adding parameters to reduce to this case, then argue that this
family collapses to a single group.

As an application of some of the ideas in the proof, we obtain a group action
version Theorem~\ref{c:BSGActionFin} of the Balog-Szemerédi-Gowers theorem,
linking non-expansion conditions of the form $|S*A| \leq O(|A|^{1+\eps})$ to
approximate subgroups (where $*$ is a faithful action of a group $G$ on a set
$X$ and $S \subseteq G$ and $A \subseteq X$ are finite subsets). This is a
qualitative generalisation of the asymmetric Balog-Szemerédi-Gowers theorem
\cite[Theorem~2.35]{TaoVu}; similar results were previously obtained in \cite{murphy}.

\subsection{Structure of the paper}
\S\S\ref{s:prelim}-\ref{s:adeq} develop the technical tools on pseudofinite dimensions we use in proving our main results:
\S\ref{s:prelim} recalls the basic definitions,
\S\ref{s:bsg} develops a version in our context of Balog-Szemerédi-Gowers,
and \S\ref{s:adeq} handles limits of types.
\S\ref{s:cohHom} then proves that any faithful algebraic group action $(G,X)$ containing asymptotically non-expanding finite subsets $(S,A)$ must, if $A$ is in suitably general position, be abelian.
\S\ref{s:CFAHS} explains how to recognise an algebraic group action from a generic algebraic configuration, specialising the group configuration theorem from model theory.
In \S\ref{s:main} we prove our main result Theorem~\ref{t:main} (introduced as Theorem~\ref{t:main_intro} above); this goes via first proving a weaker version, Proposition~\ref{p:maincgp}, which has a stronger general position assumption (this weaker version is already sufficient for the 1-dimensional asymmetric Elekes-Szabó and Elekes-Rónyai results mentioned above).
\S\ref{s:fin} spells out the consequences of Theorem~\ref{t:main} for finite sets.
Next, \S\ref{s:ER} proves our generalisation of Elekes-Rónyai, Theorem~\ref{t:higherER} (introduced as Theorem~\ref{t:higherER_intro} above).
Appendix~\ref{s:eps} repeats some of our earlier qualitative arguments in a quantitative way, perturbing our assumptions and tracing through the effect on our results. This is used to obtain the bounds on the exponents $\eta$ in some of our results, as presented above.
Finally, Appendix~\ref{s:BSGActions}, which depends only on Sections~\ref{s:prelim}-\ref{s:adeq}, proves our group action version Theorem~\ref{c:BSGActionFin} of Balog-Szemerédi-Gowers.

\subsection{Acknowledgements}
We would like to thank Emmanuel Breuillard and Udi Hrushovski for helpful 
conversations. We would also like to thank the LMS and the HCM for supporting 
research visits during which some of this work was done.
Particular thanks to the anonymous referees for their thorough reading of the 
manuscript and useful suggestions, which substantially improved the
presentation of the paper, and in particular led to cleaner arguments in
Lemma~\ref{l:adequacy}, Claim~\ref{c:energyIndep}, and the end of
Proposition~\ref{p:maincgp}.

\section{Preliminaries}
\label{s:prelim}
\subsection{Notation}
We use some abbreviations which are standard in model theory.
If $a$ and $b$ are tuples, then $ab$ denotes the concatenated tuple $a^\smallfrown 
b$, while if $A$ and $B$ are sets, $AB$ denotes $A\cup B$. We write $a \in A$ 
to mean that each element of the tuple $a$ is an element of the set $A$. All 
tuples have finite length unless otherwise mentioned. Sometimes a tuple is
treated as the set of its elements.

Let $M$ be an $\L$-structure.
If $C \subseteq M$, then $M^{\L(C)}$ is the expansion of $M$ to $\L(C) := \L 
\cup C$, adding constants interpreted as the elements of $C$.
A \emph{realisation} in $M$ of a set of formulas $\Phi$ in variables $(x_i)_{i 
\in I}$ is a (possibly infinite) sequence $\a = (a_i)_{i \in I}$ such that $M 
\vDash \phi(\a)$ for each $\phi \in \Phi$. Such a set $\Phi$ is 
\emph{consistent} if every finite subset of $\Phi$ has a realisation in $M$.
The \emph{type} of a tuple $a$ over a set $C \subseteq M$ is the maximal set of 
$\L(C)$-formulas it realises, $\tp(a/C) = \tp^\L(a/C) = \{ \phi(x,c) : M \vDash \phi(a,c);\; 
c\in C\}$.
We write $a \equiv_C b$ to mean $\tp(a/C) = \tp(b/C)$.
When $C=\emptyset$, it is omitted.

A tuple $a$ is \emph{definable} over $C$ if $a$ is the unique realisation in $M$ of some $\L(C)$-formula $\phi(x)$.
Let $\dcl(C)\subseteq M$ be the set of those elements which are definable over $C$.
Tuples $a$ and $b$ are \emph{interdefinable} over $C$ if $\dcl(Ca)=\dcl(Cb)$.

If $R$ is a binary relation, such as $\equiv_C$ or $\dind_C$, we sometimes write
e.g.\ $aRbRc$ as shorthand for ``$aRb$ and $bRc$''.

\subsection{Setup}
Let $\L$ be a countable first-order language.
Let $M=\prod_{i\to\mathcal{U}}M_i$ be an ultraproduct of $\L$-structures $M_i$ 
over some non-principal ultrafilter $\mathcal{U}$ on $\N$.

\begin{fact}[{\cite[Exercise~4.5.37, Lemma~4.3.17]{Marker-mt}}] \label{f:sat}
  $M$ is $\aleph_1$-saturated: if $C \subseteq M$ is a countable subset, then 
  any consistent set $\Phi((x_i)_{i\in\aleph_1})$ of $\L(C)$-formulas has a 
  realisation in $M$.
\end{fact}

A subset $X \subseteq M^n$ is \defn{internal} if $X = \prod_{i\to\mathcal 
U}X_i$ for some subsets $X_i \subseteq M_i^n$.

Let $\xi \in \N^\U \setminus \N$.
Recall from \cite[\S2.1, \S2.2, Fact~2.2]{cubicSurfaces}
the following definitions of $\bdl = \bdl_\xi$ on internal sets and
on tuples.
For $X$ an internal set, we set $\bdl(X) = \bdl_\xi(X) := 
\st \log_\xi |X| \in \R \cup \{-\infty,\infty\}$,
where $\st : \R^\U \cup \{-\infty,\infty\} \to \R \cup \{-\infty,\infty\}$ is 
the standard part map.
In other words, $\bdl(\prod_{i\to\U} X_i)$ is (when finite) the metric 
ultralimit in $\R$ along $\U$ of $\log_{\xi_i}(|X_i|)$, where $\xi = 
(\xi_i)_{i\in\omega}/\U$.

For $C \subseteq M$ countable and $a$ a tuple from $M$,
let $$\bdl(a/C) = \bdl^\L(a/C) := \inf_{\phi \in \tp^\L(a/C)} \bdl(\phi(M)).$$
We include $\L$ in the notation only when we need to make it explicit.
More generally, for $p$ a partial type we define $\bdl(p):=\inf_{\phi \in 
p} \bdl(\phi(M))$. By a \defn{$\bigwedge$-definable} set over $C$, we mean
a countable intersection $X = \bigcap X_i$ of definable sets over $C$, and we
similarly define $\bdl(X) = \inf_i \bdl(X_i)$; this is well-defined, see
\cite[\S2.1.5]{BB-cohMod}.

If $\bdl(a/C), \bdl(b/C) < \infty$,
we write $a \dind_C b$ if $\bdl(a/C) = \bdl(a/Cb)$,
and $a\ndind_C b$ otherwise.
More generally, if $A,B,C$ are such that for all finite tuples $a \in A$ and
$b \in B$ we have $\bdl(a/C),\bdl(b/C) < \infty$ and $a\dind_C b$, we write $A
\dind_C B$.

We will always work in countable languages in which $\bdl$ is continuous (meaning that $\tp(b) \mapsto \bdl(\phi(M,b))$ is a well-defined continuous map for each $\L$-formula $\phi(x,y)$). Then $\dind$ is a well-behaved independence notion; we recall here \cite[Fact~2.2]{cubicSurfaces} for reference.
\begin{fact}\label{f:bdlInd}
	Let $a,b,d$ be tuples from $\K$ and $C \subseteq  \K$ a countable subset,
	and suppose $\bdl(a/C),\bdl(b/C),\bdl(d/C) < \infty$
	\begin{itemize}
		\item $\trd(a/C)=0\Rightarrow \bdl(a/C)=0$ \hfill{(Algebraicity)}
		\item $a \equiv _C b \Rightarrow  \bdl(a/C) = \bdl(b/C)$ \hfill{(Invariance)}
		\item $\bdl(ab/C) = \bdl(a/Cb) + \bdl(b/C)$ \hfill{(Additivity)}
		\item $a \ind ^{\bdl}_C b \Leftrightarrow  b \ind ^{\bdl}_C a$ \hfill{(Symmetry)}
		\item $ab \ind ^{\bdl}_C d \Leftrightarrow  (b \ind ^{\bdl}_C d \wedge a \ind ^{\bdl}_{Cb} d)$ 
		\hfill{(Monotonicity and Transitivity)}
		\item For any tuple $e$ from $\K$, there exists $a' \equiv _C a$ such that $a' 
		\ind ^{\bdl}_C e$; \hfill{(Extension)}\\
                more generally, any $\bigwedge$-definable set $X$ over $C$ 
                contains an element $x$ with $\bdl(x/C) = \bdl(X)$.
	\end{itemize}
\end{fact}

Note that from additivity, it is easy to see that $\bdl(a/C)=\bdl(a/Cc)$ for any $c$ with $\bdl(c/C)=0$. Thus, $a\ind ^{\bdl}_C b  \Rightarrow a\ind ^{\bdl}_C b c$ for any $c$ with $\bdl(c/Cb)=0$.

Note also that if $a\in\acl_{\L}(Cb)$, then $\bdl(a/Cb) = 0$,
and so by additivity $\bdl(b/C) - \bdl(b/Ca) = \bdl(a/C) - \bdl(a/Cb) = \bdl(a/C)$,
and in particular $\bdl(a/C)\leq\bdl(b/C)$.

A type $p$ is \defn{broad} if $0 < \bdl(p) < \infty$; similarly for internal sets.

\begin{definition}
  A \defn{$\bdl$-independent sequence} in a type $\tp(a/C)$ with $\bdl(a/C) < \infty$
  is a (possibly infinite) sequence $\a = (a_i)_{i<\alpha}$
  of realisations of $\tp(a/C)$
  such that $a_i \dind_C a_{<i}$ for all $i<\alpha$, where $a_{<i} = 
  (a_j)_{j<i}$.
\end{definition}

The following fact and the lemma below are direct applications of $\aleph_1$-saturation of $M$.
\begin{fact}\label{f:existsInternal}\cite[Section 1.1]{bgth}
Let $p$ be a partial type over a countable set of parameters. Then there is an internal set $Y\subseteq p(M)$ such that $\bdl(Y)=\bdl(p)$.
\end{fact}

\begin{lemma}\label{l:bdlf(p)}
        Let $f:X_1\times\cdots\times X_n\to Z$ be a definable function and 
        $p_i$ a partial type in $X_i$ for $1\leq i\leq n$ defined over a 
        countable set $C$. Then \[\bdl(f(p_1(M),\ldots, 
        p_n(M)))=\max\{\bdl(f(a_1,\ldots, a_n)/C): a_i\in p_i(M), \text{ for 
        }1\leq i\leq m\}.\]
\end{lemma}
\begin{proof}
	Note that by $\aleph_1$-saturation, \[f(p_1(M),\ldots, 
        p_n(M))=\bigcap_{\phi_i\in p_i}f(\phi_1(M),\ldots,f_n(M))\] is 
        $\bigwedge$-definable, so the result follows from the ``more generally'' clause in Fact~\ref{f:bdlInd}(Extension).
      \end{proof}

\section{Balog-Szemerédi-Gowers for types}
\label{s:bsg}
Given subsets $A$ and $B$ of a group such that a large subset $Z$ of $A\times B$ has few distinct products $a\cdot b$ for $(a,b)\in Z$,
the Balog–Szemerédi–Gowers (BSG) theorem \cite[Theorem~2.44]{TaoVu} extracts a reasonably large rectangle $A'\times B' \subseteq A\times B$ with few distinct products. Further arguments of Tao, \cite[Theorem~5.4]{Tao-BSG}, go on to extract an approximate subgroup in this situation. In Proposition~\ref{p:BSGT}, we give qualitative analogues of these results which work with complete types, passing from a non-expanding 2-type to a non-expanding product of 1-types over a $\bdl$-independently extended base, and extracting a type for which arbitrarily long products still don't expand.
For context, we then explain in \S\ref{ss:BSG} how these can be seen as incarnations of the original results, by using them to to recover qualitative versions of those results.

Let $\L$, $M$, and $\bdl=\bdl_{\xi}$ be as above, and assume that $\bdl^\L$ is continuous.
Let $(G,\cdot)$ be an $\L$-definable group. We use $g\cdot h$ to denote the product of $g$ and $h$ in $G$, while $gh$ continues to denote the tuple $(g,h)$.

\begin{definition}
  A set $D \subseteq M$ is \defn{locally $\bdl$-finite} if $\bdl(d)<\infty$ for all tuples $d\in D$.
\end{definition}
\begin{proposition} \label{p:BSGT}
  Let $g,h \in G$ be such that $\bdl(g),\bdl(h)<\infty$,
  and $$g \dind  h \dind  g\cdot h \dind  g.$$

  Then there exists a countable locally $\bdl$-finite $D$ with $D \dind gh$
  such that:
  \begin{enumerate}[(i)]
    \item If $g' \equiv_D g$ and $h' \equiv_D h$, then $\bdl(g'\cdot h') \leq \bdl(g)$.
    \item If $n\in\N$ and $g_i \equiv_D g$ for $i=1,\ldots ,n$,
      and if $g_i \dind_D  g_{i+1}$ for odd $i\leq n$,
      then $\bdl(\prod_{i\leq n} g_i^{(-1)^i}/D) \leq  \bdl(g)$.
  \end{enumerate}
\end{proposition}
\begin{remark*}
  In particular, (ii) applies to an arbitrary word in realisations of $\tp(g'^{-1}\cdot g/D)$, where $g' \equiv_D g$ and $g' \dind_D g$.
\end{remark*}

We delay the proof of Proposition~\ref{p:BSGT}, and first obtain a consequence for group actions, which is what we will actually use in the main text (specifically, in Theorem~\ref{t:cohActAb}).

\begin{corollary}\label{c:symBSGact}
Let $(G,X)$ be an $\L$-definable group action.

Let $g,h\in G$ and $a\in X$ be such that $\bdl(g),\bdl(h),\bdl(a)<\infty$ and 
\[g\dind h\dind g\cdot h\dind g \text{ and }a\dind g\dind g*a.\]
Then there exists a countable locally $\bdl$-finite set $D$ with $D\dind gha$ 
such that
for any even-length alternating product $s = \prod_{1\leq i \leq 2n}g_i^{(-1)^i}$ with $g_i \equiv_D g$ for all $i\leq 2n$ and with $g_{2m-1} \dind_D g_{2m}$ for all $m\leq n$,
and any $a'\equiv_D a$,
we have
\begin{itemize}
\item 
$\bdl(s/D) \leq  \bdl(g)$;
\item 
$\bdl(s*a'/D) \leq  \bdl(a)$.
\end{itemize}
\end{corollary}
\begin{proof}
 Since $g\dind h\dind g\cdot h\dind g$, we can apply Proposition~\ref{p:BSGT}. Let $D\dind g$ be the locally $\bdl$-finite set given by the conclusion of  Proposition~\ref{p:BSGT}. Replacing $D$ with a copy independent with $ha$ over $g$, we may assume $D\dind gha$. Add $D$ to the constants.

 Let $s$, $g_i\equiv g$, and $a'\equiv a$ be as in the statement.
 Then $\bdl(s) \leq  \bdl(g)$ by Proposition~\ref{p:BSGT}(ii); it remains to show $\bdl(s*a') \leq  \bdl(a)$.
	
	Let $g_{2n+1}a'\equiv ga$ be such that $g_{2n+1}\dind_{a'}s$.
        Since $g\dind a$ and hence $g_{2n+1}\dind a'$, we have $g_{2n+1}\dind sa'$.
        Now $s' := s\cdot g_{2n+1}^{-1}=(\prod_{i\leq 2n+1}g_i^{(-1)^i})$, so $\bdl(s') \leq \bdl(g)$ by Proposition~\ref{p:BSGT}(ii).
        Also $\bdl(g_{2n+1}*a') = \bdl(g*a) = \bdl(a)$.
        So using $s*a' = s'*(g_{2n+1}*a')$, we have
	\begin{align*}
		\bdl(s*a')&=\bdl(s',g_{2n+1}*a')-\bdl(s',g_{2n+1}*a'/s*a')\\
		&\leq \bdl(g)+\bdl(a)-\bdl(s'/s,a')\\
                &= \bdl(g)+\bdl(a)-\bdl(g_{2n+1}/s,a')\\
                &= \bdl(a).\qedhere
	\end{align*}
\end{proof}

\newcommand{\fsind}{\ind^u}
\newcommand{\nfsind}{\nind^u}
The remainder of this section is devoted to the proof of 
Proposition~\ref{p:BSGT}. The definitions and results 
introduced in the process will not be used in the rest of the paper.

We work in $M$: tuples are tuples from $M$, and sets of parameters are 
\emph{countable} subsets of $M$; we use Fact~\ref{f:sat} without mention. Let $E_0\subseteq M$ be the set of interpretations of the constants of $\L$.

Recall that a type $\tp(a/B)$ is
\begin{itemize}
  \item \emph{finitely satisfiable in $C$} if for any formula $\phi(x,b) \in 
    \tp(a/B)$, there is $c \in C$ such that $M \vDash \phi(c,b)$;
  \item \emph{$C$-invariant} if for any tuples $b,b' \in B$ with $b \equiv_C 
    b'$, we have $ab \equiv_C ab'$.
\end{itemize}
We write $a \fsind_C B$ to mean that $\tp(a/CB)$ is finitely satisfiable in 
$E_0C$.

\begin{fact}[{\cite[Example~2.17]{simon-NIP}}] \label{f:factFS}
  If $a \fsind_C B$, then $\tp(a/CB)$ is $C$-invariant,
  and for any $B'$ there is $a' \equiv_{CB} a$ such that $a' \fsind_C BB'$.
\end{fact}

We first show, in Lemma~\ref{l:semiregularisation}, that $a\dind_C b$ can always be improved to $a\fsind_C b$ by working over additional $\dind_C$-independent parameters, uniformly in $a$.

\begin{definition} \label{d:chgp}
  A type $\tp(b/C)$ with $\bdl(b/C) < \infty$ is \defn{in coheir general position} (\defn{chgp}) if 
  for any $a$ with $\bdl(a/C) < \infty$,
  $$ a \dind_C b \Rightarrow  a \fsind_C b .$$
\end{definition}

\begin{lemma} \label{l:chgpProps}
  \begin{enumerate}[(a)]\item \label{chgpdcl} If $b$ and $b'$ are interdefinable over $C$, then 
  $\tp(b/C)$ is chgp if and only if $\tp(b'/C)$ is.
  \item \label{chgpConv} The converse implication $a \ndind_C b \Rightarrow  a\nfsind_C b$ 
    holds for any $a,b,C$ with $\bdl(ab/C) < \infty$.
  \end{enumerate}
\end{lemma}
\begin{proof}
  \begin{enumerate}[(a)]\item $\tp(a/bC)$ is finitely satisfiable in $E_0C$ if 
      and only if $\tp(a/\dcl(bC))$ is, and the result follows.
  \item This follows from continuity of $\bdl^\L$.
  If $a \fsind_C b$ but $a \ndind_C b$, then $b \ndind_C a$, so say 
  $\bdl(\phi(x,a)) < \bdl(b/C) =: \alpha$ for some $\L(C)$-formula $\phi(x,a) 
  \in \tp(b/aC)$,
  hence by continuity there is some $\psi(y)\in\tp(a/C)$ such that $\bdl(\phi(x,d))<\alpha$ whenever $\psi(d)$ holds. Since $\phi(b,y)\land \psi(y)\in\tp(a/Cb)$, by finite satisfiability already $\phi(b,c)\land\psi(c)$ holds for some $c \in (E_0C)^{|a|}$, hence $\bdl(\phi(x,c)) < 
  \alpha$  and $\phi(x,c) \in \tp(b/C)$, 
  contradicting $\alpha = \bdl(b/C)$.
  \qedhere
  \end{enumerate}
\end{proof}

\begin{lemma} \label{l:semiregularisation}
  Given $b$ with $\bdl(b) < \infty$, there exists a countable locally $\bdl$-finite $D$ such that
  $D \dind  b$ and $\tp(b/D)$ is chgp.
\end{lemma}
\begin{proof}
  Recursively define countable sets $D_i$ for $i \in \omega$ as follows.
  Let $D_0 := \emptyset$.
  Given $D_i$ with $D_i \dind  b$, enumerate as $(\phi^i_j)_{j \in \omega}$ those 
  formulas $\phi(x)$ over $bD_i$ for which there exists $a$ such that $\vDash  
  \phi(a)$ and $\bdl(a)<\infty$ and $a \dind_{D_i} b$.
  If $\phi^i_j$ and $a$ are such, then let $d^i_j$ realise a $\bdl$-independent extension of $\tp(a/D_ib)$ to $D_ib(d^i_{j'})_{j'<j}$,
  so that $\vDash  \phi(d^i_j)$ and $\bdl(d^i_j)<\infty$ and $d^i_j \dind_{D_i} b(d^i_{j'})_{j'<j}$.
  Let $D_{i+1} := D_i \cup (d^i_j)_{j \in \omega}$.
  Then $D_{i+1} \dind  b$ (by monotonicity and transitivity).

  Let $D := \bigcup_{i \in \omega} D_i$. Then $D \dind  b$ by finite character,
  and $D$ is countable and locally $\bdl$-finite by construction.
  We conclude by verifying that $\tp(b/D)$ is chgp.
  So suppose $a$ is such that $\bdl(a/D)<\infty$ and $a \dind_D b$, and let $\phi(x) \in \tp(a/bD)$.
  Then $\phi(x) \in \tp(a/bD_i)$ for some $i$,
  and $a \dind_{D_i} b$ since $b \dind  Da$,
  and $\bdl(a)<\infty$ because $\bdl(a/d)<\infty$ for some tuple $d \in D$,
  so by construction there is $d \in D_{i+1}$ such that $\phi(d)$, as required.
\end{proof}


\begin{lemma} \label{l:BSG}
  Let $g,h \in G$ be such that $\bdl(g),\bdl(h)<\infty$,
  and $g \dind  h \dind  g\cdot h \dind  g$,
  and $g \fsind h$.

  Let $g' \equiv  g$ and $h' \equiv  h$.
  Then $\bdl(g'\cdot h') \leq  \bdl(g)$.
\end{lemma}
\begin{proof}
  We first find $g''$ and $h''$ such that
  $g''h'' \dind  g'h'$ and $g'h'' \equiv  g''h'' \equiv  g''h' \equiv  gh$
  (This can be seen as an abstract form of the bound on ``paths of length 
  three'' in \cite[6.20]{TaoVu}).

  Let $h''$ be such that $g'h'' \equiv  gh$ and $h'' \dind_{g'} h'$.
  Then $h'' \dind  g'h'$.
  Let (by Fact~\ref{f:factFS}) $g''$ be such that $g''h' \equiv  gh$ and $g'' 
  \fsind h'h''g'$.
  Then $g''h'' \equiv  g''h'$ by invariance of finitely satisfiable types,
  and $g'' \dind  h'h''g'$ by Lemma~\ref{l:chgpProps}(\ref{chgpConv}).
  Then $g''h'' \dind  g'h'$ (by symmetry, monotonicity, and transitivity).

  Now let \mbox{$e := (g'\cdot h'', g''\cdot h'', g''\cdot h')$}.
  We have $\bdl(g\cdot h) = \bdl(g\cdot h/h) = \bdl(g/h) = \bdl(g)$, and 
  similarly $\bdl(h)=\bdl(h/g\cdot h)=\bdl(g/g\cdot h)=\bdl(g)$.
  So $\bdl(e) \leq  3\bdl(g)$,
  and $g'\cdot h' = (g'\cdot h'')\cdot(g''\cdot h'')^{-1}\cdot (g''\cdot h') 
  \in \dcl(e)$,
  and $\bdl(e/g'\cdot h') \geq  \bdl(e/g',h') = \bdl(g'',h''/g',h') = 
  \bdl(g'',h'') = 2\bdl(g)$.

  Hence $\bdl(g'\cdot h') = \bdl(e) - \bdl(e/g'\cdot h') \leq  3\bdl(g) - 
  2\bdl(g) = \bdl(g)$.
\end{proof}

\begin{proof}[Proof of Proposition~\ref{p:BSGT}]
  By Lemma~\ref{l:semiregularisation}, there exists a countable locally finite $D'$ such that $D'\dind g$ and $\tp(g/D')$ is chgp. Let $D\equiv_g D'$ be such that $D\dind_g h$. Then $\tp(g/D)$ is chgp and $D\dind gh$ by $D\dind g$ and transitivity.
  So we may expand the language by constants for $D$, since this does not affect the assumptions. Recalling our convention about constants in the definition of $\fsind$, we then have that $\tp(g)$ is chgp. We finish the proof by obtaining conclusions (i) and (ii) of the proposition with $D=\emptyset$, using that $\tp(g)$ is chgp.

  Since $\tp(g)$ is chgp, we have $g\fsind h$, so (i) follows immediately from Lemma~\ref{l:BSG}.

  For (ii), we actually prove by induction on $n$ the slightly more general statement that if $g_0=1$ or $g_0\equiv g$, and $g_i \equiv g$ for $i>0$, then $\bdl(\prod_{0\leq i \leq n} g_i^{(-1)^i}) \leq \bdl(g)$.
  \begin{claim}\label{clm:BSGT}
    \begin{enumerate}[(a)]\mbox{}
      \item $\bdl(g_1^{-1}\cdot g_2) \leq  \bdl(g)$.
      \item $\bdl(g_0\cdot g_1^{-1}) \leq  \bdl(g)$.
    \end{enumerate}
  \end{claim}
  \begin{proof}
    Let $g'h' \equiv_{g\cdot h} gh$ with $g'h' \dind_{g\cdot h} gh$.
    Then $h' \dind  g\cdot h$, so $h' \dind  gh$, and then since $g \dind  h$ we 
    have $hh' \dind  g$.
    Similarly, $hh' \dind  g'$, and also $g \dind  g'$.
    Now $g^{-1}\cdot g' = h\cdot h'^{-1}$
    since $g'\cdot h' = g\cdot h$,
    hence $g^{-1} \dind  g' \dind  g^{-1}\cdot g' \dind  g^{-1}$.
    By chgp we have $g^{-1} \fsind g'$,
    so (a) follows from Lemma~\ref{l:BSG}.

    Now note that $(g^{-1},g\cdot h)$ also satisfies the conditions on $(g,h)$ used in this argument,
    namely that $g^{-1} \dind  g\cdot h \dind  h = g^{-1}\cdot g\cdot h \dind  g^{-1}$
    and $\tp(g^{-1})$ is chgp (by Lemma~\ref{l:chgpProps}(\ref{chgpdcl})).
    So (b) holds by the same argument.
    \subqed{clm:BSGT}
  \end{proof}

  By Claim~\ref{clm:BSGT}(b), we may assume $n\geq 2$.
  Let $w := \prod_{1\leq i\leq n-2} g_{i+2}^{(-1)^i}$,
  so
  $ b := \prod_{0\leq i\leq n} g_i^{(-1)^i} = g_0\cdot g_1^{-1}\cdot g_2\cdot w $.
  By the inductive hypothesis for $n-2$, we have $\bdl(g_2\cdot w) \leq  \bdl(g)$.

  Replacing $(g_1,g_2)$ with another realisation of $\tp(g_1,g_2/g_1^{-1}\cdot 
  g_2)$ (so $g_1^{-1}\cdot g_2$ and hence $b$ don't change),
  we may assume $g_1 \dind_{g_1^{-1}\cdot g_2} g_0b$.
  But $g_1 \dind  g_1^{-1}\cdot g_2$,
  since $\bdl(g_1^{-1}\cdot g_2/g_1) = \bdl(g_2/g_1) = \bdl(g) \geq  
  \bdl(g_1^{-1}\cdot g_2)$ (using $g_1 \dind  g_2$ and Claim~\ref{clm:BSGT}(a)).
  Hence $g_1 \dind  g_0b$. Then, using also Claim~\ref{clm:BSGT}(b),

  \begin{align*} \bdl(b) &= \bdl(g_0\cdot g_1^{-1},g_2\cdot w) - \bdl(g_0\cdot g_1^{-1},g_2\cdot w/b) \\
  &\leq  \bdl(g_0\cdot g_1^{-1}) + \bdl(g_2\cdot w) - \bdl(g_0\cdot g_1^{-1}/g_0,b) \\
  &\leq  2\bdl(g) - \bdl(g_1/g_0,b) \leq  \bdl(g)
  .\qedhere\end{align*}
\end{proof}

\subsection{Recovering BSG}\label{ss:BSG}
We now show how to deduce from Proposition~\ref{p:BSGT} 
qualitative versions of the Balog–Szemerédi–Gowers theorem and Tao's 
corollary on the existence of approximate subgroups. In Appendix~\ref{appB}, 
we will generalise these results to a group action version. Neither this
subsection nor that appendix are used elsewhere in the text.

\begin{definition}
Let $G$ be a definable group in $M$. We call an internal subset $X\subseteq G$ a \textbf{coarse approximate subgroup} if $1\in X$, $X=X^{-1}:=\{g^{-1}: g\in X\}$ and there is an internal subset $Z\subseteq G$ with $\bdl(Z)=0$ such that $X\cdot X\subseteq Z\cdot X$.
\end{definition}

Approximate subgroups have small doubling, namely the size of $X\cdot X$ is roughly the same as $X$. Conversely, if a set has small tripling, it is not difficult to show one can obtain an approximate subgroup from it. The bounds involved are polynomial, which means we express this result in terms of $\bdl$ as follows.
\begin{fact}\label{f:approx}\cite[Proposition 2.40]{TaoVu}
Suppose $A$ is an internal subset of $G$ such that $\bdl(A\cdot A\cdot A)=\bdl(A)$. Let $H':=A\cup\{1\}\cup\{A^{-1}\}$ and $H:=H'\cdot H'\cdot H'$. Then $H$ is a coarse approximate subgroup containing $A$ with $\bdl(H)=\bdl(A)$.
\end{fact}

\begin{corollary}\label{c: QBSGT}(Qualitative Balog–Szemerédi–Gowers-Tao theorem)
	Let $A,B\subseteq G$ and $Z\subseteq A\times B$ be internal subsets such that $0<\bdl(A), \bdl(B)<\infty$,  $\bdl(Z)=\bdl(A)+\bdl(B)$ and $\bdl(A\cdot^{Z}B)\leq\frac{1}{2}\bdl(A)+\frac{1}{2}\bdl(B)$, where $A\cdot^{Z}B:=\{a\cdot b: (a,b)\in Z\}$. Then the following hold.

	\begin{enumerate}
	\item (\cite[Theorem 2.44]{TaoVu})
	There are internal sets $A'\subseteq A$ and $B'\subseteq B$ such that $\bdl(A')=\bdl(A)$, $\bdl(B')=\bdl(B)$, and $\bdl(A'\cdot B')=\frac{1}{2}\bdl(A)+\frac{1}{2}\bdl(B)$.
	\item (\cite[Theorem 2.48]{TaoVu})
	There exist a coarse approximate subgroup $H$ in $G$ and $g,h\in G$ such that $\bdl(H)\leq \frac{1}{2}\bdl(A)+\frac{1}{2}\bdl(B)$ and $\bdl(A\cap (g\cdot H))\geq \bdl(H)$, $\bdl(B\cap (H\cdot h))\geq \bdl(H)$.
	\end{enumerate}
	
        In fact, we also have $\bdl(A)=\bdl(B)=\bdl(A'\cdot B')=\bdl(H)$.
\end{corollary}
\begin{remark}
  Item (2) implies item (1), but we will show how to obtain (1) directly from Proposition~\ref{p:BSGT}(i).

We also remark here that item (2) provides a proof of \cite[Lemma~2.13]{bgth} which doesn't depend on \cite{Tao-BSG}, since weak general position (wgp, see Definition~\ref{d:wgp}) and Zariski density will be preserved in the coarse approximate subgroup $H$.
\end{remark}
\begin{proof}
	Let $(g,h)\in Z$ be such that $\bdl(g,h)=\bdl(Z)=\bdl(A)+\bdl(B)$ (by Fact~\ref{f:bdlInd}(Extension)).
	Since $\bdl(g)\leq\bdl(A)$ and $\bdl(h)\leq\bdl(B)$, we get $\bdl(g)=\bdl(A)$, $\bdl(h)=\bdl(B)$, and $g\dind h$.
	
	Note that $g\cdot h\in A\cdot^{Z}B$, thus $\bdl(g\cdot h)\leq \frac{1}{2}\bdl(g)+\frac{1}{2}\bdl(h)$. Therefore, \[\bdl(g)+\bdl(h)=\bdl(g,h)=\bdl(g\cdot h, h)\leq \bdl(g\cdot h)+\bdl(h)\leq \frac{1}{2}\bdl(g)+\frac{3}{2}\bdl(h),\] and $\bdl(g)\leq \bdl(h)$. Symmetrically, $\bdl(h)\leq\bdl(g)$.
	Thus, $\bdl(g)=\bdl(h)$ and $\bdl(g\cdot h)\leq \bdl(g)$. Since $\bdl(g\cdot h)= \bdl(g,h)-\bdl(h/g\cdot h)\geq \bdl(g)+\bdl(h)-\bdl(h)=\bdl(g),$ we get $\bdl(g\cdot h)=\bdl(g)$. Thus, $g \dind  h \dind  g\cdot h \dind  g$.

        Now let $D$ be the countable locally $\bdl$-finite set, $\bdl$-independent with $gh$, given by
        Proposition~\ref{p:BSGT}. Add $D$ to the constants.

        We first show (1).
        By Proposition~\ref{p:BSGT}(i), for any $g'\equiv g$ and 
        $h'\equiv h$ we have $\bdl(g'\cdot h')\leq\bdl(g)$.
	Let $s:=\tp(g)$ and $r:=\tp(h)$.
	By Lemma~\ref{l:bdlf(p)}, \[\bdl(s(M)\cdot r(M))\leq\bdl(g).\]
    Let $A'\subseteq s(M)$ be an internal subset with 
    $\bdl(A')=\bdl(s)=\bdl(A)$ and $B'\subseteq r(M)$ with 
    $\bdl(B')=\bdl(r)=\bdl(B)$ (which exist by Fact~\ref{f:existsInternal}).
    Then we get \[\bdl(g)\leq\bdl(A'\cdot B')\leq\bdl(s(M)\cdot 
    r(M))\leq\bdl(g)=\frac{1}{2}\bdl(A)+\frac{1}{2}\bdl(B).\]

    Now we prove (2). Let $g'\equiv g$ with $g'\dind gh$, and set $t:=\tp(g^{-1}\cdot g')$. Then by Proposition~\ref{p:BSGT}(ii), $\bdl(h_0\cdot h_1\cdot h_2)\leq\bdl(g)$ for any $h_i\models t$ with $0\leq i\leq 2$. Therefore, by Lemma~\ref{l:bdlf(p)}, \[\bdl(t(G)\cdot t(G)\cdot t(G))\leq \bdl(g).\]  Clearly, we also have $\bdl(t)\geq\bdl(g^{-1}\cdot g'/g')=\bdl(g)$.  Note that $\tp(g'^{-1}\cdot g/g') \vDash x \in g'^{-1}\cdot A$. Let $h':=(g^{-1}\cdot g')^{-1}\cdot h$. Then $B\cdot h'^{-1}\in \tp(g^{-1}\cdot g'/h')$. Since $\bdl(g,g',h)=3\bdl(g)=\bdl(g^{-1}\cdot g',h',g')$ by $\bdl$-independence, $\bdl(h',g')=\bdl(h'/g')+\bdl(g')=\bdl(g\cdot h/g')+\bdl(g')=2\bdl(g)$ and $\bdl(g^{-1}\cdot g'/g',h')= \bdl(g^{-1}\cdot g',g',h')-\bdl(g',h')\geq\bdl(g)$.
    Let $t':=\tp(g^{-1}\cdot g'/g',h')$, and let $S\subseteq t'(G)\subseteq t(G)$ be an internal set with $\bdl(S)=\bdl(t')=\bdl(g)$ (which exists by Fact~\ref{f:existsInternal}), and set $H':=S\cup\{1\}\cup S^{-1}$ and $H:=H'\cdot H'\cdot H'$. Then $H$ is a coarse approximate subgroup with $\bdl(H)=\bdl(g)$
    by Fact~\ref{f:approx}.

        Since $\tp(g'^{-1}\cdot g/g') \vDash x \in g'^{-1}\cdot A$ and $\tp(g^{-1}\cdot g'/h') \vDash x \in B\cdot h'^{-1}$, we have $A\supseteq g'\cdot (t'(G))^{-1}$ and $B\supseteq t'(G)\cdot h'$.
            Therefore, $(g'\cdot H)\cap A\supseteq (g'\cdot S^{-1})\cap A=g'\cdot S^{-1}$ and $(H\cdot h')\cap B\supseteq (S\cdot h')\cap B=S\cdot h'$. We conclude:
                \[\bdl((g'\cdot H)\cap A)\geq \bdl(g'\cdot S^{-1})=\bdl(g)=\bdl(H)\text{ and similarly }\bdl((H\cdot h')\cap B)\geq \bdl(H).\qedhere\]
\end{proof}


\section{Adequate structures}
\label{s:adeq}
This technical section exists to justify some limit constructions we use
later. In proving Theorem~\ref{t:cohActAb}, we will construct a sequence of
tuples which satisfy $\dind$-independence relations up to an error which tends
to zero, and Lemma~\ref{l:adequacy}(\ref{adequacy-limit}) will allow us to
extract a ``limit'' tuple which satisfies the independence relations
precisely. This involves expanding the language, and the work of this section
is to provide a setup in which this can be done harmlessly.

Let $M=\prod_{i\to\mathcal{U}}M_i$ be an ultraproduct of sets $M_i$ over some 
non-principal ultrafilter $\mathcal{U}$ on $\N$. If $\L$ is a first-order 
language, we write $M^\L$ to denote an $\L$-structure on $M$. Call a structure 
$M^\L$ \emph{internal} if every $\L$-definable set is internal; equivalently, 
$M^\L = \prod_{i\to\mathcal{U}}M_i^\L$ for some $\L$-structures $M_i^\L$ on 
$M_i$.

\begin{lemma} \label{l:adequacy}
  There exists a class of internal structures on $M$ in countable languages,
  which we call \defn{adequate} structures,
  such that:
  \begin{itemize}
    \item If $M^{\L_0}$ is an internal structure in a countable language 
      $\L_0$, then it admits an adequate expansion $M^{\L}$ for some $\L \supseteq 
      \L_0$.
    \item If $M^\L$ is an adequate structure and $A_0$ is an infinite 
      $\L$-definable set, then the following hold with $\bdl:=\bdl_{|A_0|}$:
  \begin{enumerate}[(a)]
  \item \label{adequacy-continuous}
    $\bdl^\L$ is continuous.
  \item \label{adequacy-constants}
    If $M^\L$ is adequate and $C \subseteq  M$ is countable, then the 
    expansion by constants $M^{\L(C)}$ is adequate.
  \item \label{adequacy-limit}
    Let $n \in \N$,
    let $a_i \in M^n$ for $i \in \omega$,
    and let $\U'$ be an ultrafilter on $\omega$.

    Then there exist $a \in M^n$, a countable language $\L' \supseteq \L$, and 
    an adequate expansion $M^{\L'}$ of $M^{\L}$, such that for any 
    $\L$-definable function $h : M^n \to M^{<\omega}$:
    \begin{enumerate}[(i)]
    \item
    $\tp^{\L}(a)=\lim_{i\to\U'}\tp^{\L}(a_i)$, 
    \item
    $\bdl^{\L'}(h(a)) = \lim_{i \rightarrow  \U'} \bdl^\L(h(a_i))$,
    \item Suppose $(Z_b)_{b \in B}$ is an $\L$-definable family,
    and for $\U'$-many $i$, we have $\bdl(\tp^\L(h(a_i)) \cup \{ x \in Z_b \}) 
    \leq  \alpha \in \R$ for all $b \in B$.
    Then $\bdl(\tp^{\L'}(h(a)) \cup \{ x \in Z_b \}) \leq  \alpha$ for all $b \in B$.
    \item The constants of $\L'$ are just the constants of $\L$.
    \end{enumerate}
  \end{enumerate}
  \end{itemize}
\end{lemma}
\begin{remark}
  To indicate the purpose of this notion of adequacy, consider trying to obtain (\ref{adequacy-limit})(ii) without passing to an expansion. If $\bdl(a_i)\leq\alpha$ for all $i$, then we should also have $\bdl(a)\leq\alpha$. But the internal sets $X_i\ni a_i$ witnessing $\bdl(a_i)\leq\alpha$ might vary with $i$, so simply taking $a$ to realise $\lim_{i\to\U'}\tp^{\L}(a_i)$ might not be enough. The proof will go via adding limits of such $X_i$ to the language so as to obtain the upper bound. Note that no one countable expansion will suffice for all such limit constructions.

  We avoid adding new constant symbols in (\ref{adequacy-limit}) because 
  constants will play a special role in our formalism -- see 
  Assumption~\ref{a:aclbase} and the definition of $\acl^0$.
\end{remark}
\begin{proof}
\newcommand{\Int}{{\operatorname{Int}}}
  Given a set $M_0$, we define the following first-order structure in a 
  countable language $\L_\calS$.
  The sorts are: $M_0$,
  and for each $m \in \N$ a sort $P_m := \pow(M_0^m)$ for the set of all 
  subsets of $M_0^m$. The language consists of $(m+1)$-ary membership 
  relations $\in_m$ on $M_0^m\times P_m$ for each $m \in \N$. We use $\in$ to 
  refer to any of the $\in_m$, treated as a binary relation in the natural 
  way.

  Write $\calS_i$ for the resulting structure obtained from $M_i$. 
  Consider the ultraproduct $\calS:=\prod_{i\to\U}\calS_i$.

  Let $\Int := \bigcup_m\prod_{i\to\U}\pow(M_i^m) = \bigcup_m P_m(\calS)$, the 
  set of internal subsets of powers of $M$. Consider $\Int$ with the structure 
  induced from $\calS$, i.e. such that the definable sets are those defined by  
  $\L_\calS$-formulas with no free variable of sort $M$. (In fact $\Int$ is 
  bi-definable with $\calS$, interpreting $M$ by the singleton sets in $P_1$, 
  but we will not use this.)

  Say $L \subseteq  \Int$ is \emph{adequate} if $L$ is a countable elementary 
  substructure of $\Int$.

  We define an internal structure $M^\L$ to be \emph{adequate} if $\L$ is 
  countable and the set of $\L$-definable sets is adequate.

  We associate to an adequate $L \subseteq  \Int$ an adequate structure 
  $M^{\L_L}$, with a relation for each $R \in L$ interpreted tautologically, 
  i.e.\ $M^{\L_L} \vDash  R(a) \Leftrightarrow  \calS \vDash  a \in R$.
  Since boolean combinations and existential quantification can be expressed 
  through definable functions in $\Int$, every $\L_L$-definable set is indeed 
  defined by a basic relation $R \in L$.

  We conclude by verifying the properties in the statement.

  Let $M^{\L_0}$ be an internal structure in a countable language; we show 
  that it has an adequate expansion.
  Let $L_0 \subseteq  \Int$ be the set of $\L_0$-definable sets.
  By Löwenheim-Skolem there exists an adequate $L$ with $L_0 \subseteq L$, and 
  then $M^{\L_L}$ is as required.

  Now suppose $M^\L$ is adequate, and let $L$ be the (adequate) set of all 
  $\L$-definable sets. Fix a pseudofinite set $A_0 \in L$ and set $\bdl:=\bdl_{A_0}$.
  First observe:
    \begin{claim} \label{c:cardComp}
    Given $q \in \Q_{\geq 0}$ and $m \in \N$, there is an 
    $\L_\calS(A_0)$-formula $\theta_q(Z)$ with free variable of sort $P_m$ such that $\Int \vDash \theta(Z) \Leftrightarrow |Z| \leq |A_0|^q$.
    \end{claim}
    \begin{proof}
      If $q = \frac kn$, we may take a formula expressing the existence
      of an injective function $Z^n \rightarrow A_0^k$
      (as an element of $P_{nm+km'}$, where $A_0 \in P_{m'}$).
      \subqed{c:cardComp}
    \end{proof}
  We verify the claims (a)-(c).
  \begin{enumerate}[(a)]
  \item 
    If $R \in L$, then by Claim~\ref{c:cardComp} so is the internal set
    $\{y : \exists^{\leq |A_0|^q} x\. (x,y) \in R\} = \{y : \theta_q(R_y)\}$,
    since this is the value at $R$ of an $\L_\calS(A_0)$-definable function.
    So $\L$ is closed under such cardinality quantifiers,
    and hence $\bdl^\L$ is continuous (as in \cite[\S2.1.6]{BB-cohMod}).

  \item 
    Let $L' \subseteq \Int$ be the set of $\L(C)$-definable sets.
    We show $L' \preceq \Int$.
    So let $Z \subseteq P_n$ be a non-empty $\L_\calS(L')$-definable set; we 
    aim to show $Z \cap L' \neq \emptyset$.
    We can write the parameters from $L'$ used to define $Z$ as 
    $R_1(\c),...,R_k(\c)$ where $\c \in C^m$ and $R_i \in L$.
    So $Z=Z'(\c)$ where $Z' \subseteq P_n\times M^m$ is $\L_\calS(L)$ 
    definable.
    Let $\pi(Z') \subseteq M^m$ be the projection.
    Then $Y := \{ X \in P_{n+m} : \forall\x \in \pi(Z')\; X(\x) \in Z'(\x) \} 
    \subseteq P_{n+m}$ is also $\L_\calS(L)$ definable and non-empty,
    and $L \preceq \Int$,
    so say $X_0 \in Y \cap L$.
    Then $X_0(\c) \in Z \cap L'$, as required.

  \item 
    Let $a_i \in M^n$ and $\U'$ be as in the statement.

\newcommand{\Lt}{\vec L}
    Given $a\in M^n$ and $L \subseteq L' \subseteq \Int$, write $(L',\Lt,a)$ 
    for the structure with universe $L'$ and basic predicates for the trace on 
    $L'$ of $\L_\calS(\Lt,a)$-formulas with all free variables of sorts from 
    $\Int$, where $\Lt$ denotes a tuple enumerating $L$.

  \begin{claim}\label{c:limitType}
    There exists a countable $L' \preceq \Int$ with $L \subseteq L'$
    and $a \in M^n$
    such that $(L',\Lt,a) \vDash \lim_{i \to \U'} \Th(L,\Lt,a_i)$.
  \end{claim}
  \begin{proof}
    By Löwenheim-Skolem, say $\calN$ is a countable model of this theory.
    By $\aleph_1$-saturation of $\calS$, let $a \in M^n$ realise $\lim_{i \to 
    \U'} \tp^{\L_\calS}(a_i/\Lt)$.
    Then by $\aleph_1$-saturation again, the diagram of $\calN$ is realised by 
    some $L' \subseteq \Int$ such that $\calN \cong (L',\Lt,a)$, since 
    consistency of finite parts of the diagram is ensured by existential 
    statements in $\tp^{\L_\calS}(a/\Lt)$.
    Finally, $L' \preceq \Int$ since in particular $L' \vDash 
    \Th^{\L_\calS}(L) = \Th(\Int)$, and $L'\subseteq \Int$ is a substructure 
    in a language with predicates for $\L_\calS$-formulas.
  \end{proof}

  Consider the corresponding adequate expansion $M^{\L'}$ of $M^\L$.
  Now condition (iv) from the statement holds by construction, and we conclude 
  by observing that the language of $(L',\Lt,a)$ is expressive enough for 
  (i)-(iii) to follow from Claim~\ref{c:limitType}:

  \begin{enumerate}[(i)]
  \item
    If $X \in L$, then $a \in X$ is expressed by a sentence in the language of 
    $(L',\Lt,a)$.\footnote{This is immediate if when defining $(L',\Lt,a)$ we 
    include 0-ary relation symbols interpreted as $\L_\calS(\Lt,a)$-sentences; 
    if instead we use a logical formalism which does not include 0-ary 
    relation symbols, we can make do with sentences of the form $\exists Y (Y 
    = Y \land a \in X)$.}
  \item
    $\bdl^{\L'}(h(a)) \leq \alpha$ is expressed in the 
    language of $(L',\Lt,a)$ via Claim~\ref{c:cardComp} by $\{ \exists R \in L'\; 
    (h(a) \in R \land \theta_q(R)) : q \in \Q_{>\alpha} \}$.
  \item
    Let $q > \alpha$.
    By $\aleph_1$-saturation of $\calS$,
    for $\U'$-many $i$ we have $\exists R \in L\; (h(a_i) \in R \land \forall b 
    \in B\. \theta_q(R \cap Z_b))$;
    indeed, otherwise we obtain (for this $i$) finite consistency of the countable 
    $\calS$-type in $y$
    $$\{ \neg\theta_q(R' \cap Z_y) : R' \in L,\; h(a_i) \in R' \} \cup \{y \in 
    B\},$$
    and a realisation $b \in B$ contradicts the assumption, since $q>\alpha$.
    So the corresponding statement holds for $(L',\Lt,a)$, and we conclude
    since $q>\alpha$ was arbitrary.
    \qedhere
  \end{enumerate}
  \end{enumerate}
\end{proof}


\section{Non-expansion in homogeneous spaces}
\label{s:cohHom}
In this section, we prove as Theorem~\ref{t:cohActAb} a variant of 
\cite[Theorem~A.4]{cubicSurfaces} on nilpotence of ``non-expanding'' homogeneous 
spaces, strengthening the general position assumption on the space (from weak 
to coarse general position) and strengthening the conclusion to abelianity.

\subsection{Further setup and preliminaries}
\label{ss:prelims2}
From this section on, we work in a non-principal ultrapower $K:=\C^{\U}$ of
the complex numbers.
We consider $K$ as an adequate expansion of the field structure in a language
$\L\supseteq \Lring$,
with $\bdl=\bdl_{|A_0|}$ for some infinite $\L$-definable set $A_0$. We also 
assume:
\begin{assumption}\label{a:aclbase}
  The set $E_0$ of interpretations of the constants in $\L$ is an algebraically closed subfield of $K$.
\end{assumption}

Unless otherwise specified, all sets of parameters are \emph{countable} 
subsets of $K$, and all tuples are finite tuples from $K$.

As in \cite[Definition~2.4]{cubicSurfaces}, we use superscript-0 to refer to
the reduct of $K$ to the ring language expanded by the constants of $\L$. In particular, $\acl^0(C)$ is the smallest algebraically closed subfield of $K$ containing $CE_0$, and $\dim^0(a/C)$ denotes the transcendence degree of $\acl^0(Ca)$ over $\acl^0(C)$, while $a\ind^0_C b$ means $\dim^0(a/Cb) = \dim^0(a/C)$, and $\loc^0(a/C)$ is the smallest algebraic subvariety over $CE_0$ containing $a$.
We write $a \sim^0_C b$ to mean $\acl^0(aC) = \acl^0(bC)$.

Our emphasis on expansions of the ring language by constants is motivated by the following situation. Suppose a pseudofinite set \(X\), which is Zariski dense in an algebraic group \(G\), is nevertheless concentrated on a coset \(cH\) of an algebraic subgroup \(H \leq G\), in the sense that \(\bdl(X \cap cH) = \bdl(X)\). In this case, it is important to isolate the relevant coset. To do so, we expand the language by naming the element \(c\) as a constant and work with the set \(X \cap cH\) instead of \(X\). This reduction will be implemented by passing to types in weak general position after adding constants (see Definition~\ref{d:wgp} and Lemma~\ref{l:wgpification}).

If $V$ is an algebraic variety over a subfield $k \leq K$,
by a \defn{generic} element over $k$ of $V$ we
mean an element $a \in V(K)$ not in any lower-dimensional subvariety of $V$ defined
over $\acl^0(k)$. By a \emph{generic element of $V$} we mean a generic element
over some field $k$ over which $V$ is defined. Unless otherwise stated, by
a \emph{subvariety} of a variety $V$, we mean a closed subvariety.

We will use a little of the model theoretic terminology of canonical bases.
See \cite[\S2.1.12]{BB-cohMod} for a relevant exposition; briefly, we use 
$\cb^0(a/C)$ to refer to the $\ACF$-canonical base of $\tp^0(a/\acl^0(C))$. That is, the field generated over the constants of $\L$ by $\cb^0(a/C)$ coincides with the field
generated over the constants by the field of definition of $\loc^0(a/\acl^0(C))$.
We often use this in the context of the following fact, whose proof is a straightforward exercise.
\begin{fact}\label{f:cb}
  Suppose $a \sim^0_d b$ and $a \ind^0 d \ind^0 b$.
  Then $\acl^0(\cb^0(ab/d))$ is the smallest $\acl^0$-closed subset of 
  $\acl^0(d)$ over which $a$ and $b$ are $\acl^0$-interalgebraic:
  $$\acl^0(\cb^0(ab/d)) = \bigcap \{ \acl^0(d') : d' \in \acl^0(d),\; a 
  \sim^0_{d'} b \}.$$
\end{fact}

We also use the following standard result.
\begin{fact}[{\cite[Lemma~1.2.28]{GST}}] \label{f:cbMS}
  If $C = \acl^0(C)$, then $\cb^0(a/C) \in \acl^0((a_i)_{i\in\omega})$ (and even $\dcl^0((a_i)_{i<\omega})$) whenever $(a_i)_{i\in\omega}$ is a 
  \emph{Morley sequence} in $\tp^0(a/C)$,
  meaning $a_i \equiv^0_C a$ and $a_{i+1} \ind^0_C a_0\ldots a_i$.
\end{fact}

\subsection{General position}
We use the following definitions of general position:
\begin{definition}\label{d:wgp}
    For a (countable) set $C \subseteq K$ and a tuple $a \in K$,
  \begin{itemize}
    \item $\tp(a/C)$ is \defn{in weak general position} (\defn{wgp}) if $\bdl(a/C)<\infty$ and for any tuple $b \in K$ such that $a \nind^0_C b$, we have $\bdl(a/Cb) < \bdl(a/C)$.
    \item $\tp(a/C)$ is \defn{in coarse general position} (\defn{cgp}) if $\bdl(a/C)<\infty$ and for any tuple $b \in K$ such that $a \nind^0_C b$, we have $\bdl(a/Cb) < \bdl(a/C)$ and $\bdl(a/Cb) = 0$.
  \end{itemize}
\end{definition}
\begin{remark}\label{r:cgpWgp}
  It follows immediately from these definitions that any cgp type is also wgp.
  To ensure this holds, our definition of cgp is slightly different from that in \cite[Definition~5.6]{BB-cohMod}, but the definitions agree when $\bdl(a/C)>0$.
\end{remark}
\begin{remark}\label{r:wgpAlg}
 If $\tp(a/C)$ is wgp and $\bdl(a/C) = 0$, then no $b$ as in the definition exists, so $a \in \acl^0(C)$.
\end{remark}
\begin{remark}\label{r:dim1cgp}
  If $\tp(a/C)$ is broad and $\dim^0(a/C) = 1$, for example if $a$ is an element of $K$, then $\tp(a/C)$ is cgp.
  Indeed, if $a \nind^0_C b$ then $a \in \acl^0(Cb)$, so $\bdl(a/Cb) = 0 < \bdl(a/C)$.
\end{remark}

Part (i) of the following lemma, which is essentially a reformulation of 
\cite[Lemma~2.13]{Hr-psfDims}, shows that weak general position is a harmless 
assumption from our coarse perspective, in the sense that we can always 
ensure it by working over parameters which do not affect $\bdl$.

\begin{lemma} \label{l:wgpification}
  Let $\tp(a/C)$ be a broad type.
  \begin{enumerate}[(i)]
    \item There exists $d$ such that $\bdl(d/C) = 0$ and $\tp(a/Cd)$ is wgp, in particular $\bdl(a/Cd)=\bdl(a/C)$.
    \item There exists $d$ such that $\bdl(d/C) < \infty$ and $\tp(a/Cd)$ is cgp and broad.
  \end{enumerate}
\end{lemma}
\begin{proof}


  \begin{enumerate}[(i)]
  \item Let $b\in K$ be a finite tuple such that $\bdl(a/Cb) = \bdl(a/C)$ and $\dim(a/Cb)$ achieves the minimal value for such $b$.
  Then $\tp(a/Cb)$ is wgp;
  indeed, if $a\nind^0_{Cb}b'$ for some tuple $b'\in K$, then by the minimality of $\dim(a/Cb)$, we must have $\bdl(a/Cbb')<\bdl(a/Cb)$.
  Let $d:= \cb^0(a/Cb)$.
  Let $(a_i)_{i\in\omega}$ be a $\bdl$-independent sequence of realisations of $\tp(a/\acl^0(Cb))$.
  Since $\tp(a/Cb)$ is wgp, $(a_i)_{i\in\omega}$ is a Morley sequence in $\tp^0(a/\acl^0(Cb))$.
  Then by Fact~\ref{f:cbMS}, $d\in \acl^0(Cb) \cap \acl^0(\a)$ where 
  $\a$ is a finite subsequence of $(a_i)_{i\in\omega}$.
  But then $\a \dind_C b$ and hence $\bdl(d/C)=0$; we can see this explicitly as follows.
  Note that \[\bdl(\a/C)\geq\bdl(\a/\acl^0(Cb))=|\a|\bdl(a/\acl^0(Cb))=|\a|\bdl(a/C)\geq\bdl(\a/C).\]
  Hence, $\bdl(\a/Cd)=\bdl(\a/C)$ (by $d\in\acl^0(Cb)$) and $\bdl(d/C) = \bdl(\a d/C) - \bdl(\a/Cd)= \bdl(\a /C) - \bdl(\a/Cd)= 0$, where $\bdl(\a d/C)=\bdl(\a/C)$ by $d\in\acl^0(\a)$.

  \item Let $b\in K$ be a finite tuple such that $\bdl(a/Cb) > 0$ and $\dim(a/Cb)$ achieves the minimal value for such $b$. Such a $b$ exists since $\bdl(a/C) > 0$.
  Then $\tp(a/Cb)$ is cgp;
 indeed, if $a\nind^0_{Cb}b'$ for some tuple $b'\in K$, then by the minimality of $\dim(a/Cb)$, we must have $\bdl(a/Cbb')=0$.
 Let $d:= \cb^0(a/Cb)$.
 Then by the same argument as before, $d\in\acl^0(\a)$, where $\a$ is a 
 $\bdl$-independent sequence of realisations of $\tp(a/\acl^0(Cb))$. Hence, 
 $\bdl(d/C)<\bdl(\a/C)\leq |\a|\bdl(a/C)<\infty$.
  \end{enumerate}
\end{proof}

\begin{fact}[{\cite[Lemma~2.13]{cubicSurfaces}}]\label{f:wgpFacts}
  Suppose $\bdl(a/C),\bdl(b/C) < \infty$.
  \begin{enumerate}[(i)]
    \item If $\tp(a/C)$ is wgp and $d \in \acl^0(aC)$, then $\tp(d/C)$ is wgp.
    \item If $\tp(b/C)$ and $\tp(a/Cb)$ are wgp, then so is $\tp(ab/C)$.
    \item If $a \ind ^{\bdl}_C b$ and $\tp(a/C)$ is wgp, then $\tp(a/Cb)$ is wgp.
  \end{enumerate}
\end{fact}

\begin{remark}\label{r:tuplewgp}
 From (ii) and (iii), it is immediate that if $a \ind ^{\bdl}_C b$ and both $\tp(a/C)$ and $\tp(b/C)$ are wgp, then $\tp(ab/C)$ is wgp (however, this does not hold for cgp).
\end{remark}
\begin{remark}\label{r:aclbroad}
		If $\tp(a/C)$ is wgp and $d \in \acl^0(aC)\setminus \acl^0(C)$, then $\tp(d/C)$ is broad;
		this can be seen from (i) and Remark~\ref{r:wgpAlg}.
\end{remark}

\begin{lemma}\label{l:cgppreserve}
Suppose $\bdl(a/C),\bdl(b/C) < \infty$.
\begin{enumerate}[(i)]
    \item  If $\tp(a/C)$ is cgp and $d\in\acl^0(aC)$, then $\tp(d/C)$ is cgp.
    \item  If $a \ind ^{\bdl}_C b$ and $\tp(a/C)$ is cgp, then $\tp(a/Cb)$ is 
      cgp.
    \item  If $a \ind ^{\bdl}_C b$ and $\tp(b/C)$ is wgp and $\tp(a/Cb)$ is 
      cgp, then $\tp(a/C)$ is cgp.
  \end{enumerate}
\end{lemma}
\begin{proof}
  \begin{enumerate}[(i)]
  \item
  If $d\in\acl^0(C)$, then $\tp(d/C)$ is clearly cgp.
  Otherwise, by Remark~\ref{r:aclbroad}, $\tp(d/C)$ is broad.
  Suppose $d\nind^0_{C} e$ for some $e$.
  Then $a\nind^0_{C} e$ (as $d\in\acl^0(aC))$.
  Thus, $\bdl(a/Ce)=0$ by $\tp(a/C)$ being cgp,
  and so $\bdl(d/Ce)\leq \bdl(a/Ce)=0<\bdl(d/C)$.

  \item
  Suppose $\tp(a/C)$ is cgp and $a \nind^0_{Cb} d$ for some $d$.
  Then $a \nind^0_C bd$, and hence
  $\bdl(a/Cbd)=0<\bdl(a/C)=\bdl(a/Cb),$
  by $\tp(a/C)$ being cgp and $a\ind^{\bdl}_C b$.
  Therefore, $\tp(a/Cb)$ is cgp.

  \item
  Suppose $a \nind^0_C d$ for some $d$.
  Replacing $d$ with a canonical base as in the proof of Lemma~\ref{l:wgpification},
  we may assume $\bdl(d/C) < \infty$.
  Choose $d'\equiv_{Ca} d$ such that
  $d' \ind^{\bdl}_{Ca} b$ (such a $d'$ exists by extension).
  Since $a \ind^{\bdl}_C b$, by transitivity we obtain $d'a \ind^{\bdl}_C b$, and hence
  $a \ind^{\bdl}_{Cd'} b$.

  Now $a \ind^0_C b$ since $\tp(b/C)$ is wgp,
  and $a \nind^0_C d'$ since $ad \equiv_C ad'$ and $a \nind^0_C d$,
  so $a \nind^0_{Cb} d'$.
  Thus, $0=\bdl(a/Cbd')<\bdl(a/Cb)=\bdl(a/C)$,
  by the assumption that $\tp(a/Cb)$ is cgp.
  Using $a \ind^{\bdl}_{Cd'} b$, we conclude that $0=\bdl(a/Cd)=\bdl(a/Cd')=\bdl(a/Cd'b)<\bdl(a/C)$.
  Hence $\tp(a/C)$ is cgp, as desired.
  \end{enumerate}
\end{proof}

\begin{lemma}\label{l:gpSeq}
  Suppose $\a = (a_1,\ldots,a_k)$ is a finite sequence such that $\bdl(a_i) = \bdl(a_j) < \infty$ and $\dim^0(a_i)=\dim^0(a_j)$ for all $i,j$, and $a_{i+1} \dind a_1\ldots a_i$ for all $i$.
  Let $b \in \acl^0(\a)$. Then:
  \begin{enumerate}[(i)]
    \item If $\tp(a_i)$ is wgp for all $i$, then $\tp(b)$ is wgp.
    \item If $\tp(a_i)$ is cgp for all $i$, then
      $\bdl(\a)\cdot\dim^0(b) \leq \bdl(b)\cdot\dim^0(\a)$.
  \end{enumerate}
\end{lemma}
\begin{proof}
  \begin{enumerate}[(i)]
    \item This follows by induction on $k$ from Fact~\ref{f:wgpFacts}.
    \item This is immediate if $\bdl(a_i) = 0$. Otherwise, rescaling, we may assume $\bdl(a_i) = \dim^0(a_i)$.
      The statement is then an instance of \cite[Lemma~5.13]{BB-cohMod}.
      \qedhere
  \end{enumerate}
\end{proof}

\subsection{Algebraic homogeneous spaces}
\label{ss:CFAHS}
\begin{definition}
	An \defn{algebraic homogeneous space} is a pair $(G,X)$, where $G$ is an algebraic group and $X$ an algebraic variety, along with a transitive action of $G$ on $X$ such that the 
  corresponding map $* : G\times X \rightarrow  X$ is a regular map.
  We say $(G,X)$ is a \defn{faithful algebraic homogeneous space} if in addition 
  the action is faithful,
  and we say $(G,X)$ is \defn{connected} if $G$ is connected.
  We abbreviate ``connected faithful algebraic homogeneous space'' to \defn{CFAHS}.
  We write $g\cdot h$ for the group operation on $G$, and $g*x$ for the action.
  We say $(G,X)$ is \defn{defined over} a field $C$ if the algebraic group
  $G$, the variety $X$, and the action $*$ are all defined in the algebraic
  sense over $C$ (i.e.\ by $\Lring(C)$-formulas).
\end{definition}

\begin{lemma} \label{l:trivStabZ}
  Let $(G,X)$ be a faithful transitive (abstract) group action, let $a \in X$,
  and let $G_a$ be the stabiliser of $a$ in $G$. Then $G_a\cap Z(G)=\{1\}$.
\end{lemma}
\begin{proof}
  Let $g\in G_a\cap Z(G)$. By transitivity of the action, 
  for any $b\in X$, there is $h\in G$ with $b=h*a$. Then 
  $g*b=g\cdot h*a=h\cdot g*a=h*a=b$. Hence, $g=1$ by faithfulness.
\end{proof}

\begin{lemma} \label{l:smallCentreOrbit}
  Let $(G,X)$ be a CFAHS. Suppose $G$ is nilpotent but non-abelian.
  Then no orbit of the centre is generic in $X$,
  i.e.\ $\dim(Z(G)*a) < \dim(X)$ for all $a\in X$.
\end{lemma}
\begin{proof}
  Let $a \in X$.
  By Lemma~\ref{l:trivStabZ}, $G_a\cap Z(G) = \{1\}$, and $Z(G)$ acts
  regularly on $Z(G)*a$, and so $\dim(Z(G))=\dim(Z(G)*a)$.

  Considering orbits and stabilisers, we have
  $\dim(G)=\dim(G*a)+\dim(G_a)=\dim(X)+\dim(G_a)$.
  Suppose towards a contradiction that $\dim(Z(G)*a) = \dim(X)$,
  so $\dim(G)=\dim(Z(G))+\dim(G_a)$.
  Consider the map $\eta: Z(G)\times G_a\to G$ sending $(g,h)\mapsto g\cdot h$. It 
  is a group homomorphism with trivial kernel, so its image has dimension 
  $\dim(Z(G))+\dim(G_a)=\dim(G)$. Thus, $\eta$ is a group isomorphism by 
  connectedness of $G$. Now $G_a$ is non-trivial since $G\neq Z(G)$, and $G_a$ 
  is nilpotent since $G$ is, so $Z(G_a)$ is non-trivial. However, we must have 
  $Z(G_a)\leq Z(G)$, as $G\cong Z(G)\times G_a$. This contradicts the fact 
  that $G_a\cap Z(G)$ is trivial.
\end{proof}

Recall that the \emph{nilpotency class} of a nilpotent group $G$ is the least $k$ such that
$[g_0,\ldots ,g_k] := [\ldots [[g_0,g_1],g_2],\ldots g_k] = 1$ for all $g_0,\ldots,g_k \in G$.
\begin{lemma} \label{l:itCommCentral}
  Let $G$ be a connected nilpotent algebraic group over a 
  field $K_0$, and say $G$ has nilpotency class $k$.
  Let $(g_0,\ldots ,g_{k-1}) \in G^k$ be generic over $K_0$,
  and let $t := [g_0,\ldots ,g_{k-1}]$.

  Then $t$ is central in $G$, and $\trd(t/K_0) > 0$.
\end{lemma}
\begin{proof}
  First, $t$ is central as $[g_0,\ldots,g_{k-1},g_{k}]=1$ for any $g_{k} \in G$.
  Suppose for a contradiction that $t \in \acl(K_0)$. Then $[h_0,\ldots,h_{k-1}] = t$ for all $h_i \in G$ by the genericity (since $G$ is connected and $\{(x_0,\cdots,x_{k-1}):[x_0,\ldots,x_{k-1}] = t\}$ is a Zariski-closed set in $G^k$), so $t=1$ and we contradict $G$ having nilpotency class $k$.
\end{proof}

\subsection{Abelianity}
\begin{theorem} \label{t:cohActAb}
  Let $(G,X)$ be a CFAHS defined over $\acl^0(\emptyset)$.

  Let $g \in G$ be generic, let $a \in X$ be generic,
  and suppose $\tp(g)$ is wgp and $\tp(a)$ is cgp,
  and $a \dind g \dind g*a$.

  Then $G$ is abelian, and so $(G,X)$ is isomorphic over $\acl^0(\emptyset)$ 
  with the homogeneous space $(G,G)$ in which $G$ acts on itself by addition.
\end{theorem}

Let us summarise the proof strategy before going into the details of the proof. Consider first the special case that the group action is $G$ acting on itself, and suppose moreover that in addition to $h \dind g \dind g\cdot h$, we have $g\cdot h\dind h$, or equivalently $\bdl(g)=\bdl(h)$. Then we obtain a (coarse) approximate subgroup by Balog-Szemerédi-Gowers, and it was observed in \cite[Proposition~2.21]{bgth} (and implicitly already in \cite{Breuillard-Wang}) that $\tp(h)$ being cgp then forces $G$ to be abelian.
To get from this special case to the generality of our statement, there are two obstacles. First, we may have a non-trivial group action. If the stabilisers are trivial, we are essentially in the above case of the group acting on itself. Since we work with algebraic group actions, we can easily reduce to the case of generically trivial stabilisers by considering the diagonal action of $G$ on $X^n$ for sufficiently large $n$. (In fact our arguments also go through for arbitrary group actions with suitable conditions on the stabilisers, and we formulate a result of this form in Appendix~\ref{appB}.) Secondly, we may not be in the symmetric situation that $\bdl(g)=\bdl(h)$. We handle this by iterating, essentially as in \cite[Theorem~A.4]{cubicSurfaces}. Replacing $g$ by $g_1 := g'^{-1}\cdot g$ where $g'$ is an independent copy of $g$ over $g\cdot h$, we obtain a new triple $(g_1, h, g_1\cdot h)$. We may have $\bdl(g_1)>\bdl(g)$, but $h \dind g_1 \dind g_1\cdot h$ and so $\bdl(g_1)\leq\bdl(h)$. Iterating this process, we obtain $g_i$ with $\bdl(g_i)\leq\bdl(h)$ and such that the independence relations $g_i \dind g_i'\dind g_i'^{-1}\cdot g_i=g_{i+1}$ are almost satisfied for large $i$. There may be some error for all $i$ (as in Example~\ref{e:span}), but using the limit construction in \S\ref{s:adeq} we reach exact independence, putting us in the group action BSG situation of Corollary~\ref{c:symBSGact}. It follows as in \cite{cubicSurfaces} that $G$ is nilpotent, and we then use the group theoretic results of \S\ref{ss:CFAHS} and the conclusion of Corollary~\ref{c:symBSGact} to see that if $G$ is not abelian, then $a$ concentrates on its orbit under the centre in a way which contradicts the cgp assumption.

\begin{proof}
  It suffices to show that $G$ is abelian. Indeed,
  since the action is faithful, it would then be free. So $(G,X)$ is then a 
  principal homogeneous space, and if $x_0 \in X(\acl^0(\emptyset))$,
  then $g+x_0 \mapsto g$ establishes an isomorphism $(G,X) \cong (G,G)$.

  \newcommand{\tda}{\widetilde a}
  \begin{claim} \label{c:cohActAb}
    The result holds if we also assume that there exists
    $h \in G$ such that $g \dind h \dind g\cdot h \dind g$.
  \end{claim}
  \begin{proof}
  \providecommand{\g}{\overline g}
  $G$ is nilpotent by \cite[Theorem~A.4]{cubicSurfaces} using the assumptions
  of the Theorem, or alternatively by \cite[Remark~2.12]{bgth} using the
  assumption of the Claim. Here we use that the underlying field has
  characteristic 0.

  Suppose for a contradiction that $G$ has nilpotency class $k > 1$.

  Let $D$ be as in Corollary~\ref{c:symBSGact}.
  We may assume $D \dind ga$ (by replacing $D$ with a realisation of $\tp(D/g)$ $\bdl$-independent from $a$ over $g$).
  Adding $D$ as constants to the language then preserves the assumptions of the Theorem, so we may assume $D=\emptyset$.

  Let $\g = (g_0,\ldots ,g_{2k-1})$ be a $\bdl$-independent 
  sequence in $\tp(g)$ such that $\g \dind ag$.
  Let $t := [g_0^{-1}\cdot g_1,\ldots ,g_{2k-2}^{-1}\cdot g_{2k-1}]$.
  Since $g$ is generic in $G$ and $\tp(g)$ is wgp, $\g$ is 
  generic in $G^{2k}$, so 
  $(g_0^{-1}\cdot g_1,\ldots ,g_{2k-2}^{-1}\cdot g_{2k-1})$ is 
  generic in $G^k$ (since it is the product of two independent generics $(g_0,\ldots,g_{2k-2})^{-1}$ and $(g_1,\ldots,g_{2k-1})$ in $G^k$).
  By Lemma~\ref{l:itCommCentral}, $t$ is central and $\dim^0(t) > 0$.

  Now $t \in Z(G)$, so $Z(G)*(t*a) = Z(G)*a$. Hence $\dim^0(a/t*a) \leq \dim(Z(G)*(t*a)) < \dim(X) = \dim(a)$, where the inequality holds by Lemma~\ref{l:smallCentreOrbit}, 
  and so $\bdl(a/t*a)=0$ since $\tp(a)$ is cgp.

  Meanwhile, $t$ can be written out as an alternating word in
  realisations of $\tp(g)$ with suitable adjacent elements $\dind$-independent
  as in the conclusion of Corollary~\ref{c:symBSGact},
  so $\bdl(t*a) \leq  \bdl(a)$.

  Then $0 = \bdl(a/t*a) = \bdl(a) - (\bdl(t*a) - \bdl(t*a/a)) \geq  
  \bdl(t*a/a)$.
  But $\bdl(t*a/a) = \bdl(t/a) = \bdl(t)$,
  since $t \in \dcl(t*a,a)$,
  since $G_{a} \cap Z(G) = \{1\}$ by Lemma~\ref{l:trivStabZ}.
  So $\bdl(t) = 0$; but $t$ is wgp since $t \in \acl^0(\g)$ (using Fact~\ref{f:wgpFacts}(i)),
  which contradicts $\dim^0(t)>0$ via Remark~\ref{r:wgpAlg}.

  Therefore, $G$ is abelian.
  \subqed{c:cohActAb}
  \end{proof}

  We now perform a limit construction to replace $(g,a)$ with a pair satisfying the
  assumption in the Claim as well as those in the Theorem.

  As in the proof of \cite[Theorem~A.4]{cubicSurfaces}, there is $n \in \N$ 
  such that the stabiliser of a generic point of $X^n$, under the diagonal action of $G$ on $X^n$, is trivial.
  So, extending $a$ to a $\bdl$-independent $n$-tuple 
  $\a=(a_1,\ldots ,a_n)$ with $a_1=a$ and $a_i \equiv_g a$ and $a_i \dind_g 
  a_1\ldots a_{i-1}$, we have $g \dind \a$ and $g \dind g*\a$ and $\a$ 
  has trivial stabiliser, $G_\a = 1$.

  Define a sequence $(g_i)_{i \in \omega}$ of elements in $G$ and simultaneously another sequence $(\bar{h}_i)_{i\in\omega}$ of tuples of realisations of $\tp(g)$ as follows:
  $g_0 := g=:h_0$, and given $g_i$, let $g_i'\bar{h}'_i \equiv_{g_i*\a} g_i\bar{h}_i$ with 
  $g_i'\bar{h}'_i  \dind_{g_i*\a} g_i\bar{h}_i \a$,
  and define $g_{i+1} := {g_i'}^{-1}\cdot g_i$ and $\bar{h}_{i+1}:=\bar{h}_i\bar{h}_i'$.

  \begin{claim}\label{c:gi}\footnote{This claim does not rely on the assumption that $G_\a = 1$, nor on the fact that $(G,X)$ is algebraic. Although the conclusion $\bdl(g_i/g_i\cdot g^{-1})=\bdl(g)$ will not be used in the present proof, we will apply this claim later in Appendix~\ref{appB}. There, this conclusion will be needed, and it will also be important to keep track of $\bar{h}_i$.}
 For all $i\in\omega$, $\bar{h}_i$ is a $\bdl$-independent sequence in $\tp(g)$, and $g_i\in\langle \bar{h}_i\rangle$ where $\langle \bar{h}_i\rangle$ denotes the subgroup generated by $\bar{h}_i$, and
  \[\a \dind \bar{h}_i \dind g_i*\a, \quad \bar{h}_i' \dind \bar{h}_i\a,\]
 and $\bdl(g_i/g_i\cdot g^{-1})=\bdl(g)$.
  \end{claim}
  \begin{proof}
    First note that the condition $\bar{h}_i' \dind \bar{h}_i\a$ follows from
    $\bar{h}_i \dind g_i*\a$. Indeed,
    by $\bar{h}_i' \equiv_{g_i*\a} \bar{h}_i$ we obtain $\bar{h}_i' \dind g_i*\a$. Since $\bar{h}_i' \dind_{g_i*\a} \bar{h}_i \a$, transitivity yields $\bar{h}_i' \dind \bar{h}_i \a(g_i*\a)$.

    We prove the result by induction on $i$. For $i=0$, the result is clear. Suppose for $i$. Then $g_{i+1}\in\langle \bar{h}_{i+1}\rangle$ by construction, and $\bar{h}_{i+1}$ is a $\bdl$-independent sequence in $\tp(g)$ since $\bar{h}'_i\dind \bar{h}_i$ and $\bar{h}_i$ (and hence also $\bar{h}'_i$) is a $\bdl$-independent sequence, and by transitivity of $\bdl$-independence.

  By the induction hypothesis, $\bar{h}_i' \dind_{\bar{h}_i} \a$ and $\bar{h}_i \dind \a$. By transitivity, we conclude that $\bar{h}_i' \bar{h}_i \dind \a$, and therefore $\bar{h}_{i+1} \dind \a$.

 The induction hypothesis $\bar{h}_i \dind g_i*\a$ and $\bar{h}_i' \dind_{g_i*\a} \bar{h}_i$ yields $\bar{h}_i'(g_i*\a) \dind \bar{h}_i $. Thus, $g_{i+1}*\a \dind_{\bar{h}'_i} \bar{h}_i $ (as $g_i'\in \langle\bar{h}'_i\rangle$). Note that $\bar{h}_ig_i(g_i*\a)\a\equiv \bar{h}'_ig_i'(g_i*\a)(g_{i+1}*\a)$ (as $g_i'^{-1}*(g_i*\a)=g_{i+1}*\a$). Since $\bar{h}_i\dind \a$, we have $\bar{h}_i'\dind g_{i+1}*\a$. By transitivity and $g_{i+1}*\a\dind_{\bar{h}_i'}\bar{h}_i$, we get $g_{i+1}*\a\dind \bar{h}_i'\bar{h}_i$, thus $g_{i+1}*\a\dind \bar{h}_{i+1}$ as desired.

 Finally, note that by the construction, $g$ is an element of the sequence $\bar{h}_i$ and each $g_i$ can be written as $u\cdot g$ where $u$ is a product of elements of $\bar{h}_i \setminus \{g\}$ and their inverses. Since $\bar{h}_i$ is $\bdl$-independent, we have $g \dind u$. Then
 \[\bdl(g_i/u)=\bdl(g/u)=\bdl(g).
   \subqedinmath{c:gi}
 \]
  \end{proof}


In particular, we have  for all $i\in \omega$:
\begin{equation} \label{e:gind}
	\a \dind g_i \dind g_i*\a, \quad g_i' \dind g_i\a.
\end{equation}

Then the sequence $(\bdl(g_i))_{i\in\omega}$ is non-decreasing,
  since
  $$\bdl({g_i'}^{-1}\cdot g_i) \geq \bdl({g_i'}^{-1}\cdot g_i/g_i') = \bdl(g_i/g_i') = \bdl(g_i).$$
  Using $G_\a = 1$, we have $ \bdl(g_i/\a) = \bdl(g_i*\a/\a)$.
  Thus,
  $$\bdl(g_i) = \bdl(g_i/\a) = \bdl(g_i*\a/\a) \leq \bdl(g_i*\a) = \bdl(\a).$$
  So let $\alpha := \lim_{i\to\infty}\bdl(g_i) \leq \bdl(\a) < \infty$.

  We return to the action of $G$ on $X$, rather than $X^n$, for the remainder of the proof.
  Let $(\gamma,\gamma',\tda)\in G\times G\times X$
  and an adequate expansion to $\L' \supseteq \L$ be as in 
  Lemma~\ref{l:adequacy}(\ref{adequacy-limit}) applied to $(g_i,g_i',a)_{i \in 
  \omega}$ and a non-principal ultrafilter on $\omega$. Work in $\L'$, so in
  particular we redefine $\bdl := \bdl^{\L'}$.
  By Lemma~\ref{l:adequacy}(\ref{adequacy-limit})(iv), passing to $\L'$ does
  not affect the meaning of $\acl^0$ or $\ind^0$.

  From Lemma~\ref{l:adequacy}(\ref{adequacy-limit})(ii) and (\ref{e:gind}), we obtain
  $\tda \dind \gamma \dind \gamma*\tda$ and $\gamma\dind\gamma'$.
  Indeed, $\bdl(\tda)=\bdl^\L(a)$ and $\bdl(\gamma')=\bdl(\gamma)=\lim_{i\to\infty}\bdl^\L(g_i)=\alpha$. Thus \[\bdl(\gamma \tda)=\lim_{i\to\infty}\bdl^\L(g_ia)=\lim_{i\to\infty}\bdl^\L(g_i)+\bdl^\L(a)=\bdl(\gamma)+\bdl(\tda).\] The other independences follow in the same manner.

  By \cite[Lemma~A.3]{cubicSurfaces}, for any constructible family $(V_b)_{b 
  \in B}$ over $\acl^0(\emptyset)$ of  proper subvarieties of $G$,
  there is $\eta > 0$ such that $\bdl(\tp^{\L}(g_i) \cup \{ x \in V_b \}) \leq  
  \alpha-\eta$ for all $b \in B$ and $i\in\omega$.
  Then by Lemma~\ref{l:adequacy}(\ref{adequacy-limit})(iii),
  $\gamma$ satisfies the same bounds,
  so $\tp^{\L'}(\gamma)$ is wgp.
  Now by induction, each $g_i$ is generic in $G$; indeed, $g_0=g$ is generic,
  and if $g_{i-1}$ and hence $g_{i-1}'$ are generic, then, since $g_{i-1} \ind^0 g_{i-1}'$, so is $g_i=(g'_{i-1})^{-1}\cdot g_{i-1}$.
  Hence $\gamma$ is generic in $G$ by Lemma~\ref{l:adequacy}(\ref{adequacy-limit})(i).

  Similarly, $\tp^{\L'}(\tda)$ is cgp, and $\tda$ is generic in $X$.
  So $(\gamma,\tda)$ satisfies the assumptions of this Theorem.

  Furthermore, by Lemma~\ref{l:adequacy}(\ref{adequacy-limit})(ii) we have \[\bdl(\gamma'^{-1}\cdot\gamma)=\lim_{i\to\infty}\bdl^\L(g_i'^{-1}\cdot g_i)=\lim_{i\to\infty}\bdl^\L(g_{i+1})=\alpha = \bdl(\gamma) = \bdl(\gamma').\] Note that $\bdl(\gamma'^{-1}\cdot\gamma/\gamma)=\bdl(\gamma'/\gamma)=\bdl(\gamma')=\alpha$ (because $\gamma\dind\gamma'$), and similarly $\bdl(\gamma'^{-1}\cdot\gamma/\gamma')=\alpha$,
  and so $\gamma \dind \gamma'^{-1}\cdot\gamma \dind \gamma'$.
  Setting $h := \gamma^{-1}\cdot\gamma' = (\gamma'^{-1}\cdot\gamma)^{-1}$, we have $\gamma \dind h \dind \gamma' = \gamma\cdot h \dind \gamma$, as required by Claim~\ref{c:cohActAb}, and we conclude.
\end{proof}

\begin{remark}
  The limiting process in this proof is necessary, we can not take $\gamma =
  g_i$ for any $i$; see Example~\ref{e:span}.
\end{remark}

As in \cite[Theorem~A.4]{cubicSurfaces}, we can deduce a version for 
non-generic $g$:
\begin{corollary}\label{c:cohActAbSub}
  Let $(G,X)$, $g$, $a$ be as in Theorem~\ref{t:cohActAb}, but do 
  not assume that $g$ is generic in $G$.
  Then $g$ is generic in a left coset over $\acl^0(\emptyset)$ of a connected 
  abelian algebraic subgroup $G' \leq G$.
\end{corollary}
\begin{proof}
  Let $(g_i)_{i \in \omega}$ be as in the proof of Theorem~\ref{t:cohActAb}.
  For large enough $i$, $g_i$ is generic in an $\acl^0(\emptyset)$-definable
  connected algebraic subgroup $G' \leq G$; indeed, if $Y$ is the locus over
  $\acl^0(\emptyset)$ of $g$, then $g_i$ is generic in $(Y^{-1}Y)^i$ for
  $i>0$, which (e.g.\ by Zilber's Indecomposability Theorem) is eventually a
  connected algebraic subgroup. Fix such an $i_0$.

  Now let $X' := G'*a \subseteq X$. We claim that $\dim(X') = \dim(X)$.
  Suppose $\dim(X') < \dim(X)$. Let $a' := g_{i_0}*a$.
  Then $a \nind^0 a'$, since $\dim^0(a/a') \leq \dim(X')$ as $X' = G'*a'$,
  hence $\bdl(a/a') = 0$.
  It follows that $g_{i_0}\dind aa'$, since $g_{i_0}\dind a'$, and hence $g_{i_0} \ind^0 aa'$.
  Then $g'*a = a'$ for any $g' \in G'$, so $a'=a$. Now the set $\bigcap_{g\in G'}\{x\in X: g*x=x\}$ is Zariski-closed and contains $a$. Since $a$ is generic,
  $G'$ is in the kernel of the action, contradicting faithfulness.


  Hence $X'$ is a dense constructible subset of $X$, so contains a Zariski-open dense subset $U \subseteq X$. Then $X' = \bigcup_{g \in G'} g*U$ is an open subvariety of $X$. Also, if $g \in G'$ acts trivially on $X'$ then by density it acts trivially on $X$ (since the set of fixed points of $g$ is closed), so $g = 1$. Hence $(G',X')$ is a CFAHS, and Theorem~\ref{t:cohActAb} applies to show that $G'$ is abelian and 
  $(G',X') \cong (G',G')$.

  Now let $h_0 \in Y(\acl^0(\emptyset))$, where $Y=\loc^0(g)$. Then $g' := h_0^{-1}\cdot g \in G'$.
  We conclude by showing that $g'$ is generic in $G'$.

  First, note that $g' \sim^0_{g'*a} a$ since $\Stab_{G'}(a) = 1$.
  We have $g' \dind g'*a$ since $g \dind g*a$, and similarly $g' \ind^0 g'*a$.
  Then $a \ind^0 g'*a$, since $\tp(a)$ is cgp and $\bdl(a/g'*a) = 
  \bdl(g'/g'*a) = \bdl(g') = \bdl(g) > 0$.
  So $\dim^0(g') = \dim^0(g'/g'*a) = \dim^0(a/g'*a) = \dim^0(a) = \dim(G')$,
  hence $g'$ is generic in $G'$.
\end{proof}

\section{Recognising algebraic homogeneous spaces}
\label{s:CFAHS}
In this section, we apply the group configuration theorem from model theory to 
obtain a recognition result for CFAHSs adapted to our needs.

We first recall the group configuration theorem, specialising the statement in \cite[\S5.4]{GST} for a complete stable theory $T$ to our case, where $T := \Th^0(K)$ is a countable expansion by constants of $\ACF_0$.

\begin{definition}[{\cite[Definition~5.4.3]{GST}}] \label{d:quad}
  A \defn{strict algebraic partial quadrangle} in $K$ is a $6$-tuple $(a,b,c,x,y,z)$ of tuples from $K$ such that:
\begin{enumerate}
\item
$\acl^0(ab)=\acl^0(ac)=\acl^0(bc)$.
\item $x\sim^0_ay$, and $a\sim^0 \cb^0(xy/a)$.
\item $y\sim^0_bz$, and $b\sim^0 \cb^0(yz/b)$.
\item $z\sim^0_cx$, and $c\sim^0 \cb^0(zx/c)$.
\item
Any triple of non-collinear points in the diagram below is $\ind^0$-independent, i.e.\ if $s,t,u$ are points on the diagram not all lying on one of the lines drawn in the diagram, then $s \ind^0 tu$.
\end{enumerate}
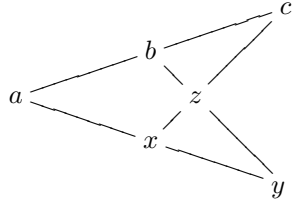
\grpconf abcxyz
  A strict algebraic partial quadrangle in $K$ \defn{over} $C \subseteq K$ is a strict algebraic partial quadrangle after adding constants for $C$ to the language.
\end{definition}

\begin{fact} \label{f:grpConf}
Suppose $(a,b,c,x,y,z)$ is a strict algebraic partial quadrangle in $K$.
Then there is a CFAHS $(G,X)$ defined over $\acl^0(\emptyset)$,
and a tuple $(a',b',c',x',y',z')\in G^3\times X^3$, such that:
\begin{enumerate}[(i)]
  \item each primed element is $\acl^0$-interalgebraic with the corresponding non-primed element;
  \item $a'*x'=y'$, $b'*y'=z'$, $c'*x'=z'$, and $c'=b'\cdot a'$;
  \item each of $a'$, $b'$, and $c'$ is generic in $G$;
  \item each of $x'$, $y'$, and $z'$ is generic in $X$.
\end{enumerate}
\end{fact}
\begin{proof}
  Let $T$ be the expansion of $\ACF_0$ by the constants of $\L$.
  Let $K_0 \leq K$ be an $\aleph_1$-saturated submodel such that $K_0 \ind^0 abcxyz$; such exists by $\aleph_1$-saturation of $K$.
  By the general model theoretic group configuration theorem \cite[Theorem~5.4.5]{GST}, we obtain a $\bigwedge$-definable connected homogeneous space $(G,X)$ over $K_0$ and primed elements in $K$ as in the statement, with the interalgebraicities being over $K_0$. It remains to see that we can replace $K_0$ with $\acl^0(\emptyset)$, and that $(G,X)$ can be taken to be a CFAHS.

  By results of Hrushovski and van den Dries, following Weil, any
  $\bigwedge$-definable group in $\ACF_0$ is a definable group which is
  definably isomorphic to an algebraic group \cite[Theorem~7.4.14, Theorem~7.5.3]{Marker-mt}.
  So we may assume that $G$ is an algebraic group. Let $H \leq G$ be the
  stabiliser of a point in $X(K_0)$, which is a definable and hence closed
  algebraic subgroup. Then $(G,X)$ is definably isomorphic to the left-regular
  action on the left coset space $(G,G/H)$. By standard results,
  $G/H$ is a variety and this action is regular, so $(G,G/H)$ is a CFAHS. So
  we may assume that $(G,X)$ is a CFAHS.\footnote{
  Alternatively, we can cite Weil's original results on generically defined
  groups, which apply also to the present case of homogeneous spaces.
  Let $(g,h,x) \in G\times G\times X$ be $\ACF_0$-generic over $\acl^0(\emptyset )$.
  Then, recalling that any generically defined $\ACF_0$-definable partial map 
  of varieties $V \rightarrow  W$ agrees generically with a rational map,
  we see that $g,h,x$ (as elements of $\loc^0(g)$ and $\loc^0(x)$) satisfy 
  (G1), (G2), (TG1), (TG2) in \cite{Weil-transformations}. By the main theorem 
  of that paper
  \cite[p375]{Weil-transformations}, we may therefore assume that $G$ is an 
  algebraic group and $(G,X)$ is an algebraic homogeneous space.
  Since $G$ is connected and the action is transitive, it follows also that 
  $X$ is irreducible.
  We actually obtain a faithful algebraic homogeneous space by this process, 
  since the original definable action was faithful. We can argue for this as 
  follows. First replace $G$ by a definably isomorphic algebraic group. By DCC 
  for algebraic subgroups of $G$, there are finitely many points $x_1,\ldots ,x_n 
  \in X$ such that their simultaneous stabiliser $G_{x_1,\ldots ,x_n} = \bigcap_i 
  G_{x_i}$ is trivial. Let $g \in G$ be generic over $x_1,\ldots ,x_n$ and the 
  parameters of the $G$-equivariant birational map $\theta : X \rightarrow  X'$ given 
  by Weil's result, where $(G,X')$ is an algebraic homogeneous space. Then 
  each $gx_i$ is generic in $X$ over these parameters, so $\theta(gx_i)$ is
  defined, and then $\bigcap_i G_{\theta(gx_i)} = \bigcap_i G_{x_i}^{g^{-1}}$
  is trivial, so $G$ acts faithfully on $X'$.}

  Let $d_G := \dim(G)$ and $d_X := \dim(X)$. Note that $\dim^0(a) =
  \dim^0(a/K_0) = \dim^0(a'/K_0) = \dim(G) = d_G$ by genericity of $a'$ in $G$
  over $K_0$, and similarly $\dim^0(c)=\dim^0(b)=\dim^0(a)=d_G$ and
  $\dim^0(x)=\dim^0(y)=\dim^0(z)=d_X$.

  Finally, say $e \in K_0$ is such that $G$, $X$, and the interalgebraicities are
  all defined over $e$. The properties we require of $e$ -- that $(G,X)$ is a
  CFAHS with $\dim(G) = d_G$ and $\dim(X) = d_X$,
  and that there are primed elements which satisfy the equalities in (ii)
  and are interalgebraic over $e$ with the corresponding non-primed elements
  -- can be expressed by a formula over $abcxyz$.
  Since $e \ind^0 abcxyz$, this formula is realised by some $e' \in
  \acl^0(\emptyset)$. Indeed, we can see this by the coheir 
  property of non-forking: since $q:=\tp(e/abcxyz)$ does not fork over $\acl(\emptyset)$, it is finitely satisfiable in $\acl(\emptyset)$, namely every formula in $q$ has a realisation in  $\acl(\emptyset)$. Finally, (iii) and (iv) hold because
  $\dim^0(a') = \dim^0(a) = d_G = \dim G$, and similarly for the others.
\end{proof}

We give a name to some of the conditions satisfied by the lines in the above diagram.
\begin{definition} \label{d:algCorrTri}
  A triple $(a,d,b)$ is an \defn{algebraic correspondence triangle} if:
  \begin{itemize}
    \item $a \sim^0_d b$;
    \item $d \sim^0 \cb^0(ab/d)$;
    \item $a \ind^0 d \ind^0 b$.
  \end{itemize}
\end{definition}

\begin{proposition} \label{p:recogCFAHS}
  Suppose $(d,a,b,d',a')$ are tuples from $K$ satisfying the following:
  \begin{enumerate}[(a)]
  \item $(a,d,b)$ and $(a',d',b)$ are algebraic correspondence 
    triangles;
  \item $d'a'\ind^0_bda$.
  \end{enumerate}
  Then setting $e := \cb^0(aa'/dd')$, we have
  $d \sim^0_e d'$.

  Suppose furthermore:
  \begin{enumerate}
    \item[(c)] $\dim^0(d')=\dim^0(d)$;
    \item[(d)] $a\ind^0 b$ and $a'\ind^0 b$;
    \item[(e)] $d \ind^0 e$.
  \end{enumerate}

  Then $(d,d',e,a,b,a')$ is a strict algebraic partial quadrangle in $K$,
  and so by Fact~\ref{f:grpConf} there is a CFAHS $(G,X)$ defined over $\acl^0(\emptyset)$ such that
  $d$ is $\acl^0$-interalgebraic with a generic $d'$ of $G$, and
  $a$ is $\acl^0$-interalgebraic with a generic $a'$ of $X$,
  and $b$ is $\acl^0$-interalgebraic with $d'*a'$.

  \grpconf {d}{d'}{e\;}{a}{b\;}{a'}
\end{proposition}
\begin{proof}
  First, assume only (a) and (b).
  \begin{claim} \mbox{} \label{c:energyIndep}
    $a \ind^0 dd'$ and
    $a' \ind^0 dd'$.
  \end{claim}
  \begin{proof} \mbox{}
    We have $d' \ind^0 dab$, since
    $d' \ind^0 b$ and $d' \ind^0_b da$.
    Thus $a \ind^0_d d'$,
    and since $a \ind^0 d$,
    we conclude $a \ind^0 dd'$.
    By a symmetric argument, also $a' \ind^0 dd'$.
    \subqed{c:energyIndep}
  \end{proof}

  Hence $a \ind^0 dd'e$, since $e\in\acl^0(dd')$. But $a \sim^0_{d'e} b$ since $a \sim^0_e a' \sim^0_{d'} b$, so by Fact~\ref{f:cb} we have $d \sim^0 \cb(ab/d) \sim^0 \cb(ab/dd'e) \in \acl^0(d'e)$.
  Symmetrically, $d'\in\acl^0(de)$, so $d\sim^0_e d'$,
  completing the proof of the first part of the statement.

  Now suppose (c)-(e) also hold. We show that $(d,d',e,a,b,a')$ is a strict algebraic partial 
  quadrangle in $K$. Item (1) in the definition follows from $d\sim^0_e d'$ and $e\in\acl^0(dd')$. Items (2) and (3) follow from (a). Item (4), that $a\sim^0_ea'$ and $e \sim^0 \cb^0(aa'/e)$, follows from $a\sim^0_{dd'}a'$ and the definition of $e$. So it remains to check item (5), that any non-collinear triple is $\ind^0$-independent.

  We have $d\ind^0 e$ by (e), and then by (d) and $d\sim^0_e d'$ we have
  $\dim^0(d'/e)=\dim^0(d/e) = \dim^0(d) = \dim^0(d')$,
  so also $d'\ind^0e$.

  Independence of $(a,d',e)$ follows from Claim~\ref{c:energyIndep} and $d'\ind^0e$; similarly for $(d,a',e)$ (indeed, $a'\ind^0 dd'$ by Claim~\ref{c:energyIndep}, hence $a'\ind^0 de$ and $d\ind^0 e$ by (e).)
  For the remaining non-collinear triples, two of the points are in one of the 
  lines $(d,a,b)$ or $(d',a',b)$, so by the symmetry between primed and unprimed elements (noting that we also have $d'\ind^0 e$, which corresponds under this symmetry to (e)) and since $d \ind^0 
  b$ and $a \in \acl^0(db)$ and $a \ind^0 b$, it suffices to show
  $d' \ind^0 db$ and $a' \ind^0 db$ and $e \ind^0 db$.
  We have $d' \ind^0 db$ since $d' \ind^0 b$ and $d' \ind^0_b d$.
  Then $b \ind^0 dd'$ so $b \ind^0_d e$, and $e \ind^0 d$ by (e), so 
  $e \ind^0 db$.
  Finally, $a' \ind^0_b d$, and $a' \ind^0 b$ by (b) and (d),
  so $a' \ind^0 db$.
\end{proof}

\begin{remark}
  We have stated Proposition~\ref{p:recogCFAHS} in the form we will use it, but note that it implies a direct recognition result for CFAHSs.
  Namely, it follows that
  an algebraic correspondence triple $(a,d,b)$ is co-ordinatewise 
  $\acl^0$-interalgebraic with a generic triple $(g,x,g*x)$ of a CFAHS if and 
  only if: $a \ind^0 b$, and
    when we take $d'a' \equiv^0_b da$ with $d'a'\ind^0_bda$, we have $d \ind^0 \cb^0(aa'/dd')$.
\end{remark}

\section{A homogeneous space version of Elekes-Szabó}
\label{s:main}
In this section, we bring the ingredients together to prove our main result.
We continue to work in the context introduced in \S\ref{ss:prelims2}.

\subsection{Correspondence triangles}
\begin{definition}
  Let $C \subseteq K$.
  A \defn{coarse correspondence triangle} (or just \defn{correspondence 
  triangle}) over $C$ is a triple $(a,d,b)$ such that
  \begin{itemize}
    \item $\tp(a/C)$ and $\tp(d/C)$ are each wgp and broad;
    \item $a \sim^0_{Cd} b$;
  \item $d \sim^0_C \cb^0(ab/Cd)$;
  \item $a \dind_C d \dind_C b$.
  \end{itemize}

  If $\tp(a/C)$ is cgp, we call $(a,d,b)$ a \defn{semi-cgp correspondence 
  triangle} over $C$, and if also $\tp(d/C)$ is cgp, we call it a \defn{cgp 
  correspondence triangle} over $C$.

  In the case $C=\emptyset$, we omit mention of it.
\end{definition}

The triple $(a,g,g*a)$ in Theorem~\ref{t:cohActAb} is a semi-cgp
correspondence triangle,
and the ``wES-triangles'' from \cite[Definition~3.3]{cubicSurfaces} are
correspondence triangles.

\begin{lemma}\label{l:wgp}
  If $(a,d,b)$ is a correspondence triangle over $C$, then $\tp(b/C)$ is also wgp and broad. If $\tp(a/C)$ is cgp, then $\tp(b/C)$ is also cgp.

  In particular, a correspondence triangle is an algebraic correspondence triangle in the sense of Definition~\ref{d:algCorrTri}.
\end{lemma}
\begin{proof}
$\tp(a/Cd)$ is wgp by Fact~\ref{f:wgpFacts}(iii),
so $\tp(b/Cd)$ is wgp by Fact~\ref{f:wgpFacts}(i),
and then $\tp(bd/C)$ is wgp by Fact~\ref{f:wgpFacts}(ii),
and so $\tp(b/C)$ is wgp by Fact~\ref{f:wgpFacts}(i).

Broadness of $\tp(b/C)$ follows from broadness of $\tp(a/C)$,
since $\bdl(b/C) = \bdl(b/Cd) = \bdl(a/Cd) = \bdl(a/C)$.

If $\tp(a/C)$ is cgp,
then $\tp(b/Cd)$ is cgp by Lemma~\ref{l:cgppreserve}(i) and (ii), and then we
conclude by Lemma~\ref{l:cgppreserve}(iii) that $\tp(b/C)$ is cgp.

For the ``in particular'' clause, just observe that using wgp we have
$a \ind^0_C d \ind^0_C b$.
\end{proof}

\begin{definition}
  Let $C \subseteq K$.
  A correspondence triangle $(a,d,b)$ over $C$ is \defn{abelian} over $C$ if
  there exist a connected commutative algebraic group $G$ over $\acl^0(C)$ and 
  $g_a,g_d,g_b \in G$ each generic in $G$ over $\acl^0(C)$ such that
  $g_a \sim^0_C a$ and $g_d \sim^0_C d$ and $g_b \sim^0_C b$,
  and $g_b = g_d+g_a$.

  Again, we omit ``over $C$'' when $C = \emptyset$.
\end{definition}

Our main result, Theorem~\ref{t:main} below, is that every semi-cgp 
correspondence triangle is abelian.

\begin{remark}
  Note that the asymmetric interalgebraicity and algebraic independence 
  properties of correspondence triangles become symmetric in the case of an 
  abelian correspondence triangle $(a,d,b)$: we have $d \sim^0_a b$ and $d 
  \sim^0_b a$, and $\dim^0(a) = \dim^0(d) = \dim^0(b) = \dim(G)$, and $a \ind^0 
  b$ since $\dim^0(a/b) = \dim^0(d/b) = \dim^0(d) = \dim^0(a)$.
\end{remark}

\begin{question}
  One could analogously define a correspondence triangle $(a,d,b)$ to be 
  \emph{nilpotent} if
  there exist a connected nilpotent algebraic group $G$ acting on a variety 
  $W$, all over $\acl^0(\emptyset)$, a possibly non-generic element $g_d \in 
  G$, and generics $w_a,w_b \in W$, such that
  $g_d \sim^0 d$ and $w_a \sim^0 a$ and $w_b \sim^0 b$,
  and $w_b = g_d*w_a$.

  Is every correspondence triangle nilpotent?
\end{question}

\begin{remark} \label{r:dcl}
  If in an abelian correspondence triangle one adds the condition $b \in 
  \dcl^0(ad)$ as opposed to merely $b \in \acl^0(ad)$, one can correspondingly 
  strengthen the interalgebraicities witnessing abelianity to rational maps.
  Although we do not explicitly use this in this paper, we take the
  opportunity to sketch an argument based on unpublished work of Hrushovski.
  We use an elementary form of this argument in the proof of 
  Theorem~\ref{t:unbalancedERMany}.

  Assume (adding parameters) that $\acl^0(\emptyset) = \dcl^0(\emptyset)$, so 
  $\tp^0(e)$ is stationary for any $e$.
  Suppose $g_a' \equiv^0_a g_a$, and set $\alpha := g_a' - g_a$.
  Then $g_a' \equiv^0_{adg_d} g_a$, since $a \ind^0 d$ and $\tp^0(g_aa)$ is 
  stationary.
  Since $b \in \dcl^0(ad)$, it follows that $g_b \equiv^0_b g_b + \alpha$.
  But $\alpha \ind^0 b$ since $a \ind^0 b$,
  and it follows (via stationarity) that $\tp(g_b/b)$ is invariant under 
  addition by $\alpha$.
  So, replacing $G$ by its quotient by the finite subgroup stabilising the 
  finite set $\tp(g_b/b)$, we reduce to $g_a \in \dcl^0(a)$, and similarly to 
  $g_d \in \dcl^0(d)$.

  In particular, if either $a$ or $d$ is generic in a rational curve, then $G$ 
  is also a rational curve, so isomorphic to either the additive or the 
  multiplicative group.

  If we also replace $b$ by a $b' \in \acl^0(b)$ such that $\dcl^0(b') = 
  \acl^0(b) \cap \dcl^0(ad)$, then we also obtain $g_b \in \dcl^0(b)$ (since 
  $g_b \in \acl^0(b) \cap \dcl^0(g_ag_d) \subseteq \acl^0(b) \cap 
  \dcl^0(ad)$).
\end{remark}

As with all results of Elekes-Szabó type, our arguments are eventually based 
on incidence bounds of Szemerédi-Trotter type. The following lemma 
encapsulates this in the form we need. We can see this as a kind of 
``$\bdl$-pseudolinearity'' -- see Remark~\ref{r:pseudolin}.
\begin{lemma} \label{l:STCb}
  Let $(a,d,b)$ be a semi-cgp correspondence triangle.
  Then $\bdl(d) \leq \bdl(a)$.
\end{lemma}
\begin{proof}
  We may assume $d = \cb^0(ab/d)$.
  Exactly as in the proof of \cite[Proposition~5.14]{BB-cohMod}, applying 
  Szemerédi-Trotter bounds to the $\bigwedge$-definable binary relation 
  $(\tp(ab)\times\tp(d)) \cap \loc^0(ab,d)$ and using that $\tp(a)$ is cgp,
  we obtain
  \begin{equation} \label{e:ST}
    \bdl(abd) \leq  \max(\frac12 \bdl(ab) + \bdl(d) - 
    \max(0,\eps_0(\bdl(d) - \frac12 \bdl(ab))), \bdl(ab), 
    \bdl(d))
  \end{equation}
  where $\eps_0 > 0$.

  Suppose $\bdl(d) > \bdl(a)$.
  We have $\bdl(b) = \bdl(b/d) = \bdl(a/d) = \bdl(a)$,
  so $\frac12 \bdl(ab) \leq  \bdl(a) < \bdl(d)$,
  and $\max(0,\eps_0(\bdl(d) - \frac12 \bdl(ab)))>0$.
  But $\bdl(abd) = \bdl(a/d) + \bdl(d) = \bdl(a) + \bdl(d) \geq  \frac12 \bdl(ab) 
  + \bdl(d)$,
  and $\frac12 \bdl(ab) + \bdl(d) > \frac12\bdl(ab) + \frac12\bdl(ab) = 
  \bdl(ab)$,
  and $\frac12\bdl(ab) + \bdl(d) > \bdl(d)$,
  so this contradicts \eqnref{ST}.
\end{proof}

\begin{lemma}\label{l:tribase}
  If $(a,d,b)$ is a (semi-cgp) correspondence triangle, and $\tp(c)$ is wgp,
  and $c \dind_d ab$, and $\tp(d/c)$ is wgp and broad,
  then $(a,d,b)$ is a (semi-cgp) correspondence triangle over $c$.
\end{lemma}
\begin{proof}
  We have $a \dind dc$ (since $d\dind a$ and $c \dind_d a$) and hence $\tp(a/c)$ is wgp (resp.\ cgp) and broad (by Fact~\ref{f:wgpFacts} resp.\ Lemma~\ref{l:cgppreserve})
  and $a \dind_c d$, and similarly for $b$.

  It remains to show that $d \sim^0_c \cb^0(ab/dc)$.
  By Fact~\ref{f:cb} (or directly from the definition of $\cb^0$), it suffices 
  to show $ab\ind^0_d c$, or equivalently $a\ind^0_d c$.
  But indeed, $a \ind^0 dc$ since $\tp(a)$ is wgp and $a\dind dc$.
\end{proof}


\begin{lemma}\label{l:vee}
  Suppose $(a,d,b)$ and $(a',d',b)$ are correspondence triangles such that
  $a'd' \dind_b ad$.
  Let $e:=\cb^0(aa'/dd')$.
  Then $(a,e,a')$ is a correspondence triangle, and moreover $d\sim^0_ed'$ and $\dim^0(e)\geq\dim^0(d)$.
\end{lemma}
\begin{proof}
  We first verify the conditions for $(a,e,a')$ to be a correspondence triangle.
  \begin{itemize}
    \item $\tp(a)$ is wgp and broad by assumption,
      $e \in \acl^0(dd')$ is wgp by $d\dind d'$ and Lemma~\ref{l:gpSeq}(i),
      and we verify below that $\tp(e)$ is broad.
    \item $a \sim^0_d b \sim^0_{d'} a'$,
      hence $a \sim^0_{dd'} a'$, and so also $a \sim^0_e a'$ by choice of 
      $e$.
    \item $e = \cb^0(aa'/e)$ by definition.
    \item
      By transitivity and monotonicity (as in Claim~\ref{c:energyIndep})
      we have
      $a \dind dd'$ and $a' \dind dd'$,
      so $a \dind e \dind a'$
      since $e \in \acl^0(dd')$.
    \end{itemize}
    Now we claim that $(d,a,b,d',a')$ satisfies (a) and (b) in Proposition~\ref{p:recogCFAHS}. Condition (a) follows directly from the assumptions. For (b), we need to prove $a'd'\ind^0_bad$. By assumption, $a'd'\dind_bad$, hence $d'\dind_bd$ and $d'\dind bd$ (as $d'\dind b$). Since $\tp(d')$ is wgp, we have $d'\ind^0bd$ and $d'\ind^0_bd$. Since $a\in\acl^0(bd)$ and $a'\in\acl^0(bd')$, we conclude $a'd'\ind^0_bad$ as desired.

  So by Proposition~\ref{p:recogCFAHS}, $d'\sim^0_e d$.
  Using $e\in \acl^0(dd')$, we deduce $e\sim^0_{d'} d$.
  Since $d\dind d'$, we obtain
  $\bdl(e)\geq \bdl(e/d') = \bdl(d/d') = \bdl(d) > 0$,
  and hence $\tp(e)$ is broad.
  Finally, the analogous calculation using $d \ind^0 d'$ yields
  $\dim^0(e) \geq \dim^0(d)$.
\end{proof}

\begin{lemma}\label{l:semicgpind}
If $(a,d,b)$ is a semi-cgp correspondence triangle, then $a\ind^0 b$.
\end{lemma}
\begin{proof}
Suppose not. Then $\bdl(a/b)=0$ since 
  $\tp(a)$ is cgp, and so $d\dind_ba$. By $d\dind b$ and 
  transitivity, we obtain $d\dind ab$, and hence 
  $d\ind^0ab$ since $\tp(d)$ is wgp.
  But $d = \cb^0(ab/d)$, so we 
  conclude $d\in\acl^0(\emptyset)$, contradicting broadness of $d$.
\end{proof}

\subsection{Abelianity of cgp correspondence triangles}
Our main theorem that semi-cgp correspondence triangles are abelian will be 
bootstrapped via first proving this for cgp correspondence triangles. In 
outline, we prove this weaker result as follows. We iteratively compose 
independent copies of the correspondence of the original triangle. By 
Proposition~\ref{p:recogCFAHS}, we get a group if the algebraic dimension of 
the correspondences stabilises. The cgp assumption bounds algebraic dimension 
by $\bdl$ (via Lemma~\ref{l:gpSeq}(ii)), and $\bdl$ is bounded throughout the 
iteration by the Szemerédi-Trotter bound of Lemma~\ref{l:STCb}. So we 
eventually get a group, which is commutative by Theorem~\ref{t:cohActAb}. From 
this, it turns out that the iteration was not in fact required, and the 
original correspondence triangle is already abelian.

\begin{remark}\label{r:pseudolin}
  This argument is analogous to the proof of the result of Buechler and Hrushovski that 
  ``pseudolinear'' strongly minimal sets are locally modular 
  (\cite[Proposition~5.3.2]{GST}, \cite[Theorem~B]{buechler}).
  There, one similarly composes and uses a 
  uniform upper bound on dimension of correspondences (this is the definition 
  of pseudolinearity) to reach a group configuration, then uses minimality to 
  see that the group must be abelian and no iteration was necessary. The 
  Szemerédi-Trotter bound of Lemma~\ref{l:STCb} directly bounds the 
  $\bdl$-dimension of correspondences, and in our asymmetric situation where 
  we may have $\bdl(d)<\bdl(a)$, this gives rise to pseudolinearity in the 
  sense of a bound $\frac{\bdl(a)}{\bdl(d)}$ on relative dimensions of 
  correspondences obtained by composing the original one.
  Contrast with the symmetric situation of \cite{BB-cohMod}, where the 
  Szemerédi-Trotter bound directly gives local modularity. See also the 
  discussion after \cite[Theorem~5.4]{Hr-psfDims}, which interprets another 
  formulation of Szemerédi-Trotter as a statement of pseudolinearity.
\end{remark}

\begin{proposition} \label{p:maincgp}
  Every cgp correspondence triangle is abelian.

 Moreover, suppose $(a,d,b)$ and $(a',d',b)$ are cgp correspondence triangles such that $a'd'\dind_b ad$, $\bdl(d)=\bdl(d')$ and $\dim^0(d)=\dim^0(d')$. Let $e:=\cb^0(aa'/dd)$. Then $(a,d,b)$ is an abelian correspondence triangle, and $(d,d',e,a,b,a')$ is a strict algebraic partial quadrangle.
\end{proposition}
\begin{proof}
  It is enough to prove the moreover part, since given a cgp correspondence triangle $(a,d,b)$ we can take $a'd'\equiv_bad$ such that $a'd'\dind_b ad$.

  Let $(a,d,b)$, $(a',d',b)$, and $e$ be as in the statement.
  Rescaling, we may assume $\dim^0(d) = \bdl(d)=\bdl(d')=\dim^0(d')$.

  Note that $d\dind a'd'$, as $d\dind_b a'd'$ and $d\dind b$.

  \newcommand{\bh}{\bar h}
 We recursively define a sequence $(d_i,d_i',b_i)_{i \in \omega}$ such that, for all $i>0$,
  \begin{enumerate}[(i)]
    \item $(a,d_i,b_i)$ is a semi-cgp correspondence triangle,
    \item $\dim^0(d_i)\geq\dim^0(d_{i-1})$,
    \item $b_id_{i-1}' \dind_{b_{i-1}} ad_{i-1}$,
      and $\dim^0(d_{i-1}') = \dim^0(d_{i-1})$,
      and $(b_i,d_{i-1}',b_{i-1})$ is a semi-cgp correspondence triangle,
    \item $d_{i-1} \sim^0_{d_i} d_{i-1}'$,
    \item $d_i \in \acl^0(\bh)$ for some finite $\bdl$-independent sequence $\bh$ of realisations of $\tp(dd')$ with $\bh \dind b_i$ and $\bh \dind a$.
  \end{enumerate}

  First, let $(d_0,d_0',b_0,d_1,b_1) := (d,d',b,e,a')$. Then (i),(ii),(iv) hold for $i=1$ by Lemma~\ref{l:vee}, (iii) holds by assumption, and (v) holds since $e\in\acl^0(dd')$ and $dd'\dind a'$ (since $d\dind a'd'$ and $a'\dind d'$) and similarly $dd'\dind a$.

  Suppose $i\geq 1$ and we have defined $(d_{i-1},d_{i-1}',b_{i-1},d_i,b_i)$ satisfying (i)-(v). We define $d_i'$, $d_{i+1}$, and $b_{i+1}$ such that $(d_i,d_i',b_i,d_{i+1},b_{i+1})$ satisfy (i)-(v).

  Let $\bh$ be as in (v) for $d_i$.
  Let $b_{i+1}d_i'\bh' \equiv_{b_i} ad_i\bh$ with $b_{i+1}d_i' \bh'\dind_{b_i} ad_i\bh$,
  and set $d_{i+1} := \cb^0(ab_{i+1}/d_id_i')$.
  Then (i),(ii),(iv) hold for $i+1$ by Lemma~\ref{l:vee}, and (iii) holds by construction. By definition, $d_{i+1}\in\acl^0(d_id_i') \subseteq \acl^0(\bh\bh')$. Note that $\bh\dind \bh'$, since $\bh\dind b_i$ by assumption. Therefore, $\bh\bh'$ is a $\bdl$-independent sequence in $\tp(dd')$. Moreover, we have $b_{i+1}\bh'\dind \bh$ (again by $\bh\dind b_i$).
  Now $\bh' \dind b_{i+1}$ since $\bh \dind a$,
  so $\bh\bh'\dind b_{i+1}$.
  Also, $\bh' \dind b_i$ since $\bh \dind b_i$, so $\bh' \dind \bh a$, and since also $\bh \dind a$, we have $\bh\bh' \dind a$.
  So (v) holds, witnessed by $\bh\bh'$.

  This concludes our construction of the sequence 
  $(d_i,d_i',b_i)_{i\in\omega}$.

  Now for all $i\in\omega$, we have:
  \begin{itemize}
    \item $a\ind^0d_i\ind^0b_i\ind^0a$ and 
      $b_{i+1}\ind^0d'_i\ind^0b_i\ind^0b_{i+1}$ by (i) (which also holds for 
      $i=0$), (iii), and Lemma~\ref{l:semicgpind};
    \item $d_i'\ind^0 b_id_i$, since $d_i'\dind b_id_i$ (by (iii) $d_i'\dind_{b_i}d_i$ and (i) $d_i'\dind b_i$) and $d_i'$ is wgp,
      and hence $b_{i+1}d_i'\ind^0_{b_i}ad_i$ (by $b_{i+1}\in\acl^0(d_i'b_i)$ and $a\in\acl^0(b_id_i)$).
  \end{itemize}

  \begin{claim} \label{c:cbDimBd}
    For all $i \in \omega$,
    $\dim^0(d_i) \leq \dim^0(d_{i+1}) \leq \bdl(d_{i+1}) \leq \bdl(a)$.
  \end{claim}
  \begin{proof}
    The first inequality is by (ii).

    By (v), $d_{i+1} \in \acl^0(\bh)$ for some $\bdl$-independent sequence 
    $\bh$
    of realisations of $\tp(dd')$. But $d \dind d'$, so $\bh$ is a
    $\bdl$-independent sequence of realisations of the types $\tp(d)$ and
    $\tp(d')$, which are cgp. So Lemma~\ref{l:gpSeq}(ii) yields the second
    inequality (since we assumed $\dim^0(d) = \bdl(d)$).

    Finally, (i) and Lemma~\ref{l:STCb} give the third inequality.
    \subqed{c:cbDimBd}
  \end{proof}

  By Claim~\ref{c:cbDimBd}, $\dim^0(d_{i+1}) = \dim^0(d_i)$ for some $i\in\omega$.
  Then by $d_i \sim^0_{d_{i+1}} d_i'$, we have $d_i \ind^0 d_{i+1}$ (since 
  $\dim^0(d_i) = \dim^0(d_{i+1}) \geq \dim^0(d_{i+1}/d_i) \geq 
  \dim^0(d_i'/d_i) = \dim^0(d_i')=\dim^0(d_i)$ (using $d_i\ind^0d_i'$)), so 
  $\dim^0(d_{i+1}) = \dim^0(d_{i+1}/d_i)$).

  Hence, as illustrated in the picture below, Proposition~\ref{p:recogCFAHS} applies to $(d_i,a,b_i,d_i',b_{i+1})$ 
  and $d_{i+1}=\cb^0(ab_{i+1}/d_id_i')$ (conditions (a)-(d) have been shown 
  for all $i$, while (e) is the condition $d_i\ind^0d_{i+1}$ just shown for 
  our choice of $i$), so $(d_i,d_i',d_{i+1},a,b_i,b_{i+1})$ is a strict 
  algebraic quadrangle yielding a CFAHS $(G,X)$ defined over 
  $\acl^0(\emptyset)$.

  \grpconf {d_i}{d_i'}{d_{i+1}}a{b_i}{b_{i+1}}

 Since $\tp(a)$ is cgp and $\tp(d_i)$ is wgp, $G$ is commutative by
  Theorem~\ref{t:cohActAb}, and up to interalgebraicity, $(a,d_i,b_i)$ is
  generic in the graph of addition of $G$.
  So $(a,d_i,b_i)$ is abelian.

  It follows that $d_i \sim^0_{b_i} a$, hence $\dim^0(d_i) = \dim^0(d_i/b_i) = 
  \dim^0(a/b_i) = \dim^0(a)$. But $\dim^0(d_0) = \dim^0(d_0/b_0) \geq 
  \dim^0(a/b_0) = \dim^0(a)$, and so by (ii) we conclude that already 
  $\dim^0(d_1)=\dim^0(d_0)$. In other words, we may actually take $i=0$ in the 
  above argument to show that $(d_0,d_0',d_1,a,b_0,b_1) = (d,d',e,a,b,a')$ is a strict 
  algebraic quadrangle and $(a,d,b)$ is abelian.
\end{proof}


\subsection{Abelianity of semi-cgp correspondence triangles}
We now proceed to the proof of our main theorem, Theorem~\ref{t:main}, which 
improves on Proposition~\ref{p:maincgp} by proving abelianity of a semi-cgp 
correspondence triangle $(a,d,b)$ without assuming that $d$ is cgp. In 
outline, this will go as follows. We find $c$ such that $\tp(d/c)$ is cgp and 
apply Proposition~\ref{p:maincgp} over $c$, yielding a family of groups 
$(G_c)_c$. Using the ``moreover'' clause in Proposition~\ref{p:maincgp}, we find 
that these groups are pairwise isogenous. Composing isogenies, we see that $d$ 
can be identified with an element of an algebraic group of affine linear 
transformations of such a group, allowing us to conclude via 
Corollary~\ref{c:cohActAbSub}.

We begin with some preliminary definitions and results.

\begin{definition}
  If $G$ is a commutative algebraic group over $C \subseteq K$ and $x,y \in 
  G(K)$, we write $x \approx_C y$ to mean that $y-x \in G(\acl^0(C))$.
  We abbreviate $x\approx_{\emptyset} y$ to $x\approx y$.
\end{definition}

\begin{fact}\label{zieglerIsog}
  Suppose $G$ and $H$ are connected commutative algebraic groups over $\acl^0(\emptyset)$,
  and $(x,a) \in G^2$ and $(x',a') \in H^2$ are generic,
  and $x \sim^0 x'$ and $a \sim^0 a'$ and $x+a \sim^0 x'+a'$.
  Then there exist $n \in \N_{>0}$ and an isogeny $\alpha : G \to H$ over 
  $\acl^0(\emptyset)$ such that $\alpha x \approx nx'$.
\end{fact}

\begin{proof}
  This is essentially the ``moreover'' clause of \cite[Fact~2.13]{BB-cohMod}; 
  that $\alpha$ is over $\acl^0(\emptyset)$ was not specified there, but this 
  follows from the proof, since the graph of $\alpha$ is obtained from the 
  stabiliser of $\tp^0(xx'/\acl^0(\emptyset))$ (as explained in 
  \cite[Lemme~2.4]{BMP-beauxGroupes}).
\end{proof}

\begin{fact}\label{f:aut}
  Let $G$ be a commutative algebraic group over $\acl^0(\emptyset)$. Denote by $ \Aut(G)$ the group of algebraic group automorphisms of $G$.
  Then there exists a connected algebraic group $\Aut^0(G) \leq \Aut(G)$ over 
  $\acl^0(\emptyset)$ such that any $\sigma \in 
  \Aut(G)$ is of the form $\sigma = \sigma'\circ\sigma_0$ where $\sigma_0 \in 
  \Aut(G)$ is defined over $\acl^0(\emptyset)$ and $\sigma' \in \Aut^0(G)$.
\end{fact}
\begin{proof}
  We first recall from \cite[Lemma~3.1]{bgth} the following facts about the 
  structure of commutative algebraic groups in characteristic 0.
  \begin{itemize}
    \item There exists an exact sequence $0 \to V \to G \to S \to 0$
      where $V$ is a vector group and $S$ is a semiabelian variety,
      all over $\acl^0(\emptyset)$.
    \item $G = G_0 \oplus V_0$ where $G_0 := G[\infty]^{\operatorname{Zar}}$
      is the Zariski closure of the torsion
      and $V_0 \leq V$ is a vector subgroup defined over $\acl^0(\emptyset)$ 
      such that $V = (V\cap G_0) \oplus V_0$.
  \end{itemize}

  Let $\Aut^0(G)$ be the kernel of the restriction homomorphism $\restriction_{G_0} : \Aut(G) \to \Aut(G_0)$.
  If $\sigma_0 \in \Aut(G_0)$ then $\sigma_0$ is determined by its action on
  $G[\infty] \subseteq \acl^0(\emptyset)$ by the density and the fact that $\sigma_0$ is an algebraic automorphism, so $\sigma_0$ is over
  $\acl^0(\emptyset)$.
  Now if $\sigma\in\Aut(G)$, then $\sigma\negmedspace\restriction_{G_0}$ extends uniquely to $\sigma_0\in\Aut(G)$ over $\acl^0(\emptyset)$ with $\sigma_0\negmedspace\restriction_{V_0}=\operatorname{id}_{V_0}$, and then $\sigma':=\sigma\circ\sigma_0^{-1}\in\Aut^0(G)$ as required.
  So it remains only to see that $\Aut^0(G)$ is a
  connected algebraic group.

  Let $L \leq \GL(V)$ be the stabiliser of $V\cap G_0$ in $\GL(V)$.
  Now $V$ is the (unique) maximal vector subgroup of $G$, since no non-trivial
  homomorphism $\G_a \to S$ exists, so any element of $\Aut^0(G)$ restricts to
  an element of $L$.
  This yields an isomorphism $\restriction_V : \Aut^0(G) \to L$, since
  any $\theta \in L$ extends uniquely to $\sigma \in \Aut^0(G)$ defined on $G
  = G_0 \oplus V_0$ by $\sigma(g_0+v_0) := g_0 + \theta(v_0)$.
  Finally, $L$ is a connected algebraic group over $\acl^0(\emptyset)$, as required.
\end{proof}

\begin{lemma}\label{l:coheir}
  If $(a,d,b)$ is a correspondence triangle, and $f \ind^0 adb$ and $(a,d,b)$ 
  is abelian over $f$, then $(a,d,b)$ is abelian.
\end{lemma}
\begin{proof}
  $\tp^0(f/adb)$ is then finitely satisfiable in $\acl^0(\emptyset)$, and so 
  the existence of $G$ over $\acl^0(f)$ of a certain dimension and $g_a 
  \sim^0_f a$, $g_d \sim^0_f d$, and $g_b \sim^0_f b$ with $g_b = g_a + g_d$ 
  implies the corresponding existence over $\acl^0(\emptyset)$, and the 
  genericity follows since $\dim^0(a) = \dim^0(a/f)$ etc.
\end{proof}

Given $a$, $B$, and $C$, we denote by $\cb^0_C(a/B)$ the value of $\cb^0(a/B)$ 
in the language $\L(C)$ obtained by adding constants for $C$.

\begin{lemma}\label{l:cbbase}
  For all $a$, $B$, and $C$, we have $\cb^0_C(a/B) \sim^0_C \cb^0(a/BC)$.
\end{lemma}
\begin{proof}
  We use the description above of $\cb^0$ in terms of fields of definition (see the discussion before Fact~\ref{f:cb}).
  Let $E_0$ be the constants of $\L$ (which form a field by assumption).
  Let $k$ be the field of definition of $\loc^0(a/\acl^0(BC)) = 
  \loc^0(a/\acl^0_C(B))$.
  Then $E_0(\cb^0(a/BC)) = E_0(k)$
  and $E_0(C,\cb^0_C(a/B)) = E_0(C,k)$,
  so $\acl^0(C,\cb^0_C(a/B)) = \acl^0(C,k) = \acl^0(C,\cb^0(a/BC))$, as 
  required.

  Alternatively, in stability theoretic terms: $\acl^0_C(\cb^0_C(a/B))$ is the 
  smallest $\acl^0_C$-closed set over which $\tp(a/\acl^0(BC))$ doesn't fork, 
  and so is $\acl^0(C,\cb^0(a/BC))$.
\end{proof}

\begin{lemma}\label{l:triover}
  If $(a,d,b)$ is a correspondence triangle over $C=\acl^0(C)$, possibly semi-cgp or cgp,  then $(a,d,b)$ is a 
  correspondence triangle (respectively semi-cgp or cgp) in the language $\L(C)$ obtained by adding 
  constants for $C$.
\end{lemma}
\begin{proof}
  The only property which is not quite immediate is that $d \sim^0_C
  \cb^0_C(ab/d)$, but by Lemma~\ref{l:cbbase} this follows from $d \sim^0_C 
  \cb^0(ab/Cd)$.
\end{proof}

\begin{theorem}\label{t:main}
  Every semi-cgp correspondence triangle is abelian.
\end{theorem}
\begin{proof}
  Let $(a,d,b)$ be a semi-cgp correspondence triangle.

  If $\tp(d)$ is cgp, we are done by Proposition~\ref{p:maincgp}. So suppose not.

  By Lemma~\ref{l:wgpification}(ii),
  let $c$ be such that $\tp(d/c)$ is cgp and broad and $\bdl(c) < \infty$.
  Suppose for a contradiction that $\bdl(c) = 0$. Then $d \ind^0 c$ since $\tp(d)$ is wgp, and then if $e$ is such that $d \nind^0 e$, then also $d\nind^0_c e$, so $\bdl(d/ec)=0$, hence also $\bdl(d/e)=0$. This contradicts $\tp(d)$ not being cgp.
  So $\tp(c)$ is broad.
  We may assume $c \dind_d ab$.

  By Lemma~\ref{l:wgpification}(i), there exists $f$ with $\bdl(f)=0$ such 
  that $\tp(c/f)$ is wgp. Then $f\dind ad$ and $f \ind^0 adb$ since $\tp(ad)$ is wgp by Remark~\ref{r:tuplewgp}. Adding $\acl^0(f)$ as constants to the language
  preserves the assumptions that $(a,d,b)$ is a semi-cgp correspondence triangle and that $\tp(d/c)$ is cgp. Indeed, since $d\dind_c f$ and $\tp(d/c)$ is broad and cgp, we get $\tp(d/cf)$ is broad and cgp by Lemma~\ref{l:cgppreserve}(ii) and independence.
  Similarly, $\tp(a/f)$ is cgp and $\tp(d/f)$ is wgp (by Fact~\ref{f:wgpFacts}(iii)).
 Thus by Lemma~\ref{l:coheir}, we may assume that $\tp(c)$ is wgp.

  By Lemma~\ref{l:tribase}, $(a,d,b)$ is a cgp correspondence triangle over 
  $c$. So by Proposition~\ref{p:maincgp}
  applied via Lemma~\ref{l:triover} in the expansion $\L(\acl^0(c))$,
  there is a 
  connected commutative algebraic group $(G_c,+)$ over $\acl^0(c)$ such that $d,a,b$
  are interalgebraic over $c$ with generics
  $d_c,a_c,b_c \in G_c$ respectively, and $b_c = a_c + d_c$.
  In particular, $d \in \acl^0(cab)$.

  Let $G := G_c$.

  \begin{claim}\label{c:cind}
  $a\dind dc$ and $b\dind dc$.
  \end{claim}
  \begin{proof}
  Since $c\dind_d a$ and $d\dind a$, by transitivity, we have $cd\dind a$. Similarly, we get $b\dind dc$.
  \subqed{c:cind}
  \end{proof}

  Take $d'a'c'$ such that $d'a'c' \equiv_{\acl^0(b)} dac$ and $d'a'c' 
  \dind_b dac$.
  \begin{claim}\label{c:cc'inde}
  $a\dind cc'dd'$, $b\dind cc'dd'$, and $a'\dind cc'dd'$. Thus, since $\tp(a)$ and $\tp(b)$ are wgp (by Lemma~\ref{l:wgp}), the corresponding independences also hold for $\ind^0$.
  \end{claim}
  \begin{proof}
  By symmetry, it is enough to prove $a\dind cc'dd'$. We have $d'c' 
  \dind_b dac$, hence  \begin{equation}\label{eq:dcind}
  d'c' \dind dcab
  \end{equation} by $d'c'\dind b$ (Claim~\ref{c:cind} and $d'c'\equiv_bdc$). Thus, 
  $d'c'\dind_{dc}a$
   and since $dc\dind a$ by Claim~\ref{c:cind}, we conclude $dcd'c'\dind a$.
  \subqed{c:cc'inde}
  \end{proof}

  Let $e := \cb^0_{cc'}(aa'/dd')$.
  \begin{claim}\label{c:quadrangle}
  $(d,d',e,a,b,a')$ is a strict algebraic partial quadrangle over $\acl^0(c,c')$, in particular $a\sim^0_{ecc'}a'$ and $d\sim^0_{ecc'}d'$.
  \end{claim}
  \begin{proof}
We claim that the conditions for the moreover part in Proposition~\ref{p:maincgp} hold for $(a,d,b)$ and $(a',d',b)$ when we pass to the expansion $\L(\acl^0(c,c'))$. Since $(a,d,b)$ is a cgp correspondence triangle over $c$ and $c'\dind_{dc}ab$ by (\ref{eq:dcind}), and $\tp(d/cc')$ is broad and cgp by $d\dind_c c'$ (which again follows from (\ref{eq:dcind})) and Lemma~\ref{l:cgppreserve}(ii), we conclude that $(a,d,b)$ is a cgp correspondence triangle over $\acl^0(c,c')$ by Lemma~\ref{l:tribase}.

  By symmetry, $(a',d',b)$ is also a cgp correspondence triangle over $\acl^0(c,c')$. Moreover, $\bdl(d/cc')=\bdl(d/c)=\bdl(d/bc)=\bdl(d'/bc')=\bdl(d'/c')=\bdl(d'/cc')$ by $dc\dind c'$, $dc\dind b$ and symmetrically $d'c'\dind c$ and $d'c'\dind b$. Since $\tp(c),\tp(c'),\tp(b)$ are each wgp, an analogous calculation shows $\dim^0(d/cc')=\dim^0(d'/cc')$. We only need to show $a'd'\dind_{bcc'}ad$. By (\ref{eq:dcind}), we have $d'c'\dind dcab$, hence $d'\dind_{bcc'}da$ and $a'd'\dind_{bcc'}da$ as $a'\in\acl^0(bd')$. We conclude by applying Proposition~\ref{p:maincgp} in $\L(\acl^0(c,c'))$, via Lemma~\ref{l:triover}.
  \subqed{c:quadrangle}
  \end{proof}

  \begin{claim}\label{c:isoms}
    There exist a connected algebraic group $H$ over $\acl^0(c')$
    and isomorphisms $\theta_a,\theta_b : G \to H$ over $\acl^0(cc')$
    and $a_H,b_H \in H(K)$
    such that
    $$\theta_a(a_c) \approx_{cc'} a_H \sim^0_{c'} a
    \text{ and } \theta_b(b_c) \approx_{cc'} b_H \sim^0_{c'} b.$$
  \end{claim}
  \begin{proof}
    Take $a'_{c'},b_{c'},d'_{c'}$ such that $a'_{c'}b_{c'}d'_{c'}a'd'c' 
    \equiv_{\acl^0(b)} a_cb_cd_cadc$, so $b_{c'} = a'_{c'} + d'_{c'}$ in 
    $G_{c'}$, a connected algebraic group over $\acl^0(c')$ in which $b_{c'}$ is 
    generic over $\acl^0(c')$,
    and $b_{c'} \sim^0_{c'} b$.

    We have $a'_{c'}\sim^0_{c'}a'\sim^0_{ecc'} a$ by Claim~\ref{c:quadrangle}, and similarly, $d'_{c'}\sim^0_{c'}d'\sim^0_{ecc'} d$.
    Now $\tp(e/\acl^0(abdcc'))$ is finitely satisfiable in $\acl^0(cc')$, 
    since $e \ind^0_{cc'} abd$ by Claim~\ref{c:quadrangle}.
    so there exist $a'_0 \sim^0_{cc'}a$ and $d'_0\sim^0_{cc'}d$ such that 
    $a'_0,d'_0 \in G_{c'}$ and $b_{c'} = a'_0 + d'_0$,
    since this existence statement holds over $ecc'$ witnessed by $a'_{c'}$ and $d'_{c'}$.
    By Fact~\ref{zieglerIsog},
    $\alpha b_{c'} \approx_{cc'} nb_c$ for some $n \in \N_{>0}$ and isogeny $\alpha 
    : G_{c'} \to G$ over $\acl^0(cc')$. Let $H := G_{c'}/\ker\alpha$. Note that $H$ is divisible since $G_{c'}$ is (by the fact that $G_{c'}$ is connected).  Taking $b_H$ to be an $n$-division point of $b_{c'}/\ker\alpha$, we obtain 
    an isomorphism $\theta_b$ as required,
    noting that $H$ is over $\acl^0(c')$ since $\ker\alpha\leq G_{c'}$ is 
    finite and hence consists of finitely many torsion elements of $G_{c'}$,
    each of which is in  $G_{c'}(\acl^0(c'))$, since each $n$-torsion group $G_{c'}[n]$ is finite.

    Now by symmetry of $a$ and $b$, an identical argument yields that we may 
    find $c'' \equiv_{\acl^0(a)} c$ with $c'' \dind_a c$,
    and $H''$ over $\acl^0(c'')$ and an isomorphism $\theta_a'' : G \to H''$ 
    over $\acl^0(cc'')$ with $\theta_a''(a_c) \approx_{cc''} a_{H''} 
    \sim^0_{c''} a$.

    Now $c' \dind c$, so $c' \ind^0 c$ since $\tp(c)$ is wgp, and similarly 
    $c'' \ind^0 c$.
    Also $c' \ind^0 ac$ (since $c'\dind ac$ by (\ref{eq:dcind}) and $\tp(c')$ being wgp), and similarly $c'' \ind^0 ac$. But $c' \equiv_{\acl^0(\emptyset)} c 
    \equiv_{\acl^0(\emptyset)} c''$,
    so
    $c' \equiv^0_{\acl^0(ac)} c''$ by stationarity of $\tp^0(c/\acl^0(\emptyset))$.

    So we obtain $H'$ over $\acl^0(c')$ and an isomorphism $\theta_a' : G \to 
    H'$ over $\acl^0(cc')$ with $\theta_a'(a_c) \approx_{cc'} a_{H'} 
    \sim^0_{c'} a$.

    Now $\theta_b \circ (\theta_a')^{-1}$ is an isomorphism $H' \cong H$ over 
    $\acl^0(cc')$, so by elementarity of the field embedding $\acl^0(c') 
    \preceq \acl^0(cc')$ there is an isomorphism $\phi : H' \cong H$ over 
    $\acl^0(c')$.
    Setting $\theta_a := \phi \circ \theta_a' : G \cong H$ and $a_H := 
    \phi(a_{H'})$, we conclude.
    \subqed{c:isoms}
  \end{proof}

  Now let $\sigma := \theta_b \circ \theta_a^{-1} \in \Aut(H)$.
  Then
  $$\sigma a_H + \theta_b(d_c) \approx_{cc'} \theta_b(a_c) + \theta_b(d_c) = 
  \theta_b(b_c) \approx_{cc'} b_H,$$
  so $h := b_H - \sigma a_H \in \acl^0(dcc')$.

  By Fact~\ref{f:aut},
  $\sigma = \sigma'\sigma_0$ where $\sigma_0$ is over $\acl^0(c')$ and 
  $\sigma' \in \Aut^0(H)$. Let $a_H':=\sigma_0(a_H)$, note that $a_H\sim^0_{c'}a_H'$.
  Let $\alpha := (h,\sigma') \in H \rtimes \Aut^0(H)$.
  Then $\alpha \in \acl^0(dcc')$ and $\alpha*a_H' = b_H$.
  Since the action is faithful, $\cb^0_{c'}(ab/\alpha) \sim^0_{c'} \alpha$. We claim that $\alpha \sim^0_{c'} d$. Note that $a\ind^0 dcc'$ and $b\ind^0 dcc'$ (by Claim~\ref{c:cc'inde}), and $\cb^0(ab/dcc')\sim^0 d$ by $ab\ind^0_d cc'$ (as $a\ind^0_d cc'$ and $b\in\acl^0(ad)$) and $\cb^0(ab/d)\sim^0 d$. Since $\alpha\in \acl^0(dcc')$ and $a\sim^0_{\alpha c'}b$, we get $d\in\acl^0(\alpha c')$ by Fact~\ref{f:cb}. On the other hand $a\ind^0_{c'}\alpha$ (by $a\ind^0\alpha c'$ as $\alpha\in\acl^0(dcc')$), $b\ind^0_{c'}\alpha$ and $\cb^0_{c'}(ab/\alpha) \sim^0_{c'} \alpha$. Since $d\in\acl^0(\alpha c')$ and $a\sim^0_{dc'} b$ we have $\alpha\in\acl^0(dc')$ as required.

  Note that $a_H'$ is generic in $H$ over $\acl^0(c')$, that $\tp(\alpha/\acl^0(c'))$ is wgp because $\tp(d/\acl^0(c'))$ is, that similarly $\tp(a_H'/\acl^0(c'))$ is cgp, and that $a_H'\dind_{c'}\alpha\dind_{c'}b_H$ as $a\dind_{c'}d\dind_{c'}b$ and $a\sim^0_{c'}a_H\sim^0_{c'}a_H'$, $d\sim^0_{c'}\alpha$, $b_H\sim^0_{c'}b$. By Corollary~\ref{c:cohActAbSub}, $\alpha$ is generic in a left coset over 
  $\acl^0(c')$ of an abelian subgroup $G'$ of $H\rtimes\Aut^0(H)$,
  so say $\alpha_0 \in (H\rtimes\Aut^0(H))(\acl^0(c'))$ and 
  $\alpha_0\cdot\alpha \in G'$. Then $(\alpha_0\cdot\alpha)*a_H' = 
  \alpha_0*b_H$, and $a_H' \sim^0_{c'} a$ and $\alpha_0\cdot\alpha \sim^0_{c'} 
  d$ and $\alpha_0*b_H \sim^0_{c'} b$,
  and so the triangle $(a,b,d)$ is abelian over $c'$. Note that $c'\ind^0 abd$ by (\ref{eq:dcind}) and $\tp(abd)$ being wgp.
  We conclude by Lemma~\ref{l:coheir}.
\end{proof}

\section{Asymmetric Elekes-Szabó}
In this section, we translate our main results from the pseudofinite setting to obtain elementary finitary statements generalising the Elekes-Szabó theorem.

First, we recall some definitions from \cite{BB-cohMod}. Bear in mind our convention that ``subvariety'' means ``closed subvariety''.
\label{s:fin}
\begin{definition}
  A \defn{generically finite algebraic correspondence} between 
  irreducible algebraic varieties $V$ and $V'$ is an irreducible 
  subvariety of the product $Z \subseteq V \times V'$ such that the 
  projections $\pi_V(Z) \subseteq V$ and $\pi_{V'}(Z) \subseteq V'$ are 
  Zariski dense, and $\dim(Z) = \dim(V) = \dim(V')$.

  Note that if $V$, $V'$, and $Z$ are defined over $\acl^0(\emptyset)$, the
  condition on $Z$ is equivalent to being of the form $Z=\loc^0((v,v'))$ for
  some generic $v \in V$ and $v'\in V'$ with $v \sim^0 v'$.
\end{definition}

\begin{definition}\label{d:corr}
  If $W_1,W_2,W_3$ are irreducible algebraic varieties over an algebraically closed field,
  we say an irreducible subvariety $U \subseteq W_1 \times W_2 \times W_3$ is in \defn{co-ordinatewise correspondence} with the graph $\Gamma_+ = \{(x,y,z) \in G^3 : x+y=z\}$ of addition in a commutative algebraic group $G$ if there exists a generically finite algebraic correspondence $Z \subseteq U \times \Gamma_+$ such that the Zariski closures of the projections $(\pi_i\times\pi'_i)(Z) \subseteq W_i\times G$ are generically finite algebraic correspondences between $W_i$ and $G$ for $i=1,2,3$. (Technically, this is a property of the tuple $(U,W_1,W_2,W_3,G)$, so we always write explicitly $U \subseteq W_1\times W_2\times W_3$.)

  If $U,W_1,W_2,W_3,(G;+)$ are defined over $\acl^0(\emptyset)$, this is equivalent to there existing generic elements $(w_1,w_2,w_3) \in U$ and $(x,y) \in G^2$ such that $w_1\sim^0x$ and $w_2\sim^0y$ and $w_3\sim^0x+y$, and $w_i$ is generic in $W_i$.
\end{definition}

Let $V$, $W_1$, and $W_2$ be irreducible complex varieties,
and let $U \subseteq W_1 \times V \times W_2$ be an irreducible subvariety.

\begin{definition}
  Let $D\subseteq V$ be a finite subset and $\eps>0$. We call $D$ \textbf{$\eps$-abelian} in $U$ if there are irreducible complex subvarieties $V'\subseteq V$ and $U'\subseteq U$ such that $U'\subseteq W_1\times V'\times W_2$ is in co-ordinatewise correspondence with the graph of addition of a commutative algebraic group $G$, and $|D\cap V'|\geq |D|^{1-\eps}$.\footnote{Note that if $v\in V'$ and the fibre $U_v$ is a generically finite algebraic correspondence between $W_1$ and $W_2$, we have $U'_v = U_v$. But in general, $U'$ could be a proper subset of $U \cap W_1\times V' \times W_2$.}
\end{definition}

In these terms, our main result Theorem~\ref{t:main} translates to the 
following finitary statement:

\begin{theorem}\label{t:main_fin}
Let $\tilde{V}$ be a constructible subset of $V$.
Suppose for all $v\in \tilde{V}$,
  $U_v := \{ (x,y) : (x,v,y) \in U \} \subseteq W_1\times W_2$ is a 
  generically finite algebraic correspondence,
  and $U_{v'} = U_v$ for only finitely many $v' \in V$.

Then for all $\eps>0$  there exists $\eta > 0$
such that the following holds.

Let $D\subseteq \tilde{V}$ be a finite subset which is not $\eps$-abelian in $U$,
and let $A \subseteq W_1$ and $B \subseteq W_2$ be finite subsets such that:
  \begin{enumerate}[(i)]
    \item $|A|=|B| \geq \frac1\eta$.
    \item $|A|^{\eps} \leq |D| \leq |A|^{\frac1\eps}$.
    \item If $W_1' \subsetneq W_1$ is a proper irreducible complex subvariety of 
      complexity\footnote{
      The notion of complexity can be any which takes values in $\N$ and is such that for any $\tau$, the family of subvarieties of $W_1$ of complexity $\leq\tau$ is a constructible family; see \cite[\S2.1.10]{BB-cohMod} for an explicit definition.} $\leq\frac1\eta$ then
      $|A \cap W_1'| \leq |A|^\eta$.
      Similarly for $(B,W_2)$ in place of $(A,W_1)$.
  \end{enumerate}

  Then $|U \cap (A\times D \times B)| \leq (|A| |D|)^{1-\eta}$.

\end{theorem}
\begin{proof}
Let $\eps>0$, suppose no such $\eta>0$ exists.
For $\eta=\frac{1}n$, let $D_n,A_n,B_n$ be a counterexample.
Let $D$ be the ultraproduct $D= \prod_{n\to\U} D_n \subseteq V(K)$, where $K = \C^\U$, and define $A$ and $B$ similarly.

Then $|A| \in \N^\U \setminus \N$ by (i), so let $\bdl := \bdl_{|A|}$.

We have $\bdl(U\cap (A\times D\times B))\geq\bdl(A)+\bdl(D)$ (since $\bdl(U\cap (A\times D\times B))\geq(1-\frac1n)(\bdl(A)+\bdl(D))$ for all $n\in\N_{>0}$ by assumption), and
$\bdl(A)=\bdl(B)$. Note that $0 < \bdl(D) < \infty$ by (ii).

Work in an adequate expansion of $(K;+,\cdot)$ satisfying Assumption~\ref{a:aclbase} in which $A,B,D$ are $\emptyset$-definable and $W_1,V,W_2,U,\tilde{V}$ are defined over $\acl^0(\emptyset)$, which exists by Lemma~\ref{l:adequacy}.

Find $(a,d,b)\in U\cap (A\times D\times B)$ such that $\bdl(adb)=\bdl(U\cap (A\times D\times B))$.
If $\bdl(a/d)>0$ then
it follows from (iii) that $\tp(a)$ is cgp, and $\loc^0(a/\acl^0(d)) = W_1$; indeed,
if $W_1'=\prod_{m\to\U}(W_1')_m$ is a proper irreducible subvariety of $W_1$ over $K$, then by (iii), $\{ m : |A_m \cap (W_1')_m| \leq |A_m|^{\frac1m}\} \in \U$, and so $\bdl(A\cap W_1')=0$. This shows that $\tp(a)$ is cgp, and then from $\bdl(a/d)>0$ we conclude $\loc^0(a/\acl^0(d)) = W_1$.
Then $a\sim^0_d b$ since $U_d$ is a generically finite algebraic correspondence,
so $\bdl(b/d)=\bdl(a/d)>0$.

By symmetry, we obtain that $\bdl(a/d) > 0$ if and only if $\bdl(b/d) > 0$.
If $\bdl(a/d) = 0 = \bdl(b/d)$ then $\bdl(d) = \bdl(adb) \geq \bdl(A) + \bdl(D) \geq \bdl(A) + \bdl(d)$, contradicting $\bdl(A) > 0$.
So we conclude that $\tp(a)$ and $\tp(b)$ are cgp and broad, and $a \sim^0_d b$.

Now since
$\bdl(a)\leq\bdl(A)$ and 
$\bdl(d)\leq\bdl(D)$,
we have $\bdl(A)+\bdl(D)\leq \bdl(U\cap(A\times D \times B)) = \bdl(abd) = \bdl(ad) \leq \bdl(a)+\bdl(d)\leq \bdl(A)+\bdl(D)$,
so $a\dind d$, and $\bdl(d) = \bdl(D) > 0$.
Similarly, $b \dind d$.

We may assume $\tp(d)$ is wgp by Lemma~\ref{l:wgpification}(i); this may increase $\acl^0(\emptyset)$, but we still have $\loc^0(a/\acl^0(d))=W_1$ by (iii) since $\bdl(a/d)>0$.
Hence also $\loc^0(a)=W_1$, and similarly $\loc^0(b)=W_2$.

Now $d$ is interalgebraic with $cb^0(ab/d)$ by the assumption that $U_{v'} = U_d$ for only finitely many $v' \in V$.

Let $V'=loc^0(d)$. Then $\bdl(D\cap V') \geq \bdl(d)$.

By Theorem~\ref{t:main}, $(a,d,b)$ is abelian. By the remark in the last paragraph of Definition~\ref{d:corr}, there is a commutative algebraic group $G$ such that $\loc^0(a,d,b)\subseteq W_1\times V'\times W_2$
is in co-ordinatewise correspondence with the graph of addition of $G$.

Now $U' := \loc^0(a,d,b)$ is a subset of $U\cap (W_1\times V'\times W_2)$. But no $D_n$ is $\eps$-abelian in $U$, so we must have $\bdl(D\cap V')\leq (1-\eps)\bdl(D)$, contradicting $\bdl(D\cap V')\geq\bdl(d)=\bdl(D)$.
\end{proof}

We deduce a somewhat simpler statement by introducing additional algebraic 
assumptions:

\begin{corollary}\label{c:main_fin_sym}
  Suppose now $\dim(W_1) = \dim(W_2) = \dim(V) =: k$, and
  $\dim(U) = 2k$ and $U$ projects dominantly to each pair $W_1\times V$, $V\times W_2$, and $W_1 \times W_2$, and suppose $U$ is not in co-ordinatewise correspondence with the graph of addition of a commutative algebraic group.

Then for all $\eps>0$  there exists $\eta > 0$ such that
if $A \subseteq W_1$ and $B\subseteq W_2$ and $D \subseteq V$ are finite 
subsets such that for each $d \in D$, the fibre $U_d$ is a generically finite 
algebraic correspondence between $W_1$ and $W_2$,
and $U_d = U_v$ holds for only finitely many $v \in V$,
and (i)-(iii) of Theorem~\ref{t:main_fin} hold,
then $|U \cap (A\times D \times B)| \leq (|A| |D|)^{1-\eta}$.
\end{corollary}
\begin{proof}
  Let $\tilde V$ be the set of $v \in V$ for which $U_v$ is a generically finite algebraic correspondence between $W_1$ and $W_2$, and $U_v = U_{v'}$ holds for only finitely many $v' \in V$.
  Then $\tilde V$ is constructible. Indeed, the finiteness condition is definable by uniform finiteness, and the correspondence condition holds of $U_v$ if and only if it is irreducible and $\dim(U_v) = \dim(\pi_{W_1}(U_v)) = \dim(\pi_{W_2}(U_v)) = k$.

  The result then follows from Theorem~\ref{t:main_fin} once we observe that $D$ is not $\eps$-abelian. But indeed, suppose $U'$ and $V'$ witness that $D$ is $\eps$-abelian. Then $\dim U' = 2\dim(W_1) = \dim(U)$, so $U' = U$ (and $V' = V$), contradicting the assumption that $U$ is not in co-ordinatewise correspondence with the graph of addition of a commutative algebraic group.
\end{proof}

Specialising further, we obtain an unbalanced generalisation of the 1-dimensional 
case of the main theorem of \cite{ES-groups}, and in this case we can bound 
the exponent using the results of Appendix~\ref{s:eps}:
\begin{theorem}\label{t:main_fin_1d}
  Let $U \subseteq W_1\times V \times W_2$ be as in 
  Corollary~\ref{c:main_fin_sym} with $k=1$.

  Let $c$ be the constant from Fact~\ref{fact-JZ}.
  Replacing $c$ with $\min(1,c)$ if necessary, we may assume $c \leq 1$.
  Let $\eps_0 := \eps2^{-(\frac4{c\eps}+7)}$
  and let $\eta_0 := \frac{\eps_0}{1+\frac1\eps}$.

  Then for all $1>\eps>0$ there exists $\eta > 0$ such that the following holds.
  If $A \subseteq W_1$, $B\subseteq W_2$, and $D \subseteq V$ are finite 
  subsets satisfying conditions
  (i)-(ii) of Theorem~\ref{t:main_fin},
  then $|U \cap (A\times D \times B)| \leq (|A| |D|)^{1-\eta_0}$.
\end{theorem}
\begin{proof}
  Otherwise, as in the proof of Theorem~\ref{t:main_fin},
  we obtain $(a,d,b) \in A\times D \times B$
  with $\bdl(D) = 1 = \dim(V)$ (by setting $\bdl=\bdl_{|D|}$)
  and $\bdl(A)=\bdl(B) \in [\eps,\frac1\eps]$
  and $\bdl(adb) \geq (1-\eta_0)(\bdl(A)+1)$.

  Note $\eta_0 \leq \eps2^{-7} \leq \frac\eps4 \leq \frac\eps{2\eps+2}$,
  so
  $\bdl(adb) \geq (1-\frac\eps{2(\eps+1)})(\eps+1) = \eps+1-\frac\eps2 > 1$.

  We may also assume that $D$ omits the finitely many elements $v\in V$ for 
  which $U_v$ is not a correspondence.
  Note that for all $v\in V$, $U_v = U_{v'}$ for only finitely many $v'$; indeed, otherwise it would hold for all but finitely many $v'$, and then the projection to $W_1\times W_2$ of $U$ would be a finite union of fibres $U_{v_i}$ so has dimension $\leq 1$, contrary to assumption.
  We do not have $\bdl(a/d)=0=\bdl(b/d)$, since $\bdl(adb) > 1 \geq \bdl(d)$.
  So (as in Theorem~\ref{t:main_fin})
  $\dim^0(a/d)=1=\dim^0(b/d)$ and $a\sim^0_d b$, and $d \sim^0 \cb^0(ab/d)$.

  We proceed to show that $(a,d,b)$ is an $(\eps,1)$-correspondence triangle (see Definition~\ref{d:epsTri}).

  Now $\bdl(a/d) = \bdl(adb) - \bdl(d)
  \geq (1-\eta_0)(\bdl(A)+1) - \bdl(d)
  \geq \bdl(a) - \eta_0(\frac1\eps+1) + 1 - \bdl(d)
  $,
  so
  $\frac{\bdl(a) - \bdl(a/d)}{\bdl(d)}
  \leq 1-\frac{1-\eta_0(\frac1\eps+1)}{\bdl(d)} 
  = 1-\frac{1-\eps_0}{\bdl(d)} 
  \leq \eps_0$
  (using $\eta_0(\frac1\eps+1)=\eps_0$ and $\bdl(d)\leq1$).

  Note that
  \begin{align*}
  \bdl(d) &= \bdl(adb) - \bdl(a/d)
  \geq \bdl(adb) - \bdl(A)\\
  &\geq (1-\eta_0)(\bdl(A)+1)-\bdl(A)=1-\eta_0(\bdl(A)+1)\\
  &\geq 1-\eta_0(\frac1\eps+1)
  = 1- \eps_0 > \frac{1}{2},
  \end{align*}

  and

  \begin{align*}
  \bdl(a) &= \bdl(ad) - \bdl(d/a)
  \geq \bdl(adb) - \bdl(D)\\
  &\geq (1-\eta_0)(\bdl(A)+1)-\bdl(D)\geq \bdl(A)-\eta_0(\bdl(A)+1)\\
  &\geq \bdl(A)-\eta_0(\frac1\eps+1)
  = \bdl(A)- \eps_0 > \frac{\eps}{2}  ~(\text{since } \eps_0< \frac{\eps}{2}).
  \end{align*}

  Thus, $\frac{1}{2}<\bdl(d)\leq\bdl(D)=1$ and $\frac{\eps}{2}<\bdl(a)\leq\bdl(A)\leq\frac{1}{\eps}$ and  $\frac{\bdl(a)}{\bdl(d)}\in(\frac{\eps}{2},\frac{2}{\eps})$.

  Rescaling by working with $\bdl_{\xi'}$ instead of $\bdl=\bdl_{\xi}=\bdl_{|D|}$ where $\bdl_{\xi'}(d)=1$ (and hence $\bdl_{\xi'}(t)=\frac{1}{\bdl(d)}\bdl(t)$ for all $t\in K^n$), we may assume $\bdl(d) = 1$. Then $\bdl(a)-\bdl(a/d) \leq \eps_0$ 
  and $\bdl(a) \in (\frac{\eps}{2},\frac{2}{\eps})$.

  Let $i_0 := 1 + \lfloor\bdl(a)\rfloor \leq 1+\frac2\eps$ and note
  that $\bdl(a)2^{-(\frac{2\bdl(a)}c+7)} > \frac{\eps}{2}\cdot 2^{-(\frac4{c\eps}+7)} \geq \eps_0$.
  We also have
  $\eps_0 \leq 2^{-(\frac4{c\eps}+8)}\leq 2^{-(2i_0+4)}\leq \min(2^{-(i_0+1)},\frac1{20})$ (using $c\leq 1$).

  Then Proposition~\ref{p:epsmaincgp} case (1), with $\gamma=0$ and $\eps'=1$, 
  applies to show that $(a,d,b)$ is abelian, leading to a contradiction as in 
  Corollary~\ref{c:main_fin_sym}.
\end{proof}

\section{Expansion under polynomial families}
\label{s:ER}
In this section, we prove our generalisation of the Elekes-Rónyai theorem.

We write $\AGL_1$ for the 2-dimensional algebraic group of affine linear transformations of $\A^1$, which we identify with the semidirect product $\AGL_1 = \G_a \rtimes \G_m$ where $\G_m$ acts on $\G_a$ by multiplication.
\begin{fact}\label{f:AGLsubgroups}
  The proper non-trivial connected algebraic subgroups $H \leq \AGL_1$ are precisely $\G_a \leq \G_a \rtimes \G_m$ and the conjugates within $\AGL_1$ of $\G_m \leq \G_a \rtimes \G_m$.
\end{fact}
\begin{proof}
  The image of $H$ under the homomorphism $\pi_2 : \G_a\rtimes \G_m \to \G_m$ 
  is either trivial,
  in which case $H = \G_a$, or the whole of $\G_m$.
  In the latter case, say $(a,b) \in H$ where $b\in\G_m$ is not a root of unity.
  Setting $\gamma := (\frac a{b-1},1)$, we have $\gamma^{-1}(0,b)\gamma = (a,b)$, so
  $(a,b) \in (\G_m)^\gamma$,
  so $H$ intersects $(\G_m)^\gamma$ on the infinite subgroup generated by $(a,b)$.
  But $\dim(H) = 1$, so then $H=(\G_m)^\gamma$.
\end{proof}

For a field $k$ and a natural number $d$, we write $k[t]_{\leq d}$ for the set
of polynomials over $k$ of degree at most $d$.

\begin{lemma}\label{l:addmul}
  Let $f\in K[t]$, and suppose $(a,f,f(a))$ is a correspondence triangle for some $a \in K$ (where recall we identify $f$ with its tuple of coefficients).
Then there are $h_0,g\in \acl^0(\emptyset)[t]$ and $d \in \acl^0(f)$ such that
$f(x)=h_0(g(x)+_Gd)$, where $+_G$ is either addition or multiplication.
\end{lemma}
\begin{proof}
  Let $C_0$ be the subfield $\acl^0(\emptyset) \leq K$.
  Let $b := f(a)$. The correspondence triangle $(a,f,b)$ is semi-cgp, since $\dim^0(a) = 1$ and so $a$ is cgp by Remark~\ref{r:dim1cgp}.

  By Theorem~\ref{t:main},
  $\dim^0(f) = \dim^0(a) = 1$,
  and taking $f'b' \equiv_a fb$ with
  $f'b' \dind_a fb$ (and hence $f'b' \ind^0_a fb$ via Fact~\ref{f:wgpFacts}(ii)),
  and setting $e := \cb^0(bb'/ff')$,
  we have $\dim^0(e) = 1$ and $f\ind^0 e$.

  The remainder of this proof uses ideas and lemmas from \cite{ER}.

  \newcommand{\hf}{\hat f}
  Let $s := \deg(f) = \deg(f')$.
  Let $g' \in K[x]$ be of maximal degree such that
  $\exists l,l' \in K[x]_{\leq s}\; (f = l\o g' \wedge f' = l'\o g')$.

  \begin{claim}\label{c:descend}
    There exist $g \in C_0[x]$
    and $\hf \in \acl^0(f)[x]$
    and $\hf' \in \acl^0(f')[x]$
    such that $f = \hf \o g$ and $f' = \hf' \o g$
    and $\deg(g) = \deg(g')$.
  \end{claim}
  \begin{proof}
    \newcommand{\gmodsim}{[g']_\sim}

    Consider the definable equivalence relation on $K[x]_{\leq s}$
    defined by $h \sim h'$ if there exists a non-constant linear polynomial $\alpha \in K[x]$ such
    that $h' = \alpha \o h$.

    By \cite[Proposition~9]{ER}, $\{ h \in K[x]_{\leq s} : \exists l \in K[x]_{\leq
    s} \; f = l \o h \}$ is a union of $\leq 2^s$ $\sim$-equivalence classes,
    so $\gmodsim \in \acl^0(f)$.

    Similarly $\gmodsim \in \acl^0(f')$. But $f \ind^0 f'$, so $\gmodsim \in
    \acl^0(f) \cap \acl^0(f') = \acl^0(\emptyset) = C_0$; in other words, the constructible set $\gmodsim$ is defined over $C_0$.
    Since $C_0\leq K$ is an algebraically closed subfield, there exists $g \in \gmodsim\cap C_0[x]_{\leq s}$.

    Finally, we can take $\hf \in \acl^0(f)[x]$ and $\hf' \in \acl^0(f')[x]$
    such that $f = \hf \o g$ and $f' = \hf' \o g$,
    by elementarity of $\acl^0(f)$ resp.\ $\acl^0(f')$ in $K$.
    \subqed{c:descend}
  \end{proof}

  \newcommand{\tf}{\tilde f}
  \newcommand{\thf}{\tilde \hf}
  \newcommand{\tb}{\tilde b}
  \newcommand{\tphi}{\tilde phi}
  Now let $\tf \tf' \thf \thf' \tb \tb' \equiv^0_e ff'\hf\hf'bb'$ with $\tf \tf' \thf \thf' \tb \tb' \ind^0_e ff'\hf\hf'bb'$

  Let $\gamma$ be the Zariski closure in $\A^2$ of
  $\{ (\hf(x), \hf'(x)) : x\in\A^1 \} = (\hf\times \hf')(\A^1)$.
  Then $\gamma = \loc^0(bb'/\acl^0(ff'))$, since $\gamma$ is irreducible,
  and $\dim^0(bb'/ff') = 1 = \dim \gamma$,
  and $(b,b') = (\hf(g(a)), \hf'(g(a))) \in \gamma$.
  So by choice of $e$, also $\gamma = \loc^0(bb'/e)$, and $\gamma$ is defined
  over $\dcl^0(e)$.
  Hence $\gamma$ is also the Zariski closure of $(\thf\times \thf')(\A^1)$.

  So by \cite[Lemma~12]{ER},
  there are polynomials $h,h',\phi_f,\phi_{\tf} \in K[x]$
  such that
  \begin{align*}
    \hf &= h \o \phi_f & \hf' &= h' \o \phi_f \\
    \thf &= h \o \phi_{\tf}.
  \end{align*}

  Then $\phi_f$ is linear by maximality of $\deg(g)$,
  since $f = \hf \o g = h \o \phi_f \o g$
  and $f' = \hf' \o g = h' \o \phi_f \o g$.

  \newcommand{\simo}{\sim_{\operatorname{out}}}
  \newcommand{\hsimo}{[h]_{\simo}}
  We now proceed as in the proof of the claim.
  Note that $f \ind^0 \tf$ (since $f \ind^0 e$), and hence $\hf \ind^0 \thf$.
  Consider the definable equivalence relation on $K[x]_{\leq s}$
  defined by $h \simo h'$ if there exists a non-constant linear polynomial
  $\alpha \in K[x]$ such that $h' = h \o \alpha$.
  By \cite[Proposition~9]{ER} (see the remark at the end of \cite[Definition~8]{ER}),
  $\{ \theta \in K[x]_{\leq s} : \exists l \in K[x]_{\leq
  s} \; \hf = \theta \o l \}$ is a union of $\leq 2^s$ $\simo$-equivalence classes,
  so $\hsimo \in \acl^0(\hf)$, and similarly $\hsimo \in \acl^0(\thf)$,
  so $\hsimo \in \acl^0(\hf) \cap \acl^0(\thf) = C_0$,
  so there is a non-constant linear polynomial $\alpha \in K[x]$
  such that $h_0 := h \o \alpha \in C_0[x]$.

  So we have decomposed $f$ as $f = h_0 \o \alpha^{-1} \o \phi_f \o g$ with $h_0,g \in C_0[x]$ and
  $\alpha^{-1} \o \phi_f$ linear, and so by elementarity of $\acl^0(f)$ in
  $K$, we may find a linear polynomial $\phi_f' \in \acl^0(f)[x]$
  such that $f = h_0 \o \phi_f' \o g$.

  Now $a' := g(a) \sim^0 a$ and $\phi_f' \sim^0 f$ and $b' := \phi_f'(g(a)) \sim^0 b$. So $(a',\phi_f',b')$ is a cgp correspondence triangle,
  and satisfies the conditions of
  Corollary~\ref{c:cohActAbSub}
  for the action of $\AGL_1 = \G_a \rtimes \G_m$ on $\A^1$.
  Hence $\phi_f' \in \AGL_1$ is generic in a left coset $X$ defined over $C_0$ of a connected abelian subgroup $G' \leq \AGL_1$.
  Since $X$ is over $C_0$, also $G'$ is over $C_0$.
  So by Fact~\ref{f:AGLsubgroups}, $G'$ is either $\G_a \leq \G_a \rtimes \G_m$
  or a conjugate by some $\alpha \in \AGL_1(C_0)$ of $\G_m \leq \G_a \rtimes \G_m$.

  Therefore, after composing $h_0$ and $g$ with linear polynomials over
  $C_0$, we obtain $f = h_0 \o \phi'_f \o g$ where $\phi'_f$ is either
  addition or multiplication by a constant.
\end{proof}

For a finite set $F \subseteq \C[t]_{\leq d}$ and $\eps>0$ we call $F$ \textbf{$\eps$-additive}, if there are $g,h \in \C[t]$ such that
  \[| F \cap \{ h(g(t) + a) : a \in \C \} | \geq |F|^{1-\eps}.\]
  Similarly, we call $F$ \textbf{$\eps$-multiplicative}, if there are $g,h \in \C[t]$ such that
  \[| F \cap \{ h(g(t)\cdot a) : a \in \C \} | \geq |F|^{1-\eps}.\]
\begin{theorem}\label{t:higherER}
  For all $\eps>0$, there exists $\eta > 0$
  such that
  for any finite set $F \subseteq \C[t]$ of non-constant polynomials of degree
  $\leq 1/\eps$ which is neither $\eps$-additive nor $\eps$-multiplicative, 
  for any finite set $A \subseteq \C$
 with  $|F|\geq |A|^\eps\geq \frac{1}{\eta}$,
  we have  \[|F*A| = |\{ f(a) : f \in F,\; a \in A\}|\geq |A|^{1+\eta}.\]
\end{theorem}
\begin{proof}
  Let $\eps > 0$, and suppose there is no such $\eta$.
  Let $d:=\lceil\frac{1}{\eps}\rceil$.
  Taking an ultraproduct of counterexamples $(A_n,F_n)$ with $\eta =
  \frac1n$, we obtain infinite pseudofinite sets $A \subseteq K = \C^\U$ and 
  $F\subseteq K[t]_{\leq d}$ such that,
  with $\bdl = \bdl_{|A|}$,
  \begin{itemize}
    \item $\bdl(F) \geq \eps\bdl(A)$
    \item $\bdl(F*A) = \bdl(A)$
  \end{itemize}
  (where $\bdl(F*A) \geq \bdl(A)$ holds since if $f \in F$, we have $|A| \leq d|f*A|$, since $f(x)=b$ has $\leq d$ solutions in $x$ for any $b$),
  and for any $g,h \in K[t]$ with $\deg g + \deg h \leq d$,
  setting $V_{g,h}^+ := \{ h(g(t) + a) : a \in K \} \subseteq K[t]_{\leq d}$
  and $V_{g,h}^\cdot := \{ h(g(t) \cdot a) : a \in K \} \subseteq K[t]_{\leq 
  d}$,
  we have $\bdl(F \cap V_{g,h}^+) \leq (1-\eps)\bdl(F)$
  and $\bdl(F \cap V_{g,h}^\cdot) \leq (1-\eps)\bdl(F)$.
  Note that $\bdl(F)<\infty$. Indeed, take $Y\subseteq A$ of size $d+1$. If $f(y)=g(y)$ for all $y\in Y$, then $f=g$. Therefore, $|(F*Y)^{d+1}|\geq |F|$. Hence, $\bdl(F)\leq (d+1)\bdl(F*A)=(d+1)\bdl(A)$.

  Work in an adequate expansion of $(K;+,\cdot)$ satisfying Assumption~\ref{a:aclbase} in which $A$ and $F$ are $\emptyset$-definable, which exists by Lemma~\ref{l:adequacy}.
  Take $(a,f) \in A\times F$
  with $\bdl(af) = \bdl(A\times F) = \bdl(A) + \bdl(F)$,
  where we identify a polynomial in $K[t]_{\leq d}$ with the element of
  $K^{d+1}$ consisting of its coefficients.
  Then $a \dind f$.
  Set $b := f(a) \in F*A$,
  so (since $f$ is a non-constant polynomial) $a \sim^0_f b$,
  and $f \sim^0 \cb^0(ab/f)$ by the identification of $f$ with its 
  coefficients.
  Then $\bdl(b) \leq \bdl(F*A) = \bdl(A) = \bdl(a) = \bdl(a/f) = \bdl(b/f) \leq \bdl(b)$,
  so $b \dind f$.

  Now by Lemma~\ref{l:wgpification}(i), let $c$ be such that $\bdl(c) = 0$ and 
  $\tp(f/c)$ is wgp. Then $\bdl(ab/fc) = \bdl(ab/f) > 0$ so $\dim^0(ab/fc) > 
  0$, so $\dim^0(ab/fc) = 1 = \dim^0(ab/f)$,
  so $ab \ind^0_f c$.
  Also $a \ind^0 f \ind^0 b$ since $\bdl(a/f),\bdl(b/f)>0$.
  By Lemma~\ref{l:cbbase}, $\cb^0_c(ab/f) \sim^0_c f$.

  Add $\acl^0(c)$ to the language. Then $\tp(f)$ is wgp and broad,
  and $\tp(a)$ is broad and cgp since $\dim^0(a) = 1$ (see Remark~\ref{r:dim1cgp}),
  and by the above we still have $f \sim^0 \cb^0(ab/f)$.

  Therefore, $(a,f,b)$ is a semi-cgp correspondence triangle.

  By Lemma~\ref{l:addmul},
  $f \in V_{g,h_0}^+$ or $f\in V_{g,h_0}^\cdot$.
  But then since $V_{g,h_0}^+$ and $V_{g,h_0}^\cdot$ are over $\acl^0(\emptyset)$,
  we have
  $(1-\eps)\bdl(F) \geq \bdl(F \cap V_{g,h_0}^+) \geq \bdl(f) = \bdl(F)$,
  which is a contradiction.
\end{proof}

\begin{remark}
  The bound on degree in Theorem~\ref{t:higherER} is necessary:
  for example,
  taking $F:=\{x^{2^i}: 0<i\leq n\}$ and $A:=\{2^{2^i}: i\leq n\}$, we have
  $F*A= \{2^{2^j}:0<j\leq 2n\}$.
\end{remark}

\begin{remark}
  We leave it to the interested reader to verify a version of
  Theorem~\ref{t:higherER} in the case of rational functions $F \subseteq
  \C(x)$, where $\PSL_2$ will appear in place of $\AGL_1$ in the proof.
\end{remark}

We now proceed to deduce our arbitrarily unbalanced version of Elekes-Rónyai.

\begin{definition}
A polynomial $f \in \C[x,\bar y]$
is \defn{additive} if there
exist $g,h\in \C[x]$ and $s\in \C[\bar y]$ such that $f(x,\bar y) = g(h(x)+s(\bar y))$, and 
\defn{multiplicative} if there
exist $g,h \in \C[x]$ and $s\in \C[\bar y]$ such that $f(x,\bar y) = g(h(x)\cdot s(\bar y))$.
\end{definition}

\begin{theorem}\label{t:unbalancedERMany}
  Suppose $f\in\C[x,\y]$ is neither additive nor multiplicative. Then for any $\eps>0$, there are $C,\eta>0$ such that for any finite $A, B_0,\ldots,B_{|\y|-1}\subseteq \C$ with $|A|\geq|B_i|\geq |A|^\eps\geq C$ for each $i<|\y|$, we have $|f(A,B_0,\ldots,B_{|\y|-1})|\geq |A|^{1+\eta}$.

  Moreover, one can take $\eta:=|\y|^{-2}\cdot 2^{-\frac{c'}\eps}$ for some absolute constant $c'$; in particular, $\eta$ depends only on $|\y|$ and $\eps$.
\end{theorem}
\begin{proof}
  Let $f\in\C[x,\y]$ and $\eps>0$ be given.
  We first prove the qualitative statement that $C,\eta>0$ exist, then consider how to adapt the argument to obtain $\eta$.

  Let $m := |\y|$.
  Suppose no such $C,\eta>0$ exist.
  By taking an ultraproduct of counter-examples of growing sizes for decreasing $\eta$, we find internal subsets $A,B_0,\ldots,B_{m-1}\subseteq K$ such that for all $n\in \N$, $|f(A,B_0,\ldots,B_{m-1})| < |A|^{1+\frac1n}$ and $|A|\geq|B_i|\geq |A|^\eps\geq n$ for all $i<m$.
  Let $\bdl:=\bdl_{\max\{|B_i|:i<m\}}$.
  Take $(a,\b)\in A\times B_0\times\ldots\times B_{m-1}$ such that $\bdl(a\b)=\bdl(A\times B_0\times\ldots\times B_{m-1})$.
  Then $a\dind\b$, and $\b = (b_i)_{i<m}$ is a $\dind$-independent sequence.
  Note that $\max\{\bdl(b_i):i<m\}=1$ and $1\leq \bdl(a)\leq \frac{1}{\eps}$ and $\bdl(b_i)\geq\eps\bdl(a)$.
  In particular, $\bdl(a),\bdl(b_i)>0$. Hence $(a,\b)$ is a generic tuple in $K^{m+1}$ over $C_0:=\acl^0(\emptyset)$.

  Let $f_\b:=f(x,\b)\in K[x]$ and $b:=f_{\bar{b}}(a)$.
  Then $f_\b\in\acl^0(\b)$ (identifying $f_\b$ with the sequence of its coefficients), and $a\dind f_\b$.
  Note that $f_\b$ is non-constant, since $f(x,\bar y)$ is not additive and so depends on $x$, and $\b$ is generic in $K^m$.
  So $b \sim^0_{f_\b} a$.

  Now $\bdl(b) \leq \bdl(a)$, so $\bdl(b/f_\b) = \bdl(a/f_\b) = \bdl(a) \geq \bdl(b)$ and so $b \dind f_\b$.
  Since $\b$ is $\dind$-independent and each $b_i$ is cgp (by $\dim^0(b_i)=1$ and Remark~\ref{r:dim1cgp}), $\b$ and hence $f_\b$ are wgp by Fact~\ref{f:wgpFacts} and Remark~\ref{r:tuplewgp}.
  So $(a,f_\b,b)$ is a semi-cgp correspondence triangle,
  hence by Lemma~\ref{l:addmul}, there are $h_0,g\in C_0[t]$ and $d \in \acl^0(\b)$
  such that $f_\b(x)=h_0(g(x)+_Gd)$, where $+_G$ is either addition or multiplication.

  Now we have $a\sim^0 g(a)$, $b\sim^0 g(a)+_G d$ and $b\in\dcl^0(a,\b)$.
  We follow the argument of Remark~\ref{r:dcl} and claim that the minimal polynomial $m(t)$ of $d$ over $C_0(a,\b)$ is $t^k - d^k$ for some $k \in \N$.
  Indeed, let $\sigma \in \Gal(C_0(a,\b))$, and let $\xi := \sigma(d) -_G d \in \acl^0(\b)$.
  Let $g_b := g(a)+_G d \in \acl^0(b)$.
  Then $\sigma(g_b) = \sigma(g(a)+_G d) = g(a)+_G\sigma(d) = g_b+_G\xi$, and since $b \in C_0(a,\b)$ and $g_b \in \acl^0(b)$, we have $\sigma(g_b) \in \acl^0(b)$ and so $\xi = \sigma(g_b) -_G g_b \in \acl^0(b)$.
  But $b \ind^0 \b$; indeed, since $f_\b$ is non-constant we have $a \sim^0_\b b$, and so $\dim^0(b/\b) = \dim^0(a/\b) = 1 = \dim^0(b)$.
  So $\xi \in \acl^0(b) \cap \acl^0(\b) = \acl^0(\emptyset) = C_0$.
  Hence the roots of $m(t)$ are invariant under $(x \mapsto x+_G\xi)$.
  Considering arbitrary $\sigma$, we conclude that the roots of $m(t)$ are of the form $d+_G \Xi$ for the finite subgroup $\Xi = \{ \sigma(d) -_G d : \sigma \in \Gal(C_0(a,\b))\}$.
  If $+_G$ is addition, then $\Xi$ is trivial as $\operatorname{char}(K)=0$, so $m(t) = t-d$ as required. If $+_G$ is multiplication, then the roots are of the form $\{ d\cdot \zeta^i : i<k\}$ where $\zeta$ is a primitive $k$th root of unity, and then $m(t)=t^k-d^k$ as required.

  Continue to assume $+_G$ is multiplication.
  Say $h_0(t)=\sum_{i\leq n}a_it^i$.
  Then we can write $b=h_0(g(a)\cdot d)$ as
  \[b=\sum_{i\leq n}a_ig(a)^id^{i-\lfloor\frac{i}{k}\rfloor k}\cdot d^{\lfloor\frac{i}{k}\rfloor k}. \]
  Let $d_i:=d^{\lfloor\frac{i}{k}\rfloor k}\in C_0(\bar{b})$.
  Then \[\sum_{i\leq n}a_ig(a)^id_it^{i-\lfloor\frac{i}{k}\rfloor k}-b=0\] is a polynomial in $C_0(a,\b)$ of degree $<k$ satisfied by $d$.
  Thus, it is the constant zero polynomial, and $a_i=0$ for all $i$ not divisible by $k$ (as $g(a)\neq 0)$.
  Therefore, $h_0(t)=\hat{h}_0(t^k)$ for some $\hat{h}_0\in C_0[t]$ and $f(a,\b)=\hat{h}_0(g(a)^k\cdot d^k)$.

  Meanwhile in the additive case, we have $k=1$ and $f(a,\b) = h_0(g(a)+d)$.
  But $d^k \in C_0(a,\b) \cap \acl^0(\b) = C_0(\b)$ and $(a,\b)$ is generic,
  so in both cases we have an equality of rational functions $f(x,\y) = h'(g'(x)+_G s(\y))$ where $h'$ and $g'$ are non-constant polynomials over $C_0$ and $s(\y) \in C_0(\y)$.
  But then since $f$ and $h'$ are polynomials, so must be $g'(x)+_G s(\y)$ (as can be seen by considering $C_0(x,\y)$ as the fraction field of the unique factorisation domain $C_0[x,\y]$), so also $s \in C_0[\y]$.
  But this contradicts the assumption that $f$ is neither additive nor multiplicative.

  This concludes the qualitative argument. We now consider the ``moreover'' clause in the statement. So let $\eta := m^{-2}\cdot 2^{-\frac{c'}\eps}$, and suppose for a contradiction that we have $a,\b,b$ as above, but now with $\bdl(b) \leq (1+\eta)\bdl(a)$ rather than $\bdl(b)\leq \bdl(a)$.

  We apply Lemma~\ref{l:epsaddmul} from Appendix~\ref{s:eps}, using some notation defined in that appendix.

  First, since $\bdl(b) \leq \bdl(a) + \eta\bdl(a)$,
  we have $\bdl(b/f_\b)=\bdl(a/f_\b)=\bdl(a)\geq\bdl(b)-\eta\bdl(a)$.
  Also, since $\b$ is a $\bdl$-independent sequence with $\bdl(b_i)\geq \eps\bdl(a)$, it is easy to see that $\tp(\b)$ is $\frac{\eps\bdl(a)}{2}$-wgp.
  By Lemma~\ref{l:epswgpFacts}(i), $f_\b$ is also $\frac{\eps\bdl(a)}2$-wgp.

  Thus, $(a,f_\b,b)$ is a semi-cgp $(\eta\bdl(a),\frac{\eps\bdl(a)}{2})$-correspondence triangle over $C_0:=\acl^0(\emptyset)$.
  By Lemma~\ref{l:epsaddmul}, there are $h_0,g\in C_0[t]$ and an absolute constant $c$ such that $f_\b=h_0\circ\phi_f'\circ g$ where $\phi_f'$ is either addition or multiplication, provided $\eta\bdl(a)<\{\frac{\eps\bdl(a)}{2^{i_0+5}m},\frac{\bdl(a)}{2^{\frac{2}{c\eps}+7}}\}$ where $i_0:=1+\lfloor\frac{2}{\eps}\rfloor$.
  Choosing an appropriate absolute constant $c'$ calculated in terms of the absolute constant $c$, we get that $\eta\bdl(a)$ does satisfy this bound.
\end{proof}

\appendix
\section{Explicit powersaving}
\label{s:eps}

In this appendix, we adapt some of our qualitative arguments to obtain bounds on the exponents. The goal of this section is to establish our quantitative asymmetric Elekes–Rónyai Theorem (Theorem~\ref{t:unbalancedERMany}).
 This involves quantitatively weakening our independence and general position assumptions, and tracing how these errors compound through the proofs.
This method works in general. However, some of these results have more restrictive assumptions than the qualitative counterparts, primarily because we do not have a quantitative version of Theorem~\ref{t:cohActAb} for actions on varieties of dimension $>1$.

\newcommand{\dinde}[1]{\ind^{\bdl,#1}}
\newcommand{\ndinde}[1]{\nind^{\bdl,#1}}
\subsection{$\eps$-Correspondence triangles}
We first give a quantitative analogue of the notion of correspondence triangles. To do so, we weaken the assumption of \(\bdl\)-independence and strengthen the assumption on wgp; these correspond respectively to the expansion constant $\eta$ and the lower bound $\epsilon$ in the statement of Theorem~\ref{t:unbalancedERMany}.
\begin{definition}
  Suppose $\bdl(a/C),\bdl(b/C) < \infty$.
  Let $\eps \geq 0$.
  \begin{itemize}
    \item
      Define $a\dinde\eps_Cb$ to mean: $\bdl(ab/C) + \eps \geq \bdl(a/C) + 
      \bdl(b/C)$.
      Equivalently, $\bdl(a/C) - \bdl(a/bC) = \bdl(b/C) - \bdl(b/aC) \leq \eps$.
    \item Say $\tp(a/C)$ is \defn{$\eps$-wgp} if for any tuple $b$ such that $a \nind^0_C b$,
      we have $a \ndinde\eps_C b$.
    \item Say $\tp(a/C)$ is \defn{$\eps$-cgp} if for any tuple $b$ such that 
      $a \nind^0_C b$, we have $\bdl(a/bC) \leq \eps$.
    \item Say $\tp(a/C)$ is \defn{$\eps$-broad} if $\bdl(a/C)>\eps$.
  \end{itemize}
\end{definition}

Note that $\dinde\eps$ is symmetric, and we have the following versions of transitivity and monotonicity:
$a\dinde\eps bc$ implies $a\dinde\eps b$ and $a\dinde\eps_b c$,
and conversely $a\dinde\eps b$ and $a\dinde{\eps'}_b c$ imply $a\dinde{\eps+\eps'}bc$.

Note also that $a\dinde\eps Cb$ implies $a\dinde{\eps’} Cb$ whenever $0\leq \eps\leq \eps’$. Moreover, if $\tp(a/C)$ is $\eps’$-wgp (respectively, $\eps’$-broad), then it is also $\eps$-wgp (respectively, $\eps$-broad), while if $\tp(a/C)$ is $\eps$-cgp, then it is also $\eps’$-cgp. In particular, if $\tp(a/C)$ is $\eps$-wgp for some $\eps\geq 0$, then it is wgp.

The following Lemma generalises Fact~\ref{f:wgpFacts}. Here, $\eps$-wgp can be preserved, but sometimes at the cost of replacing \(\eps\) by a smaller value.
\begin{lemma}\label{l:epswgpFacts}
  Suppose $\bdl(a/C),\bdl(b/C) < \infty$.
  \begin{enumerate}[(i)]
    \item If $\tp(a/C)$ is $\eps$-wgp and $d \in \acl^0(aC)$, then $\tp(d/C)$ is $\eps$-wgp.
    \item If $\tp(b/C)$ and $\tp(a/Cb)$ are $\eps$-wgp, then so is $\tp(ab/C)$.
    \item If $a \dinde\eps_C b$ and $\tp(a/C)$ is $\eps'$-wgp with $\eps'\geq\eps$, then $\tp(a/Cb)$ is $(\eps'-\eps)$-wgp.
  \end{enumerate}
\end{lemma}
\begin{proof}
  \begin{enumerate}[(i)]
    \item
  Suppose $d \nind ^0_C e$.
  Let $a' \equiv _{Cd} a$ with $a' \dind_{Cd} e$.
  Then $a' \nind^0_C e$ since $d \in \acl^0(a'C)$,
  and $\tp(a'/C) = \tp(a/C)$ is $\eps$-wgp,
  so $a' \ndinde\eps_C e$,
  and so $d \ndinde\eps_C e$;
  indeed, if $d \dinde\eps_C e$, then since also $a'\dind_{Cd} e$,
  by transitivity with $\eps+0=\eps$ we have $a'd \dinde\eps_C e$,
  and so in particular $a' \dinde\eps_C e$.
    \item
  If $ab \nind^0_C d$ then
  $a \nind^0_{Cb} d$ or $b \nind^0_C d$, so
  $a \ndinde\eps_{Cb} d$ or $b \ndinde\eps_C d$ by the wgp assumptions,
  so $ab \ndinde\eps_C d$ as required.
    \item
      If $a \nind^0_{Cb} d$,
      then $a \nind^0_C bd$,
      so $a \ndinde{\eps'}_C bd$ by $\eps'$-wgp,
      so $a \ndinde{\eps'-\eps}_{Cb} d$ as required.
  \end{enumerate}
\end{proof}

\begin{definition}
  A sequence $\bar{a}=(a_i)_{i\in\omega}$ is called a $\dinde\eps$-independent 
  sequence, if $\bdl(a_i)<\infty$ and $a_{i+1} \dinde\eps a_0\ldots a_i$ for 
  all $i$.
\end{definition}


\begin{lemma}\label{l:epscgpSeq}

  Suppose $\a = (a_0,\ldots,a_k)$ is a finite $\dinde\eps$-independent sequence.
  Let $b \in \acl^0(\a)$.
\begin{enumerate}
    \item Assume $\tp(a_i)$ is $\gamma$-cgp for all $i$, and
  $\bdl(a_i) = \dim^0(a_i)$.
  Then $\bdl(b) \geq \dim^0(b) - k\eps - (k+1)\gamma$.
  \item Assume $\max\{\bdl(a_i):i\leq k\}=1=\dim^0(a_i)$ for all $i$ and $\min\{\bdl(a_i):i\leq k\}\geq r$ for some $r > \eps$. Then $\bdl(b)\geq r\dim^0(b)- k\eps$.
\end{enumerate}

\end{lemma}
\begin{proof}
  We first prove item (1).

  We have $\bdl(\a) \geq \dim^0(\a) - k\eps$.

  We now show that $\bdl(\a/b) \leq \dim^0(\a/b) + (k+1)\gamma$ (for this we do 
  not use $b \in \acl^0(\a)$).
  For $k=0$: if $\a \nind^0 b$ then $\bdl(\a/b) \leq \gamma$ by $\gamma$-cgp, 
  and if $\a \ind^0 b$ then $\bdl(\a/b) \leq \bdl(\a) = \dim^0(\a) = 
  \dim^0(\a/b)$.
  For $k\geq 1$, we argue by induction. Say $\a = \a'a$. Then 
  \begin{align*}
  \bdl(\a/b) &= \bdl(\a'/ba) + \bdl(a/b)\\
   &\leq \dim^0(\a'/ba) + k\gamma + \dim^0(a/b) + 
  \gamma = \dim^0(\a/b) + (k+1)\gamma.
  \end{align*}

  So $$\bdl(b)
  = \bdl(\a) - \bdl(\a/b)
  \geq \dim^0(\a) - \dim^0(\a/b) - k\eps - (k+1)\gamma =
  \dim^0(b) - k\eps - (k+1)\gamma.$$

  Now we prove item (2).

  Note that $\bar{a}$ is an $\acl^0$-independent sequence, since if $a_i\in\acl^0(a_j:j<i)$, then $\bdl(a_i/(a_j)_{j<i})=0<r-\eps$ contradicts $a_i\dinde\eps (a_j)_{j<i}$.

  Let $\bar{a}'=(a_i)_{i\in I'}$ be a minimal subtuple of $\bar{a}$ such that $\bar{a}'$ is a transcendence base of $\bar{a}$ over $\acl^0(b)$. Then $\dim^0(\bar{a}')+\dim^0(b)=\dim^0(\bar{a}'/b)+\dim^0(b)=\dim^0(\bar{a})$, and $\dim^0(b)=\dim^0(\bar{a})-\dim^0(\bar{a}')$. As $\bar{a}$ is $\acl^0$-independent, we have $|\{0,\ldots,k\}\setminus I'|=\dim^0(b)$. Hence,
  \begin{align*}
  \bdl(b)&=\bdl(\bar{a})-\bdl(\bar{a}'/b)\geq \bdl(\bar{a})-\bdl(\bar{a}')\\
  &\geq \sum_{ i\leq k}\bdl(a_i)-k\eps-\sum_{i\in I'}\bdl(a_i)\geq r\dim^0(b)-k\eps.\qedhere
  \end{align*}
\end{proof}

\begin{definition}\label{d:epsTri}
  Let $\eps' \geq \eps \geq 0$.
  Let $C \subseteq K$.
  A \defn{coarse $(\eps,\eps')$-correspondence triangle} (or just 
  \defn{$(\eps,\eps')$-correspondence triangle}) over $C$ is a triple $(a,d,b)$ such 
  that
  \begin{itemize}
    \item $\tp(a/C)$ is wgp;
    \item $\tp(d/C)$ is $\eps'$-wgp;
    \item $a,d \notin \acl^0(C)$;
    \item $a \sim^0_{Cd} b$;
  \item $d \sim^0_C \cb^0(ab/Cd)$;
  \item $a \dinde\eps_C d \dinde\eps_C b$.
  \end{itemize}

  If $\tp(a/C)$ is cgp, we call $(a,d,b)$ a \defn{semi-cgp 
  $(\eps,\eps')$-correspondence triangle} over $C$.

  In the case $C=\emptyset$, we omit mention of it.

  We define abelianity of $(\eps,\eps')$-correspondence triangles exactly as 
  for correspondence triangles (which are precisely $(0,0)$-correspondence 
  triangles).
\end{definition}

\begin{remark}\label{r:epsTriInd}
  If $(a,d,b)$ is an $(\eps,\eps')$-correspondence triangle over $C$ then 
  $a\ind^0_Cd\ind^0_Cb$, since $\eps'\geq\eps$.
\end{remark}

\begin{remark}\label{r:epswgp}
  If $(a,d,b)$ is an $(\eps,\eps')$-correspondence triangle over $C$, then 
  $\tp(b/C)$ is also wgp.

  Indeed, since $d\dinde\eps_C a$ and $\tp(d/C)$ is $\eps'$-wgp, by 
  Lemma~\ref{l:epswgpFacts}(iii), $\tp(d/aC)$ is $(\eps'-\eps)$-wgp, and in 
  particular is wgp. Since $\tp(a/C)$ is wgp, also $\tp(da/C)$ is wgp, hence 
  $\tp(b/C)$ is also wgp by Lemma~\ref{l:epswgpFacts}(i)(ii). Note also that 
  $\tp(b/C)$ is broad. Recall that for a wgp tuple, broad is equivalent to 
  non-algebraic. Hence, if $b\in\acl^0(C)$ then $a\in\acl^0(Cd)$ contradicting 
  $a\notin\acl^0(C)$ and $a\ind^0_C d$.

\end{remark}

Using an explicit Szemerédi-Trotter bound, we get an explicit bound nf 
$\bdl(d)$ in terms of $\bdl(a)$ and the algebraic dimension of $d$ when 
$(a,d,b)$ is a semi-cgp $(\eps,\eps')$-correspondence triangle.
\begin{lemma} \label{l:epsSTCb}
  Let $(a,d,b)$ be a semi-cgp $(\eps,\eps')$-correspondence triangle.
  Suppose also that $\tp(a/C)$ is $3\eps$-broad.
  Then $\bdl(d) \leq \bdl(a)+(6\dim^0(d) + 5)\eps$.
\end{lemma}
\begin{proof}
  We may assume $d = \cb^0(ab/d)$.
  Exactly as in the proof of \cite[Proposition~5.14]{BB-cohMod}, applying 
  Szemerédi-Trotter bounds to the $\bigwedge$-definable binary relation 
  $(\tp(ab)\times\tp(d)) \cap \loc^0(ab,d)$ and using that $\tp(a)$ is cgp,
  we obtain
  \begin{equation} \label{e:epsST}
    \bdl(abd) \leq  \max(\frac12 \bdl(ab) + \bdl(d) - 
    \max(0,\eps_0(\bdl(d) - \frac12 \bdl(ab))), \bdl(ab), 
    \bdl(d))
  \end{equation}
  where $\eps_0 > 0$.
  As in the proof of \cite[Lemma~3.11]{BB-cohMod}, by the primitive element 
  theorem we may assume that $d$ is a $(\dim^0(d)+1)$-tuple, and then the 
  semialgebraic incidence bounds proven in \cite{FoxPachEtAl}[Theorem~1.2] 
  yield that \eqnref{epsST} holds for any $\eps_0 < \frac1{4(\dim^0(d)+1)-1}$, 
  and hence also for $\eps_0 = \frac1{4(\dim^0(d)+1)-1}$.

  Suppose for a contradiction that $\bdl(d) > \bdl(a)+(6\dim^0(d)+5)\eps$.
  Rearranging, we obtain $\bdl(d) > \bdl(a) + \frac\eps2(\frac3{\eps_0}+1)$ 
  where
  $\eps_0 = \frac3{2(6\dim^0(d)+5) - 1} = \frac1{4(\dim^0(d)+1)-1}$, so 
  \eqnref{epsST} holds for $\eps_0$.

  We have $\bdl(b/a) \leq \bdl(b) \leq \bdl(b/d) + \eps = \bdl(a/d) + \eps 
  \leq \bdl(a) + \eps$,
  so \[\frac12 \bdl(ab) = \frac12 (\bdl(b/a)+\bdl(a)) \leq \bdl(a) + \frac\eps2 
  < \bdl(d) - \frac\eps2\frac3{\eps_0}\]
  and so $\eps_0(\bdl(d) - \frac12 \bdl(ab))>\frac{3\eps}2$.
  But \[\bdl(abd) = \bdl(a/d) + \bdl(d) \geq \bdl(a) - \eps + \bdl(d) \geq  
  \frac12 \bdl(ab) + \bdl(d) - \frac{3\eps}2,\]
  and \[\frac12 \bdl(ab) + \bdl(d) - \frac{3\eps}2 > \frac12\bdl(ab) + 
  \frac12\bdl(ab) = \bdl(ab)\] (using $\eps_0\leq 1$),
  and $\frac12\bdl(ab) + \bdl(d) - \frac{3\eps}2> \bdl(d)$ since 
  $\frac12\bdl(ab) \geq \frac12\bdl(a) > \frac{3\eps}2$ since $\bdl(a)>3\eps$,
  so this contradicts \eqnref{epsST}.
\end{proof}

The next two lemmas are quantitative versions of Lemma~\ref{l:vee} and 
Lemma~\ref{l:semicgpind} respectively.
\begin{lemma}\label{l:epsvee}
  Let $\eps' \geq 3\eps \geq 0$.
  Suppose $(a,d,b)$ and $(a',d',b)$ are $(\eps,\eps')$-correspondence 
  triangles such that
  $a'd' \dind_b ad$.
  Let $e:=\cb^0(aa'/dd')$.
  Then $(a,e,a')$ is a $(2\eps,\eps'-\eps)$-correspondence triangle,
  and $d\sim^0_ed'$ and $\dim^0(e)\geq\dim^0(d)$.
\end{lemma}
\begin{proof}
  We first verify the conditions for $(a,e,a')$ to be a $(2\eps,\eps'-\eps)$-correspondence triangle.
  \begin{itemize}
    \item $\tp(a)$ is wgp and broad by assumption.
      We have $d\dinde\eps d'$,
      and $\tp(d)$ is $\eps'$-wgp,
      so $\tp(d/d')$ is $(\eps'-\eps)$-wgp,
      hence so is $\tp(dd')$ and hence $\tp(e)$.
      We verify below that $\tp(e)$ is $(\eps'-\eps)$-broad, and hence $e 
      \notin\acl^0(\emptyset)$.
    \item $a \sim^0_d b \sim^0_{d'} a'$,
      hence $a \sim^0_{dd'} a'$, and so also $a \sim^0_e a'$ by choice of 
      $e$.
    \item $e = \cb^0(aa'/e)$ by definition.
    \item
      We have
      $ad \dinde\eps d'$
      so $a\dinde\eps_d d'$
      and $a\dinde\eps d$,
      so $a\dinde{2\eps} dd'$,
      so $a\dinde{2\eps} e$
      since $e \in \acl^0(dd')$,
      and similarly $a'\dinde{2\eps} e$.
    \end{itemize}
    Now we claim that $(d,a,b,d',a')$ satisfies (a) and (b) in 
    Proposition~\ref{p:recogCFAHS}. Condition (a) follows directly 
    from the assumptions. For (b), we only need to prove 
    $a'd'\ind^0_bad$. By assumption, $a'd'\dind_bad$, hence $d'\dind_bd$ and 
    $d'\dinde\eps bd$ (as $d'\dinde\eps b$). Since $\tp(d')$ is $\eps$-wgp, we 
    have $d'\ind^0bd$ and $d'\ind^0_bd$. Since $a\in\acl^0(bd)$ and 
    $a'\in\acl^0(bd')$, we conclude $a'd'\ind^0_bad$ as desired.

  So by Proposition~\ref{p:recogCFAHS}, $d'\sim^0_e d$.
  Using $e\in \acl^0(dd')$, we deduce $e\sim^0_{d'} d$.
  Since $d\dinde\eps d'$, we obtain
  $\bdl(e)\geq \bdl(e/d') = \bdl(d/d') \geq \bdl(d) - \eps > \eps'-\eps$,
  and hence $\tp(e)$ is $(\eps'-\eps)$-broad.
  Finally, the analogous calculation using $d \ind^0 d'$ yields
  $\dim^0(e) \geq \dim^0(d)$.
\end{proof}

\begin{lemma}\label{l:epssemicgpind}
  If $(a,d,b)$ is a semi-cgp $(\eps,\eps')$-correspondence triangle, then 
  $a\ind^0 b$.
\end{lemma}
\begin{proof}
  Suppose not. Then $\bdl(a/b)=0$ since $\tp(a)$ is cgp, and so $d\dind_ba$.
  By $d\dinde\eps b$ and transitivity, we obtain $d\dinde\eps ab$, and hence 
  $d\ind^0ab$ since $\tp(d)$ is $\eps$-wgp.
  But $d = \cb^0(ab/d)$, contradicting $d\notin\acl^0(\emptyset)$.
\end{proof}

\subsection{Quantitative non-expansion for actions on curves}
    We now give a quantitative version of Theorem~\ref{t:cohActAb} in the case 
    that $X$ is an algebraic curve.

    \begin{fact}\cite[Theorem 1.3]{JZ}\label{fact-JZ}
For all $n\in\mathbb{N}^{>0}$, $\gamma, \gamma'>0$ and $r,k\geq 1$, there are $M\in\mathbb{N}$ and $\delta>0$ such that the following holds.
    Suppose $G$ is a subgroup of $\mathrm{PGL}_{n}(K)$ or $\mathrm{GL}_{n}(K)$ for some field $K$ of characteristic 0 and that $G$ acts on some set $X$.
    Suppose there are finite sets $S\subseteq G$, $T\subseteq X$ with $G=\langle S\rangle$ satisfying:
\begin{enumerate}
    \item
    $M < \log|T|<r\log|S|$;
    \item
    $|S\cap gH|<|S|^{1-\gamma}$ for any nilpotent subgroup $H\leq G$ of step at most $n-1$ and any $g\in G$;
    \item
    $|\{(x_1,\ldots,x_k)\in T^k:\Stab_G(x_1,x_2,..,x_k)\neq \mathrm{id}_G\}| \leq \frac1M |T|^{k-\gamma'}$.
\end{enumerate}
Then $|\{(x,y,g):x,y\in T,\;g\in S,\;x=g*y\}|\leq |T|^{1-\delta}|S|$. In fact any $\delta<\min\{\frac{\gamma'}{k}, 2^{-\lceil\frac{2kr}{c\gamma}\rceil-2}\}$ works, where $c$ is a constant depending only on $n$.
\end{fact}

The following fact is a special case of the general model theoretic result \cite[Theorem~3.27]{poizatStableGroups}. It can alternatively be seen by reducing to the case of a smooth projective curve and considering its automorphism group, taking cases according to the genus.
\begin{fact}\label{f:CFAHSCurve}
  If $X$ is an algebraic curve, then any CFAHS $(G,X)$ over $\acl^0$ is isomorphic over $\acl^0$ to one of the following:
  \begin{itemize}
    \item $(G,G)$ for a connected commutative 1-dimensional algebraic group $G$ acting on itself by addition.
    \item $(\AGL_1,\A^1)$, the group of affine linear transformations $\AGL_1 = \G_a \rtimes \G_m$ acting on the affine line.
    \item $(\PSL_2,\mathbb{P}^1)$, the group of Möbius transformations acting on the projective line.
  \end{itemize}
\end{fact}

\begin{lemma}\label{l:eps1dimab}
There is an absolute constant $c'>0$ for which the following holds.\footnote{This constant comes from the constants for the Balog-Szemerédi-Gowers theorem and the product theorem for $\mathrm{SL}_2(\mathbb{C})$.}

  Let $(G,X)$ be a CFAHS defined over $\acl^0(\emptyset)$ with $\dim^0(X) = 1$.

  Let $a \in X(K)$ and $g \in G(K)$.
  Suppose $\tp(g)$ is $\eps'$-wgp, $g \notin \acl^0(\emptyset)$, $\tp(a)$ is broad,
  and $a \dinde\eps g \dinde\eps g*a$,
  where $\eps< \min\{\bdl(a)2^{-\frac{\bdl(a)}{c'\eps'}-5},\frac{\bdl(g)}{4},\frac{\eps'}{4}\}$.

  Then $g$ is generic in a left coset over $\acl^0(\emptyset)$ of a connected one-dimensional subgroup $H\leq G$.

  In particular, if $g\in G(K)$ is generic, then $G$ is abelian of dimension 1, and so $(G,X)$ is isomorphic over $\acl^0(\emptyset)$ 
  with the homogeneous space $(G,G)$ in which $G$ acts on itself by addition.
\end{lemma}
\begin{proof}
  We conclude immediately if $\dim^0(G) = 1$. Otherwise, by Fact~\ref{f:CFAHSCurve},
we may assume $(G,X)=(\AGL_1,\A^1)$ or $(G,X)=(\PSL_2,\mathbb{P}^1)$. We conclude by showing that $Y := \loc^0(g)$ is a left coset of a connected one-dimensional subgroup of $G$. Suppose not, then either $\dim^0(g)\geq 2$ or $\dim^0(g)=1$ and $Y\cap hH$ is finite for any one-dimensional subgroup $H\leq G$ and $h\in G$.

Note that $\bdl(g*a)\leq \bdl(g*a/g) + \eps = \bdl(a/g) + \eps \leq \bdl(a)+\eps$. Similarly, $\bdl(a)\leq\bdl(g*a)+\eps$.

   Let $\Gamma_{(G,X)}=\{(x,y,g)\in X^2\times G: x=g*y\}$. Let $\pi_i:\Gamma_{(G,X)}\to X$ for $i=0,1$ and $\pi_2:\Gamma_{(G,X)}\to G$ be the co-ordinate projections.
   By Fact~\ref{f:existsInternal}, there is an internal set $A\subseteq\tp((g*a)ag)$ such that $\bdl(A)=\bdl((g*a)ag)\geq\bdl(a)+\bdl(g)-\eps$.
   Let $T:=\pi_0(A)\cup\pi_1(A)$. Then $\bdl(T)=\max\{\pi_0(A),\pi_1(A)\}\leq\max\{\bdl(a),\bdl(g*a)\}\leq\bdl(a)+\eps$.
   Let $S:=\pi_2(A)$. Then $\bdl(S)\leq\bdl(g)$.

   Note that $\bdl(S)+\bdl(T)\geq\bdl((T^2\times S)\cap\Gamma_{(G,X)})\geq\bdl(A)\geq\bdl(g)+\bdl(a)-\eps$. Combining with $\bdl(T)\leq\bdl(a)+\eps$ and $\bdl(S)\leq\bdl(g)$, we get $\bdl(S)\geq\bdl(g)-2\eps$ and $\bdl(T)\geq\bdl(a)-\eps$.
   Then $\bdl(T)\leq\bdl(a)+\eps\leq\frac{\bdl(S)}{\bdl(g)-2\eps}(\bdl(a)+\eps)<\frac{3\bdl(a)}{\bdl(g)}\bdl(S)$, where the last inequality is by $\eps<\frac{\bdl(g)}{4}$ and $\eps<\frac{\bdl(a)}2$.
   Also, 
   \begin{align*}
  \bdl((T^2\times S)\cap\Gamma_{(G,X)})
  &\geq \bdl(a)(1-\frac{\eps}{\bdl(a)})+\bdl(g)
  \geq (\bdl(T) - \eps)(1-\frac\eps{\bdl(a)})+\bdl(S)\\
   &\geq \bdl(T)(1-\frac\eps{\bdl(a)}) - \eps + \bdl(S)
  > \bdl(T)(1-\frac{3\eps}{\bdl(a)})+\bdl(S),
   \end{align*}
   since $1 < \frac{2\bdl(T)}{\bdl(a)}$.

   Suppose $T=\lim_{n\to\mathcal{U}}T_n$ and $S=\lim_{n\to\mathcal{U}}S_n$.
   Let $M \in \N$.
   Then for ultrafilter-many $n$, we have $M < \log|T_n|\leq 3\frac{\bdl(a)}{\bdl(g)}\log|S_n|$ and \[|\{(x,y,g):x,y\in T_n,\;g\in S_n,\;x=g*y\}|\geq |T_n|^{1-\frac{3\eps}{\bdl(a)}}|S_n|.\]

   Let $G_n=\langle S_n\rangle$, then $G_n\leq G(\C)$. Since $G_n$ is a group of affine linear transformations or of elements of $\mathrm{PSL}_2(\C)$ acting on $\mathbb{P}^1(\C)$, at most three distinct points determine an element in $G_n$. Hence, 
   \begin{align*}
    &|\{(x_1,x_2,x_3)\in T_n:\Stab_{G_n}(x_1,x_2,x_3)\neq \mathrm{id}\}|\\
    \leq &|\{(x_1,x_2,x_3)\in T_n^3:x_i=x_j \text{ for some }i\neq j \}|\leq 3|T_n|^2.
   \end{align*}

   Moreover, we claim there are ultrafilter-many $n$, such that $|S_n\cap hH|<|S_n|^{1-\frac{\eps'}{2\bdl(g)}}$ for any abelian subgroup $H\leq G_n$ and $h\in G_n$. Suppose not, then for ultrafilter-many $n$, there are abelian subgroups $H_n\leq G_n$ and $h_n\in G_n$ such that $|S_n\cap h_nH_n|\geq |S_n|^{1-\frac{\eps'}{2\bdl(g)}}$. Let $H:=\prod_{n\to\mathcal{U}}H_n$ and $h:=(h_n)_{n}/\mathcal{U}$. Then $\bdl(\tp(g)\cap hH)\geq \bdl(S\cap hH)\geq(1-\frac{\eps'}{2\bdl(g)})\bdl(S)\geq\bdl(S)-\frac{\eps'}{2}\geq\bdl(g)-2\eps-\frac{\eps'}{2}\geq \bdl(g)-\eps'$, as $\eps<\frac{\eps'}{4}$ and $\bdl(g) \geq \bdl(S)$.
   Let $\bar H$ be the Zariski closure of $H$, then $\bar H$ is also an abelian subgroup of $G$. Since $G=\AGL_1$ or $G=\mathrm{PSL}_2$,  we have $1\geq \dim^0(\bar H)=\dim^0(h\bar H)$. If $\dim^0(g)\geq 2$, by $\eps'$-wgp of $\tp(g)$ we have $\bdl(\tp(g)\cap h\bar H)<\bdl(g)-\eps'$, and if $\dim^0(g)=1$, since $Y\cap h\bar H$ is finite we again have $\bdl(\tp(g)\cap h\bar H)=0<\bdl(g)-\eps'$ (since $\tp(g)$ is $\eps'$-wgp and non-algebraic, so $\bdl(g)>\eps'$). This contradicts $\bdl(\tp(g)\cap hH)\geq \bdl(g)-\eps'$.

   Now we may apply Fact \ref{fact-JZ} with $k:=3, r:=3\frac{\bdl(a)}{\bdl(g)}, \gamma:=\frac{\eps'}{2\bdl(g)}$ and $\gamma':=3/4$. We get \[\min\left\{2^{-\lceil \frac{36\bdl(a)}{c\eps'}\rceil-2},\frac{1}{4}\right\}\leq \frac{3\eps}{\bdl(a)},\]
   where $c$ is the absolute constant from the bound of Fact \ref{fact-JZ} for 
   $n=2$. Setting $c':=c/36$, then \[\eps\geq \frac{\bdl(a)}{3}2^{-\lceil 
   \frac{\bdl(a)}{c'\eps'}\rceil-2}\geq \bdl(a) 2^{-\frac{\bdl(a)}{c'\eps'}-5},\]
   contradicting the assumption on $\epsilon$. \end{proof}

\subsection{Abelianity of cgp $(\eps,\eps')$-correspondence triangles}
Recall that Proposition~\ref{p:maincgp} asserts that cgp correspondence triangles are abelian. We now present a quantitative version of this result. We also include a special case of semi-cgp correspondence triangles, which will later be used to deduce the quantitative asymmetric Elekes–Rónyai theorem.
\begin{proposition} \label{p:epsmaincgp}
  Let $(a,d,b)$ be a semi-cgp $(\eps,\eps')$-correspondence triangle.

Suppose one of the following holds
  \begin{enumerate}
      \item $d$ is $\gamma$-cgp and
        $\bdl(d) = \dim^0(d)$;
      \item 
        $d\in \acl(c_i:i< m)$ with each $\dim^0(c_i)=1$ and $(c_i)_{i<m}$ is a $\dind$-independent sequence such that $\max\{\bdl(c_i):i< m\}=1$ and $\bdl(c_i)\geq r$ for each $i< m$.
  \end{enumerate}

  In case (1), set $r:=1$ and $m:=\dim^0(d)$.
  In case (2), set $\gamma:=0$.
  Let $i_0 := 1+\lfloor\frac{\bdl(a)}{r}\rfloor$.
  Suppose $\eps \leq \min\{\frac{r-2^{i_0}\gamma}{2^{2i_0+4}m},\frac{\eps'}{2^{i_0+1}},\frac{\eps'}{20}\}$.

  Suppose further that $\dim^0(a)=1$ and $\eps<\bdl(a)2^{-\frac{2\bdl(a)}{c'\eps'}-7}$ where $c'$ is the absolute constant from Fact \ref{fact-JZ}.

  Then $(a,d,b)$ is abelian.
\end{proposition}
Before turning to the details, let us summarise the proof strategy. We follow the proof of Proposition~\ref{p:maincgp}. Starting from $(a,d,b)$, we recursively construct $(\epsilon_i,\epsilon_i’)$-correspondence triangles $(a,d_i,b_i)$, where the parameters $\epsilon_i$ and $\epsilon_i’$ are worsened exponentially with $i$. As in the proof of Proposition~\ref{p:maincgp}, the goal is to show that $\dim^0(d_i)$ stabilises after finitely many steps.

To achieve this, we use the explicit Szemerédi–Trotter bound established in Lemma~\ref{l:epsSTCb}, which yields an explicit bound $i_0$ on the number of iterations required before $\dim^0(d_i)$ stabilises. Once this happens, we conclude that $(a,d_i,b_i)$ is generic in the graph of a CFAHS $(G,X)$. Under the additional assumption that $\dim^0(a)=1$, which implies that $X$ is an algebraic curve, we use Lemma~\ref{l:eps1dimab} to conclude that such non-expansion cannot occur unless $G$ is abelian.

Note that the exponential dependence on $i_0$ arises because the error parameters are weakened by a squaring operation at each iteration.
\begin{proof}
  \newcommand{\bh}{\bar h}

 We recursively define a sequence $(d_i,d_i',b_i)_{i \leq i_0}$ such that, 
 for all $0<i\leq i_0$,
  \begin{enumerate}[(i)]
    \item $(a,d_i,b_i)$ is a $(2^i\eps,\eps' - (2^i-1)\eps)$-correspondence 
      triangle,
    \item $\dim^0(d_i)\geq\dim^0(d_{i-1})$,
    \item $b_id_{i-1}' \dind_{b_{i-1}} ad_{i-1}$,
      and $\dim^0(d_{i-1}') = \dim^0(d_{i-1})$,
      and $(b_i,d_{i-1}',b_{i-1})$ is a 
      $(2^{i-1}\eps,\eps'-(2^{i-1}-1)\eps)$-correspondence triangle,
    \item $d_{i-1} \sim^0_{d_i} d_{i-1}'$,
    \item $d_i \in \acl^0(\bh)$ for some finite $\dinde{2^i\eps}$-independent 
      sequence $\bh$ of $2^{i}$ realisations of $\tp(d)$ in case (1) and  of $\tp((c_i)_{i<m})$ in case (2), with $\bh 
      \dinde{2^i\eps} b_i$ and $\bh \dinde{2^i\eps} a$.
  \end{enumerate}

  First, let $(a,d_0,b_0):=(a,d,b)$, $\bar h:=d$ in case (1) and $\bar{h}:=(c_i)_{i<m}$ in case (2). Suppose $0 \leq i < i_0$ and we have defined 
  $(d_i,b_i)$ and $(d_{j},d_{j}',b_{j})_{j<i}$ satisfying (i)-(v).
  We define $d_i'$, $d_{i+1}$, $b_{i+1}$ such that 
  $(d_i,d_i',b_i,d_{i+1},b_{i+1})$ satisfy (i)-(v).

  Let $\bh$ be as in (v) for $d_i$.
  Let $b_{i+1}d_i'\bh' \equiv_{b_i} ad_i\bh$ with $b_{i+1}d_i' \bh'\dind_{b_i} ad_i\bh$,
  and set $d_{i+1} := \cb^0(ab_{i+1}/d_id_i')$.
  Note that $\eps'-(2^i-1)\eps\geq3\cdot2^i\eps$ since $\eps'-(2^i-1)\eps - 3\cdot2^i\eps = \eps'-(2^{i+2}-1)\eps \geq 0$.

  So (i),(ii),(iv) hold for $i+1$ by Lemma~\ref{l:epsvee}, and (iii) holds by 
  construction. By definition, $d_{i+1}\in\acl^0(d_id_i') \subseteq 
  \acl^0(\bh\bh')$. Note that $\bh\dinde{2^i\eps} \bh'$, since 
  $\bh\dinde{2^i\eps} b_i$ by assumption. Therefore, $\bh\bh'$ is a 
  $\dinde{2^{i+1}\eps}$-independent sequence in $\tp(d)$ or in $\tp((c_i))_{i<m})$. Moreover, we have 
  $b_{i+1}\bh'\dinde{2^i\eps} \bh$ (again by $\bh\dinde{2^i\eps} b_i$).
  Now $\bh' \dinde{2^i\eps} b_{i+1}$ since $\bh \dinde{2^i\eps} a$,
  so $\bh\bh'\dinde{2^{i+1}\eps} b_{i+1}$.
  Also, $\bh' \dinde{2^i\eps} b_i$ since $\bh \dinde{2^i\eps} b_i$, so $\bh' \dinde{2^i\eps} \bh a$, and since also $\bh \dinde{2^i\eps} a$, we have $\bh\bh' \dinde{2^{i+1}\eps} a$.
  So (v) holds, witnessed by $\bh\bh'$.

  This concludes our construction of the sequence 
  $(d_i,d_i',b_i)_{i\in\omega}$.

  Now for all $i<i_0$, we have:
  \begin{itemize}
    \item $a\ind^0d_i\ind^0b_i\ind^0a$ and 
      $b_{i+1}\ind^0d'_i\ind^0b_i\ind^0b_{i+1}$ by (i) (which also holds for 
      $i=0$), (iii), Remark~\ref{r:epsTriInd}, and Lemma~\ref{l:epssemicgpind};
    \item $d_i'\ind^0 b_id_i$, and hence $b_{i+1}d_i'\ind^0_{b_i}ad_i$:
      indeed, $d_i' \dinde{2^i\eps} b_id_i$ and $d_i'$ is 
      $(\eps'-(2^i-1)\eps)$-wgp and $\eps'-(2^i-1)\eps > 2^i\eps$, so $d_i' 
      \ind^0 b_id_i$.
  \end{itemize}

  \begin{claim} \label{c:epscbDimBd}
    For all $0 < i \leq i_0$,
    $$\dim^0(d_{i-1}) \leq \dim^0(d_i)
    \leq \frac{1}{r}(\bdl(a) + 2^{2i+4}m\eps + 2^i\gamma).$$
  \end{claim}
  \begin{proof}
    The first inequality is by (ii).

    Let $r=1$ in case (1) and $\gamma=0$ in case (2). Note that in both cases, we have $\bdl(d)\geq r\dim^0(d)$ (by Lemma~\ref{l:epscgpSeq} in case (2)) and $m\geq\dim^0(d)$.
    By Lemma~\ref{l:epsSTCb}, we have \[\bdl(a) \geq \bdl(d)-(6\dim^0(d)+5)\eps \geq r\dim^0(d)-(6m+5)\eps \geq \frac{r}{2}
    >3\cdot2^{i_0}\eps \geq 3\cdot2^i\eps,\] since $(6m+5)\eps\leq \frac{r}{2}$ and $\eps<\frac{r}{6\cdot 2^{i_0}}$,
    so by (i) and Lemma~\ref{l:epsSTCb} again, we obtain 
    \begin{equation}\label{eq:bdldi}
    \bdl(d_i) \leq \bdl(a) + 
    2^i\eps(6\dim^0(d_i)+5).
    \end{equation}

    In case (1), by (v), $d_i \in \acl^0(\bh)$ for some $\dinde{2^i\eps}$-independent sequence 
    $\bh$
    of $2^{i}$ realisations of $\tp(d)$.
    By Lemma~\ref{l:epscgpSeq} (and 
    since we assumed $\dim^0(d) = \bdl(d)$), we obtain
    $\dim^0(d_i) \leq \bdl(d_i) + 2^{2i}\eps + 2^i\gamma$.

    In case (2), since $d_i \in \acl^0(\bh)$ for some $\dinde{2^i\eps}$-independent sequence 
    $\bh$
    of $2^{i}$ realisations of $\tp((c_i)_{i<m})$ and $(c_i)_{i< m}$ is a $\dind$-independence sequence, we have $\bh$ 
    is a
    $\dinde{2^i\eps}$-independent sequence of $m\cdot 2^i$ realisations of the types among $\{\tp(c_i):{i<m}\}$.
    So by Lemma~\ref{l:epscgpSeq}, we obtain 
   $r\dim^0(d_i)\leq \bdl(d_i)+2^{2i}m\eps.$

    In both cases, \[\dim^0(d_i) \leq \frac{1}{r}(\bdl(d_i) + 2^{2i}m\eps + 2^i\gamma).\]
    Combining with the inequality~(\ref{eq:bdldi}), we conclude \[\dim^0(d_i)
    \leq \frac{1}{r}(\bdl(a) + (2^i(6\dim^0(d_i)+5)+2^{2i}m)\eps+2^i\gamma).\]
    Using the coarse bounds $\dim^0(d_i) \leq 2^i\dim^0(d)\leq 2^im$, this gives
    $\dim^0(d_i) \leq \frac{1}{r}(\bdl(a) + 2^{2i+4}m\eps+2^i\gamma)$.
    \subqed{c:epscbDimBd}
  \end{proof}

  Suppose for a contradiction that
  $\dim^0(d_{i+1}) > \dim^0(d_i)$ for all $i < i_0$.
  Then \[\dim^0(d) + i_0 \leq \dim^0(d_{i_0}) \leq \frac{1}{r}(\bdl(a) + 2^{2i_0+4}m\eps +   2^{i_0}\gamma).\]
By our assumptions, $\eps \leq 
  \frac{r-2^{i_0}\gamma}{2^{2i_0+4}m}$, so 
  \[\dim^0(d) + i_0 \leq \frac{1}{r}(\bdl(a) + 2^{2i_0+4}m\eps + 
  2^{i_0}\gamma)\leq  \frac{1}{r}(\bdl(a)+r- 2^{i_0}\gamma+ 
  2^{i_0}\gamma)=\frac{\bdl(a)}{r}+1.\]
  Thus $i_0 \leq \frac{\bdl(a)}{r}$ (as $\dim^0(d)\geq 1$ by the definition of $(\eps,\eps')$-correspondence triangle),
  contradicting
  $i_0 \geq \lfloor\frac{\bdl(a)}{r}\rfloor+1$.

  So $\dim^0(d_{i+1}) = \dim^0(d_i)$ for some $i < i_0$.

  As in the proof of Proposition~\ref{p:maincgp},
  Proposition~\ref{p:recogCFAHS} yields a CFAHS $(G,X)$ over 
  $\acl^0(\emptyset)$ such that up to $\acl^0$-interalgebraicity, 
  $(d_i,a,b_i)$ is generic in the graph of the action.

  Now, assume $\dim^0(a)=1$, then $\dim(X)=1$ and $\dim^0(d_i)\leq 3$, since $\dim(G)$ is bounded by 3 by Fact~\ref{f:CFAHSCurve}. Hence, already $\dim^0(d_{i+1})=\dim^0(d_i)$ for $i\leq 2$. Now, we have a $(4\eps,\eps'-3\eps)$-correspondence triangle $(a,d_2,b_2)$ which is, up to $\acl^0$-interalgebraicity, 
 generic in the graph of the action of $G$ on $X$.

  By Lemma~\ref{l:eps1dimab}, $G$ is commutative provided \[4\eps<\min\{\bdl(a) 2^{-\frac{\bdl(a)}{c'(\eps'-3\eps)}-5},\frac{\eps'-3\eps}{4},\frac{\bdl(d_2)}{4}\},\] which holds when $\eps<\bdl(a) 2^{-\frac{2\bdl(a)}{c'\eps'}-7}$ (since $\eps'-3\eps>\eps'/2$) and $\eps<\frac{\eps'}{20}$ (note that $4\eps<r/4\leq \bdl(d)/4\leq\bdl(d_2)/4$). And we then conclude exactly 
  as in the proof of Proposition~\ref{p:maincgp}.
\end{proof}

\begin{remark}
  Proposition~\ref{p:epsmaincgp} is our quantitative analog of Proposition~\ref{p:maincgp} for cgp triangles.
  We could hope the same method would yield a quantitative version of 
  Theorem~\ref{t:main} for the semi-cgp triangles, by perturbing the assumptions and trying to control the error occurring in the proof. Attempting to do this leads however to the problem that the 
  parameter $c$ in the proof of Proposition~\ref{p:maincgp}, chosen such that 
  $\tp(d/c)$ is cgp, has unbounded transcendence degree $\dim^0(c)$. If we 
  assume that the semi-cgp triangle $(a,d,b)$ is an 
  $(\eps,\eps')$-correspondence triangle, we obtain $\bdl(c) \geq \eps'$, but 
  in order to later obtain $c \ind^0 c'$ we need to be able to add parameters 
  $f$ such that $c$ becomes $\eps$-wgp. However, ensuring this appears to 
  require assuming $\eps' \geq \dim^0(c)\eps$. Since the geometry of the 
  situation gives no upper bound on $\dim^0(c)$, this prevents us from giving 
  a direct bound on the powersaving as in Proposition~\ref{p:epsmaincgp}. For 
  this reason, we content ourselves with qualitative statements in this 
  regime.
\end{remark}

\subsection{The case of a polynomial}
Here, we deduce our quantative analogue of Lemma~\ref{l:addmul}, used in the explicit Elekes-Rónyai result Theorem~\ref{t:unbalancedERMany}.
\begin{lemma}\label{l:epsaddmul}
  Let $f\in K[t]$, and suppose that for some $a \in K$, setting $b := f(a)$ and identifying $f$ with its tuple of coefficients, there exist $r$, $\eta$, and $m$ such that:
  \begin{enumerate}
    \item $(a,f,b)$ is a semi-cgp $(\eta,r)$-correspondence triangle;
    \item $0<r\leq 1$;
    \item $\eta< \min\{\frac{r}{2^{2i_0+4}m}, \bdl(a)2^{-\frac{\bdl(a)}{cr}-7}\}$, where $i_0 := 1+\lfloor\frac{\bdl(a)}{r}\rfloor$ and $c$ is an absolute constant;
    \item $f\in \acl^0((c_i)_{i<m})$ for some $\bdl$-independent sequence $(c_i)_{i<m}$ such that $\dim^0(c_i)=1$ and $\bdl(c_i) \geq r$ for all $i<m$, and $\max\{\bdl(c_i):i<m\}=1$.
  \end{enumerate}

Then there are $h_0,g\in \acl^0(\emptyset)[t]$ and $d \in \acl^0(f)$ such that
$f(x)=h_0(g(x)+_Gd)$, where $+_G$ is either addition or multiplication.
\end{lemma}
\begin{proof}
  The proof is identical to that of Lemma~\ref{l:addmul}, with the following two exceptions:
  \begin{itemize}
    \item We apply Proposition~\ref{p:epsmaincgp} rather than Theorem~\ref{t:main} at the start of the proof;
    \item We apply Lemma~\ref{l:eps1dimab} rather than Corollary~\ref{c:cohActAbSub} near the end, noting that $(a',\phi_f',b')$ remains a semi-cgp $(\eta,r)$-correspondence triangle satisfying the conditions of Lemma~\ref{l:eps1dimab} for $\AGL_1 = \G_a \rtimes \G_m$ acting on $\A^1$.\qedhere
  \end{itemize}
\end{proof}

\section{Balog-Szemerédi-Gowers for group actions}\label{appB}
\label{s:BSGActions}
In this section, we show how our techniques can be applied to obtain a Balog-Szemerédi-Gowers type theorem for arbitrary group actions $(G,X)$. Because of the existence of stabilisers, it is not true any more that from non-expansion (i.e., $\bdl(S*A)=\bdl(A)$ for some $S\subseteq G$ and $A\subseteq X$ of positive $\bdl$)  one can always find approximate subgroups. It is true once we impose some conditions on stabilisers. Similar results were previously obtained in \cite{murphy}.

First, we prove a version for types. As discussed in \S\ref{s:cohHom}, this generalises a part of the proof of Theorem~\ref{t:cohActAb}. As there, we define an increasing sequence of group elements with $\bdl$ which is bounded above and hence tends to a limit. For this argument to go through, we impose a technical condition ((3) below) on the $\bdl$ of the stabilisers.

\begin{lemma}\label{l:typeBSGAction}
Let $M^\L$ be an adequate expansion of $M=(G,X)$.
Suppose there are constants $t\geq 0, n\in\mathbb{N}$, and elements $g\in G$ and $a\in X$, such that, setting $\lambda := \bdl(g)$ and $\mu := \bdl(a)$, we have:
\begin{enumerate}
\item
$0<\lambda,\mu<\infty$;
\item
$a\dind g\dind g*a$;
\item
There exists $\bar{a}=(a_1,\ldots,a_n)$ with $a_1=a$, and $a_i\equiv_ga$, and $a_i\dind_g(a_j)_{j<i}$, such that for any $\bdl$-independent sequence $(g_i)_{i\in I}$ in $\tp(g)$ and any $h\in \langle g_i: i\in I\rangle$,
 we have $\bdl(h/\bar{a},h*\bar{a})\leq t\lambda.$
\end{enumerate}
Then there are $\L'\supseteq\L$ and an adequate expansion $M^{\L'}$,
and a complete $\L'$-type $p$ in the sort $G$ with $\bdl(p)\in[\mu, (n+t)\lambda]$, 
and a complete $\L'$-type $q\supseteq \tp^{\L}(a)$ in the sort $X$ with $\bdl(q)=\lambda$, such that the following hold.
\begin{enumerate}[(I)]
\item
There exists $h\in G$ such that $\bdl((h\cdot p(G))\cap \tp^{\L}(g))\geq\bdl(\tp^\L(g))$.
\item
For all $m$, for all $s_i\vDash p$, for all $\eps_i\in\{-1,1\}$ and for all $a'\vDash q$:
\begin{itemize}
\item
$\bdl(\prod_{i\leq m}s_i^{\eps_i})\leq\bdl(p)$;
\item
$\bdl((\prod_{i\leq m}s_i^{\eps_i})* a')\leq\lambda$.
\end{itemize}
\end{enumerate}
\end{lemma}
\begin{proof}

By Claim~\ref{c:gi} in the proof of Theorem~\ref{t:cohActAb}, we can define $(g_i,g_i')_{i \in \omega}$ and $\delta$-independent sequences $(\bar{h}_i)_{i\in\omega}$ and $(\bar{h}'_i)_{i\in\omega}$ in $\tp(g)$ such that $g_{i+1}=g_i'^{-1}\cdot g_i$, $g_i\in\langle \bar{h}_i\rangle$, $g_i\bar{h}_i\equiv g_i'\bar{h}_i'$, $\bdl(g_i/g_i\cdot g^{-1})=\mu$ and \[\a \dind g_i \dind g_i*\a, \quad g_i' \dind g_i\a.\]

We claim that $\bdl(g_i)\leq (n+t)\lambda$.
Indeed,
\begin{align*}
\bdl(g_i)&=\bdl(g_i,\bar{a})-\bdl(\bar{a}/g_i)=\bdl(g_i,\bar{a},g_i*\bar{a})-\bdl(g_i*\bar{a}/g_i)\\
&=\bdl(\bar{a},g_i*\a)+\bdl(g_i/\bar{a},g_i*\a)-\bdl(g_i*\bar{a})\leq  \bdl(\a)+\bdl(g_i/\bar{a},g_i*\bar{a}).
\end{align*} 
Since $g_i\in\langle\bar{h}_i\rangle$ and $\bar{h}_i$ is a $\bdl$-independent sequence in $\tp(g)$, by assumption, $\bdl(g_i/\bar{a},g_i*\bar{a})\leq t\lambda$.
Thus, $\bdl(g_i)\leq (n+t)\lambda$.
Moreover, $\bdl(g_{i+1})=\bdl((g_i')^{-1}\cdot g_i)\geq\bdl(g_i'/g_i)=\bdl(g_i)$.
In conclusion, the sequence $(\bdl(g_i))_{i\in\omega}$ is non-decreasing and bounded, hence tends to some value $\nu\in[\mu,(n+t)\lambda]$.

Let $(\gamma,\gamma',a',g')$ and an adequate expansion $M^{\L'_0}$ be obtained from Lemma~\ref{l:adequacy}(\ref{adequacy-limit}) applied to $((g_i,g_i',a,g)_{i\in\omega})$. We now work in $\L'_0$.
Then $\tp(a'g')\supseteq \tp^{\L}(ag)$, and 
$\nu=\bdl(\gamma)=\bdl(\gamma')=\bdl(\gamma'^{-1}\cdot\gamma)$, $2\nu=\bdl(\gamma,\gamma')$, $\bdl(a'/\gamma)=\bdl(a')=\lambda=\bdl(\gamma*a')=\bdl(\gamma*a'/\gamma)$ (since $\bdl^\L(g_i*a)=\bdl^\L(a)=\lambda$ and $a{\dind}^\L g_i{\dind}^\L g_i*a$ for all $i$) and $\bdl(\gamma/\gamma\cdot g'^{-1})=\bdl^L(g)=\mu$.
In conclusion, we have 
\[\gamma\dind\gamma'\dind\gamma'^{-1}\cdot\gamma\dind\gamma\text{ and }a'\dind \gamma\dind\gamma*a',\] and $\bdl(\gamma/\gamma\cdot g'^{-1})=\mu$.

Now we apply Corollary~\ref{c:symBSGact} and get a countable locally $\bdl$-finite set $D$ which satisfies the conclusion of the lemma. We may assume $D\dind \gamma\gamma'a'g$, hence $\bdl$ (and $\bdl$-independence) of any subtuples of $\gamma\gamma'a'g$ remains the same over $D$.

Let $\L':=\L'_0(D)$ and $M^{\L'}$ the corresponding expansion, which is adequate by Lemma~\ref{l:adequacy}(\ref{adequacy-constants}).
From now on we work in $\L'$.
Let $\beta\equiv\gamma$ be such that $\beta\dind \gamma g'$.
Let $p:=\tp(\beta^{-1}\cdot\gamma)$ and $q:=\tp(a')$. Note that conclusion (II) now follows from the conclusion of Corollary~\ref{c:symBSGact} and $\beta\dind\gamma$.
For conclusion (I), let $h:=(\beta^{-1}\cdot\gamma\cdot g'^{-1})^{-1}$. Note that \[\bdl(\beta^{-1}\cdot\gamma/\beta^{-1}\cdot\gamma\cdot  g'^{-1})\geq\bdl(\beta^{-1}\cdot\gamma/\beta^{-1}\cdot\gamma\cdot g'^{-1},\beta)=\bdl(\gamma/\gamma\cdot g'^{-1},\beta)=\bdl(\gamma/\gamma\cdot g'^{-1})=\mu\]
where the penultimate equality is because $\gamma g'\dind\beta$, hence $\gamma\dind_{\gamma\cdot g'^{-1}}\beta$.
Suppose $\tp^{\L}(g)\vDash x\in S'$, hence also $\tp(g')\vDash x \in S'$.
Then setting $S'_h:=h^{-1}\cdot S'=\{x\in G: h\cdot x\in S'\}$,
we have $\tp(\beta^{-1}\cdot\gamma/h)\vDash x\in S'_h$.
Hence, $\bdl((h\cdot p(G))\cap S')=\bdl(p(G)\cap S'_h)\geq\bdl(\beta^{-1}\cdot\gamma/h)\geq\mu,$ and $\bdl((h\cdot p(G))\cap \tp^{\L}(g))\geq\mu$ as desired.\qedhere

\end{proof}

For a subset $S\subseteq G$ and $a\in X^n$, we write $\Stab_S(a)$ for $\Stab_G(a)\cap S$.
\begin{theorem}\label{t:BSGAction}
Suppose there are $n \in \mathbb{N}^{>0}$ and $t\in \mathbb{R}^{>0}$ and 
internal sets $S\subseteq G$ with $S=S^{-1}$ and $A,B\subseteq X$, satisfying:
\begin{enumerate}
\item
$0<\bdl(A)=\bdl(B)<\infty$, $0<\bdl(S)<\infty$;
\item
$\bdl(\{(a'_1,\ldots,a'_n)\in A^n:\bdl(\mathrm{Stab}_{S^k}(a'_1,\ldots,a'_n))\geq t\bdl(A)\})<n\bdl(A)$ for all $k$;\footnote{By continuity of $\bdl$, the set $\{(a'_1,\ldots,a'_n)\in A^n:\bdl(\mbox{Stab}_{S^k}(a'_1,\ldots,a_n')\geq t\bdl(A)\}$ is $\bigwedge$-definable over $\emptyset$, hence it makes sense to talk about its $\bdl$-value.}
\item
(Large-incidence) $\bdl(\{(s,a,b)\in S\times A\times B: s*a=b\})=\bdl(S)+\bdl(A)$.
\end{enumerate}
Then there exist an internal coarse approximate subgroup $H\leq G$ and an internal subset $T\subseteq A$ such that the following holds:
\begin{enumerate}[(I)]
\item
$\bdl(T)=\bdl(A)$, and $\bdl(H)\leq (n+t)\bdl(A)$, and there is $h\in G$ with $\bdl(h\cdot H\cap S)=\bdl(S)$ for some $h\in G$;
\item
(Non-expansion) $\bdl(H*T)=\bdl(T)$.
\end{enumerate}
\end{theorem}
\begin{proof}
Work in an adequate expansion of $(G,X,S,A,B)$.
Let $(g,a,g* a)$ be an element of the definable set $\{(g,a,b)\in S\times A\times B: g*a=b\}$ such that $\bdl(g,a,g*a)=\bdl(g,a)=\bdl(g,g*a)=\bdl(S)+\bdl(A)$. Then, $g\dind a$ and $g\dind g*a$, 
and $0<\bdl(g)=\bdl(S)<\infty$ by assumption.
Let $a_1=a$ and for $1\leq i\leq n$, inductively define $a_i\equiv_g a$ with $a_i\dind_g (a_j)_{j<i}$.
Let $\bar{a}:=(a_1,\ldots,a_n)$.
Then $\bdl(\bar{a})=n\bdl(A)$.

To apply Lemma~\ref{l:typeBSGAction}, we need to show that 
\[\bdl(h/\bar{a},h*\bar{a})\leq t\bdl(a)\] for any $h\in\langle \bar{g}\rangle$ where $\bar{g}$ is a $\bdl$-independent sequence in $\tp(g)$. Since each element in $\bar{g}$ is in $S$, there is some $m\in\N$ such that $h\in S^m$.
Let $D:=\{h'\in S^{ m}: h'*\bar{a}=h*\bar{a}\}$, then $\tp(h/\bar{a},h*\bar{a})\vDash D$.
We only need to show that $\bdl(D)\leq t\bdl(a)$.
Note that $h^{-1}\cdot D\subseteq \mathrm{Stab}_{S^{2m}}(\bar{a})$.
If $\bdl(D)> t\bdl(a)$, then $\bdl(\mbox{Stab}_{S^{2m}}(\bar{a}))\geq t\bdl(A)$.
By continuity of $\bdl$, the set \[D_t:=\{(a'_1,\ldots,a'_n)\in A^n:\bdl(\mbox{Stab}_{S^{2m}}(a'_1,\ldots,a_n')\geq t\bdl(A)\}\] is $\bigwedge$-definable over $\emptyset$. Therefore, it contains all realisations of $\tp(a_1,\ldots,a_n)$.
Thus $\bdl(a_1,\ldots,a_n)\leq\bdl(D_t)<n\bdl(A)$, a contradiction.

Let $h\in G$ and complete types $p$ and $q$ in an adequate expansion $\L'$ be as given by Lemma~\ref{l:typeBSGAction}. Then together with Lemma~\ref{l:bdlf(p)}, we have
\begin{itemize}
\item $\bdl(p)\leq (n+t)\bdl(\tp^\L(a))=(n+t)\bdl(A)$, $\bdl(q)=\bdl(\tp^\L(a))=\bdl(A)$;
\item $\bdl((h\cdot p(G))\cap S)\geq\bdl(\tp^\L(g))=\bdl(S)$;
\item $\bdl(p(G)\cdot p(G)\cdot p(G))=\bdl(p(G))$;
\item $\bdl((p(G)\cdot p(G)\cdot p(G))*q(X))\leq\bdl(\tp^\L(a))=\bdl(A)$.
\end{itemize}

Let $E\subseteq p(G)$ and $T\subseteq q(X)$ such that $\bdl(E)=\bdl(p)\leq (n+t)\bdl(A)$ and $\bdl(T)=\bdl(q)=\bdl(A)$ (which exist by Fact~\ref{f:existsInternal}), then
\[\bdl((h\cdot E)\cap S)\geq\bdl(S),
\bdl(E\cdot E\cdot E)=\bdl(E) \text{ and } \bdl((E\cdot E\cdot E)* T)=\bdl(T).\]

Therefore, let $H':=E\cup\{1\}\cup E^{-1}$ and $H:=H'\cdot H'\cdot H'$, then $H$ is a coarse approximate subgroup (by Fact~\ref{f:approx}) such that $\bdl(H)=\bdl(p)$, $\bdl(h\cdot H\cap S)\geq\bdl(S)$ and $\bdl((H* T)\leq\bdl(T)$ (since $E\subseteq H$) as desired.
\end{proof}

\begin{remark}
  In Theorem~\ref{t:BSGAction}, we can not take $H$ to be a translate of a subset of any finite product $S\cdot\ldots\cdot S$, and correspondingly, in Lemma~\ref{l:typeBSGAction}, the limiting process is necessary and we can not take $g$ to be any $g_i$. The following example in the case of $(K,+)$ acting on itself by addition demonstrates this; in the notation of the example, if $g \in S$ with $\bdl(g)=\bdl(S) = 1$, then $\bdl(g_i) = \bdl(\sum_{i<2^i} S)$ so $\bdl(g_i) < \bdl(g_{i+1})$.
\end{remark}
\begin{example}\label{e:span}
  Let $(a_i)_{i\in\omega}$ be a $\Q$-linearly independent sequence of complex numbers.
  For finite $N$, define:
  \begin{itemize}
    \item $n_N := \lfloor\log_2\log_2N\rfloor$.
    \item $I_{N,i} := ([-N^{2^{-i}},N^{2^{-i}}] \cap \Z)\cdot a_i$ for $i<n_N$.
    \item $S_N := \bigcup_{i<n_N} I_{N,i}$.
    \item $A_N := \sum_{i<n_N} I_{N,i}$.
  \end{itemize}
  Let $(S,A) := \prod_{N \to \U} (S_N,A_N)$ be the ultraproduct, and let $\xi := (N)_{N \in \omega} / \U \in \N^\U$.
  Then $\bdl(S) = 1$ and $\bdl(A) = \bdl(S+A) = \bdl(A+A) = 2$,
  but $\bdl(\sum_{i<k} S) = 2 - 2^{-{k-1}}$.
\end{example}

One readily verifies the following finitary version of Theorem~\ref{t:BSGAction}.
\begin{corollary} \label{c:BSGActionFin}
For all $n\in\mathbb{N}^{>0}$ and $r,t>0$, for all $\delta>0$ there is $\eps>0$ such that the following holds.
Let $G$ be a group acting on a set $X$.
Suppose there are finite sets $S\subseteq G$ and $A\subseteq X$ satisfying:
\begin{enumerate}
\item
$\frac{1}{\eps}<|A|^{\frac{1}r}\leq|S|\leq |A|^r$;
\item
$|\{(a_1,\ldots,a_n)\in A^n:|\mathrm{Stab}_{\langle S\rangle}(a_1,\ldots,a_n)|\geq |A|^t\}|\leq |A|^{n-\delta}$, where $\langle S\rangle$ is the group generated by $S$;
\item
(Large-incidence) $|\{(s,a,b)\in S\times A\times A, s*a=b\}|\geq |S||A|^{1-\eps}$.
\end{enumerate}
Then there exist a finite $|A|^{\delta}$-approximate subgroup $H\leq G$ and a subset $T\subseteq A$ such that the following holds:
\begin{enumerate}
\item
$|H|\leq |A|^{n+t+\delta}$, $|H\cap h\cdot S|\geq |A|^{-\delta} |S|$ for some $h\in G$;
\item
$|T|\geq |A|^{1-\delta}$;
\item
(Non-expansion) $|H*T|\leq |A|^{1+\delta}$.
\end{enumerate}
\end{corollary}

\bibliographystyle{alpha}
\bibliography{homogES}
\end{document}